\documentclass[12pt]{amsart}
\usepackage{amsmath,amssymb,amscd,array, mathrsfs }

\usepackage{amsmath,amscd,amssymb,amsthm,array}
\usepackage{color}

\usepackage{amsmath,amssymb,mathrsfs,amsthm, tikz-cd,mathrsfs}
\usepackage{comment}
\def\url#1{\expandafter\s

\tring\csname #1\endcsname}

\newcommand {\foo}{{\mathfrak{oo}}}
\def\mmat #1,#2,#3,#4,{\text{\small\arraycolsep=3pt $
\begin{pmatrix}#1&#2\\#3&#4\end{pmatrix}$}}

\usepackage{hyperref}
\usepackage{capt-of}

\usepackage{multirow}

\usepackage{comments}
\newComments\SBe{Said}{blue}
\newComments\SBo{Sofiane}{blue}
\newComments\AM{Nacer}{blue}
\newComments\DL{DL}{red}
\newComments\QEh{QEh}{blue}

\def\mmat #1,#2,#3,#4,{\text{\small\arraycolsep=3pt $
\begin{pmatrix}#1&#2\\#3&#4\end{pmatrix}$}}

\usepackage{lscape}
\usepackage{tikz-cd}
\usepackage{enumerate}

\usepackage{DLdef1}
\usepackage{multicol}

\def\mmat #1,#2,#3,#4,{\text{\small\arraycolsep=3pt $
\begin{pmatrix}#1&#2\\#3&#4\end{pmatrix}$}}

\renewcommand {\ssbegin}[2][*]
 {\refstepcounter{subsection}%
\if#1*
\addcontentsline{toc}{subsection}{\thesubsection.\hskip 1pc #2}%
\else
\addcontentsline{toc}{subsection}{\thesubsection.\hskip 1pc #2. #1}%
\fi
 \def \secno {\gdef \secno {}{\ssecfont
\thesubsection.\hskip 2ex}%
 }%
 \begin{#2}}

\renewcommand {\sssbegin}[2][*]
 {\refstepcounter{subsubsection}
\if#1*
\addcontentsline{toc}{subsubsection}{\thesubsubsection.\hskip 1pc #2}%
\else
\addcontentsline{toc}{subsubsection}{\thesubsubsection.\hskip 1pc #2. #1}
\fi
 \def \secno {\gdef \secno {}{\ssecfont \thesubsubsection.\hskip 2ex}%
 }%
 \begin{#2}}

\renewcommand {\parbegin}[2][*]
 {\refstepcounter{paragraph}
\if#1*
\addcontentsline{toc}{paragraph}{\theparagraph.\hskip 1pc #2}%
\else
\addcontentsline{toc}{paragraph}{\theparagraph.\hskip 1pc #2. #1}
\fi
 \def \secno {\gdef \secno {}{\ssecfont \theparagraph.\hskip 2ex}%
 }%
 \begin{#2}}

\newcommand {\fba}{{\mathfrak{ba}}}

\newcommand {\ce}{{\text{CE}}}

\newcommand{\trr}{\triangleright}
\rmnameii{etr}{etr}
\rmnameii{evv}{ev}

\DeclareMathOperator{\onabla}{\overline{\nabla}}
\DeclareMathOperator{\K}{\mathbb{K}}

\DeclareMathOperator{\Asso}{\mathrm{Asso}}

\newcommand {\w}{\omega}

\newcommand{\black}{\color{black}}
\newcommand{\Z}{\mathbb{Z}}

\begin{document}

\title[Lagrangian extensions of Lie superalgebras in characteristic 2]{Left-symmetric superalgebras and Lagrangian extensions of Lie superalgebras in characteristic 2}

\author{Sa\"{i}d Benayadi}

\address{Laboratoire de Mathématiques IECL UMR CNRS 7502, Université de Lorraine, 3 rue Augustin Fresnel, BP 45112, F-57073 Metz Cedex 03, France.}
\email{said.benayadi@univ-lorraine.fr}

\author{Sofiane Bouarroudj}

\address{Division of Science and Mathematics, New York University Abu Dhabi, P.O. Box 129188, Abu Dhabi, United Arab Emirates.}
\email{sofiane.bouarroudj@nyu.edu}

\author{Quentin Ehret}
\address {Division of Science and Mathematics, New York University Abu Dhabi, P.O. Box 129188, Abu Dhabi, United Arab Emirates.}
\email{qe209@nyu.edu}

\thanks{S.Bo. and Q.E. were supported by the grant NYUAD-065.}

\keywords {Modular Lie superalgebras; Characteristic 2; Left-symmetric superalgebras; Lagrangian extensions.}

 \subjclass[2020]{17B50; 17B56; 17B05; 17D25}

\begin{abstract}
The purpose of this paper is twofold. First, we introduce the notions of left-symmetric and left alternative structures on  superspaces in characteristic 2.  We describe their main properties and classify them in dimension 2. We show that left-symmetric structures can be queerified if and only if they are left-alternative.

Secondly, we present a method of Lagrangian extension of Lie superalgebras in characteristic 2 with a flat torsion-free connection. We show that any strongly polarized quasi-Frobenius Lie superalgebra can be obtained as a Lagrangian extension. Further, we  demonstrate that Lagrangian extensions are classified by a certain cohomology space that we introduce.  To illustrate our constructions, all Lagrangian extensions in dimension 4 have been described.
\end{abstract}


\maketitle

\thispagestyle{empty}
\setcounter{tocdepth}{2}
\tableofcontents


\section{Introduction} \label{intro}

\subsection{Lie superalgebras in characteristic 2} Over a field of characteristic 2, a Lie superalgebra is a $\Z/2\Z$-graded vector superspace $\fg=\fg_\ev\oplus\fg_\od$ such that the even part $\fg_\ev$ is a Lie algebra, the odd part is a two-sided $\fg_\od$-module due to symmetry, a map $s:\fg_\od\rightarrow\fg_\ev$, known as \textit{squaring} that satisfies some compatibility conditions, including a modified Jacobi identity; see Definition \ref{super2}. In our knowledge, Lebedev (\cite{LeD}) gave the definition, though earlier versions can be found in \cite{BMY}, and also in \cite{NR} under the name of ``strongly anticommutative graded algebras". The presence of the squaring makes it necessary to develop new methods for studying these Lie superalgebras. For example, the concepts of double extensions and Manin triples are developed in \cite{BB1,BB2}, taking into account the squaring. Furthermore, there is a large number of Lie superalgebras that only exist in characteristic 2, with no analogs in other characteristics; see \cite{BGL1, BGLLS, BLLSq, LeD}. As an example, there are four versions of the Hamiltonian Lie superalgebra, which are mutually non-isomorphic.  

Although the classification of simple Lie superalgebras in characteristic 2 remains open, those that admit a Cartan matrix have been classified in \cite{BGL1}. Moreover, it was shown in \cite{BLLSq} that every simple Lie superalgebra can be obtained from a simple Lie algebra using either queerification or ``method 2". The classification problem of simple Lie superalgebras in characteristic 2 is thus reduced to the classification problem of Lie algebras, which is also a challenging one, see \cite{S1,S2}. 

The first instances of a cohomology theory for Lie superalgebras in characteristic 2 can be found in \cite{BGLL}, while its full description can be found in \cite{BM23}. There is only a need for 1- and 2-cocycles in the present work, and their construction is discussed in Section \ref{cohoLie2}.  

\subsection{Lagrangian extensions of Lie (super)algebras} Bordemann introduced the concept of $T^*$-extension  over a field of characteristic 0 (see \cite{Bo}). Basically, starting from a non-associative algebra $A$, he defined an algebra structure on $A\oplus A^*$ by means of an $A^*$-valued bilinear map  satisfying some cohomological conditions, and deduced a non-degenerate invariant {\bf symmetric} bilinear form on $A\oplus A^*$. Most notably, he demonstrated that every finite-dimensional nilpotent algebra over complex numbers equipped with an invariant symmetric bilinear form can be obtained as a $T^*$-extension of a smaller algebra. The term stems from the fact that the Lie algebra of the cotangent bundle of a Lie group is a $T^*$-extension of its Lie algebra.  Later, a superization of Boredmann's approach was given in \cite{BBB}.
 
Baues and Cort\'es examined this concept later (see \cite{BC}), naming it Lagrangian extension. In contrast to Bordemann's approach, it was applied to Lie algebras with flat and torsion-free connections, resulting in Lie algebras admitting a non-degenerate {\bf antisymmetric} closed form. In addition, they introduced Lagrangian cohomology and showed that it captures equivalence classes of Lagrangian extensions. Furthermore, they showed that every quasi-Frobenius Lie algebra admitting a Lagrangian ideal can be obtained through this process.

Baues and Cort\'es studied in depth the structure theory of symplectic Lie groups. Recall that a Lie group $G$ is called symplectic if it is equipped with a left-invariant closed 2-form whose rank is equal to $\dim(G)$. In this case, its Lie algebra $\fg = \mathrm{Lie}(G)$ would be a quasi-Frobenius Lie algebra, that is, there exists a closed non-degenerate 2-form $\omega$ on $\fg$. Another contribution of \cite{BC} is showing that all symplectic Lie groups of dimension less than or equal to 6 admit a Lagrangian subgroup and that there exists an eight-dimensional symplectic Lie group without a Lagrangian subgroup.
 
The concept of Lagrangian extension was later superized in \cite{BM} provided that the field is of characteristic not 2. Lagrangian extensions can be performed in two ways on a flat Lie superalgebra $\fh$. One is to consider either $\fh\oplus\fh^*$ or $\fh\oplus\Pi(\fh^*$), where $\Pi$ is the change in parity functor. They were referred to in \cite{BM} as $T^*$-extension and $\Pi T^*$-extension, respectively, and the anti-symmetric form  they yielded are either even or odd. All Theorems stated in \cite{BC} pertaining to Lagrangian extensions were superized in \cite{BM}; for instance, it was shown that each strongly polarized Lie superalgebra can be obtained as a $T^*$-extension or a $\Pi T^*$-extension of a smaller superalgebra. 

The study of four-dimensional real Lie superalgebras, which are Lagrangian extensions of smaller superalgebras, is presented in \cite{BR}, utilizing Backhouse's classification of real Lie superalgebras (see \cite{Ba}).

All of the constructions in \cite{BM} are valid over any field of characteristic that is not $2$. The purpose of this paper is to explore Lagrangian extensions when the base field has characteristic 2, and to adapt the results of \cite{BM} to that context while taking into account the squaring.

\subsection{Left-symmetric (super)algebras}

In order to perform Lagrangian extensions  in the vein of \cite{BC,BM}, a Lie (super)algebra $\fh$ must be equipped with a flat torsion-free connection. Having a flat torsion-free connection on $\fh$ is equivalent to endowing $\fh$ with a left-symmetric structure, see, e.g., \cite{Bu4}. Recall that a superalgebra $(A,\trr)$ is called left-symmetric if the associator
\begin{equation}\label{as}\text{as}_{\trr}:~ (a,b,c)\mapsto a\trr(b\trr c)-(a\trr b)\trr c\end{equation} is left-symmetric ``in $a$ and $b$''. Namely, $\text{as}_{\trr}(a,b,c)=(-1)^{|a||b|}\text{as}_{\trr}(b,a,c)$ for all homogeneous $a,b,c \in A$.  The first appearance of these structures dates back to 1857, when Cayley  dealt with rooted trees algebras. Vinberg and Koszul reintroduced them in the context of homogeneous cones (\cite{Vi}) and affine flat manifolds (\cite{Ko}). They also appeared under the name ``pre-Lie algebras" in the work of Gerstenhaber on Hochschild cohomology and deformations of algebras. It is common nowadays to use the term ``pre-Lie algebra" even though, in the definition of pre-Lie, the associator is \textit{right}-symmetric. In order to avoid confusion, we use the term ``left-symmetric". For complementary survey papers on left-symmetric algebras, see \cite{Bu4,M}.

Over the real numbers, the problem of determining whether a simply connected solvable Lie group admits a free transitive action on an affine space is equivalent to the algebraic problem of finding a left-symmetric structure on the corresponding solvable Lie algebra, see \cite{Bu2,Me,S}. Simple Lie algebras do not admit left-symmetric structures over the field of complex numbers, see \cite{H}. As it turns out, this is not the case in the super setting. In \cite{DZ}, it was shown that the Lie superalgebra $\mathfrak{sl}(m+1|m)$ admits left-symmetric structures. This paper provides another example involving the Hamiltonian Lie superalgebra $\fh_\Pi^{(1)}(0|4)$ in characteristic 2, see Example \ref{hamilton}.

The classification of left-symmetric algebras was initiated by Segal in \cite{S} and turns out to be a difficult problem. Burde studied simple modular algebras (\cite{Bu1}), Baues studied $\fgl(n)$ (\cite{B}), and Bai classified every 3-dimensional left-symmetric algebra over complex numbers (\cite{Bai}). Dimitrov and Zhang provided a full classification of $\mathfrak{sl}(2|1)$, and proved that $\mathfrak{sl}(m|1),$ for $m\geq 3$, does not admit any left-symmetric structure. In addition, see \cite{KB}, which explored left-symmetric structures on super Virasoro algebras, and \cite{CH}, which dealt with infinite-dimensional superalgebras. 

\subsection{Outline of the paper and main results} In Sections \ref{super2}, \ref{cohoLie2}, and \ref{bill}, we review basic constructions related to Lie superalgebras in characteristic 2 and their cohomology. In Section \ref{LSSA}, we introduce left-symmetric superalgebras in characteristic 2 and give some examples of such structures. Here, the main novelty is the presence of a squaring of the associated Lie superalgebra, as shown in Proposition \ref{LSSA-Lie}. In particular, we provide another example of a simple Lie superalgebra admitting a left-symmetric structure, see Example \ref{hamilton}. In characteristic 2, the queerification method allows to construct Lie superalgebras from Lie algebras and can be applied to any \textit{restricted} Lie algebra, see \cite{BLLSq}. We demonstrate in Section \ref{queer} that a left-symmetric superalgebra in characteristic 2 can be queerified once it is left-alternative, see Lemma \ref{LemQ2}. Moreover, the Lie superalgebra of a queerified left-alternative algebra  is equal to the queerification of the Lie algebra associated with it, see Proposition \ref{lediagrammecommute}. Section \ref{classifLSSA2} is devoted to the classification of left-symmetric superalgebras in dimension 2. We classify all left-symmetric structures compatible with each 2-dimensional Lie superalgebra.

Section \ref{sectionlagrange} develops the theory of Lagrangian extensions of Lie superalgebras in characteristic 2.  First we define connections as non-associative products on a superalgebra with purely algebraic definitions, as well as analogs of torsion and curvature (see Definition \ref{cct}), and their properties. The explicit construction of a Lagrangian extension is given in Theorem \ref{sensdirect}, and its converse is Theorem \ref{conversethm}. Namely, we show that every strongly polarized Lie superalgebra is isomorphic to a $T^*$- or $\Pi T^*$-extension of a suitable smaller superalgebra, that we construct explicitly. Furthermore, we study the functorial properties of our constructions in Section \ref{functor} and demonstrate that the Lagrangian extensions of a given superalgebra are determined by its Lagrangian cohomology, see Theorem \ref{correspondence-thm}. Several examples are given in Section \ref{fewewamples}. To conclude, we compute all 4-dimensional Lie superalgebras which are obtained as $T^*$- or $\Pi T^*$-extensions in Section \ref{sectionclassif4}. 

The appendix \ref{appendixcoho} shows the Lagrangian cohomology groups for each left-symmetric 2-dimensional superalgebra, while the appendix \ref{invariants} shows some invariants for the 4-dimensional Lie superalgebras. In Appendix \ref{postlie}, we initiate the study of post-Lie superalgebras in characteristic 2 and link them with flat parallel connections.\\

\noindent\textbf{Change of parity functor.}  Let $V=V_\ev\oplus V_\od$ a $\Z_2$-graded vector space. We denote by $\Pi$ the \textit{change of parity functor} $\Pi: V\mapsto \Pi (V)$, where $\Pi(V)$ is another copy of $V$ such that $ \Pi(V)_\ev:=V_\od;~~\Pi(V)_\od:=V_\ev$ (see \cite{L}).  Elements of $\Pi(V)$ shall be denoted by $\Pi(v),~\forall v\in V$. We also shall identify $V$ and $\Pi(\Pi(V))$ in a natural way.\\

In this paper, $\K$ stands for an arbitrary field of characteristic 2. We use the terminology ``Lagrangian extension" to mention both $T^*$-extensions and $\Pi T^*$-extensions.\\

\section{Left-symmetric superalgebras in characteristic $2$}

\subsection{Lie superalgebras in characteristic 2} \label{super2}

Let $V$ and $W$ be two vector spaces over $\K$. A map $s:V\rightarrow W$ is called a \emph{squaring} if
\begin{equation*}
    s(\lambda x)=\lambda^2s(x),~\forall \lambda\in\K,~\forall x\in V.\end{equation*}
We also suppose that the map 
\[
V\times V \rightarrow W \qquad (x,y)\mapsto s(x+y)-s(x)-s(y)
\] is bilinear. 

Following \cite{BMY, LeD}, a \emph{Lie superalgebra} over a field $\K$ of characteristic $p=2$ is a $\Z/2\Z$-graded vector space $\fg=\fg_{\ev}\oplus \fg_{\od}$ such that the even part $\fg_{\ev}$ is a Lie algebra, the odd part $\fg_{\od}$ is a $\fg_{\ev}$-module made two-sided by symmetry (the bracket on $\fg_{\ev}$ as well as the action of $\fg_{\ev}$ on $\fg_{\od}$ are given by the same symbol $[\cdot, \cdot]$), a squaring  $s:\fg_{\od}\rightarrow \fg_{\ev}$, such that the bracket of two odd elements is given by:
    \begin{equation}\label{billy}
        [x,y]:=s(x+y)-s(x)-s(y),~~ \forall x,y\in \fg_{\od}.
    \end{equation}
The Jacobi identity involving the squaring reads as follows: 
    \begin{equation}\label{jacobi2} 
        [s(x),y]=[x,[x,y]],~\forall x\in \fg_{\od},~\forall y\in \fg.    \end{equation}
Such a Lie superalgebra will be denoted by $\left(\fg,[\cdot,\cdot]_\fg,s_\fg\right)$ or by just $\fg$ if no confusion is possible.\\

Any associative superalgebra $A=A_{\ev}\oplus A_{\od}$ in characteristic $2$ leads to  a Lie superalgebra. The bracket is the usual commutator and the squaring is defined by $s(x) := x\cdot x,~\forall x \in A_\od$.\\

Let $\fg$ be a Lie superalgebra and $\fh\subseteq  \fg$ be a homogeneous linear subspace. We say that 
    \begin{itemize}
        \item $\fh$ is a \emph{Lie subalgebra} if it is closed under the bracket and under the squaring;
        \item $\fh$ is an \emph{ideal} if $[h,x]\in \fh,~\forall h\in \fh,~\forall x\in \fg$ as well as $s(h)\in \fh, ~\forall h\in \fh_{\od}$. 
    \end{itemize}

Let $\fg$ and $\fh$ be two Lie superalgebras. An even linear map $\varphi: \fg\rightarrow \fh$ is called a \emph{Lie superalgebras morphism} if
\begin{equation}
        \varphi\left([x,y]_\fg\right)=[\varphi(x),\varphi(y)]_\fh,~\forall x\in \fg_{\ev},~\forall y\in \fg, \text{ and }
        \varphi\circ s_\fg(x)=s_\fh\circ\varphi(x),~\forall x\in \fg_{\od}. 
\end{equation}

 A linear map $D: \fg\rightarrow \fg$ is called a \emph{derivation} of $\fg$ if
\begin{equation}\label{derivation}
        D([x,y])=[D(x),y]+[x,D(y)],~\forall x\in \fg_{\ev},~y\in \fg,\text{ and }
        D(s(x))=[D(x),x],~\forall x\in \fg_{\od}.  
\end{equation}
The superspace of all derivations of a Lie superalgebra $\fg$ is itself a Lie superalgebra and is denoted by $\Der(\fg)$.\\

 A \emph{representation} of $\fg$ in a $\Z/2\Z$-graded vector space $M$ is an even map $\rho: \fg\rightarrow \text{End}(M)$ satisfying
\begin{eqnarray*}
\rho([x, y])& = &[\rho(x), \rho(y)] \text{ $\forall x, y \in \fg$; and}\\[1mm]
\rho (s(x)) &=&(\rho(x))^2 \text{ $\forall x \in \fg_{\bar 1}.$}
\end{eqnarray*}
 Such a $M$ is called a \emph{$\fg$-module}.\\~

 In \cite{BGL1}, simple Lie superalgebras admitting a Cartan matrix have been classified. There are several Lie superalgebras in the list that are indigenous to the characteristic 2 and have no analogs in higher characteristic cases.  For examples of vectorial Lie superalgebras, see \cite{BLLSq} and references therein. \\

\noindent\textbf{Notation.} Let $\fg=(\fg,[\cdot,\cdot],s)$ be a Lie superalgebra in characteristic $2$. We denote $$s(\fg_\od):=\text{Span}\{s(x),x\in \fg_\od\}.$$

\sssbegin{Proposition}[\text{\cite[Theorem 3.1]{BM23}}]\label{semidirectLie} Let $(\fg,[\cdot,\cdot],s)$ be a Lie superalgebra and let $(\rho,M)$ be a representation of $\fg$. The superspace $\fg\oplus M$ is a Lie superalgebra with the bracket
\begin{equation}
    [x+v,y+w]_{\ltimes}:=[x,y]+\rho(x)(w)+\rho(y)(v),~\forall x,y\in \fg,~\forall v,w\in M
\end{equation} and the squaring
\begin{equation}
    s_{\ltimes}(x+v):=s(x)+\rho(x)(v),~\forall x\in \fg_\od,~\forall v\in M_\od.
\end{equation}
The Lie superalgebra $(\fg\oplus M,[\cdot,\cdot]_{\ltimes},s_{\ltimes})$ is denoted by $\fg\ltimes M$ and is called the \textit{semidirect product of $\fg$ by the representation $V$.}
\end{Proposition}

\sssbegin{Example}[Affine Lie superalgebra] Let $(\fg,[\cdot,\cdot],s)$ be a Lie superalgebra and let $(\pi,\fgl(\fg))$ be a representation of $\fg$. The semi-direct product $\fg\ltimes\fgl(\fg)$ is a Lie superalgebra called the \textit{affine Lie superalgebra}.   
\end{Example}

\sssbegin{Proposition}\label{classif2}
    Let $\fg$ be a $2$-dimensional Lie superalgebras over an arbitrary field of characteristic $2$. Then,  $\fg$ is isomorphic to one of the following superalgebras.\\
   \begin{itemize} 
     \item \underline{$\sdim(\fg)=(0|2)$}: $\fg={\bf L_{0|2}^1}=\left<0|e_1,e_2\right>$.\\
     
    \item \underline{$\sdim(\fg)=(1|1)$};  $\fg=\left<e_1|e_2\right>$.

    \begin{multicols}{2}
        
        \begin{enumerate}[$(1)$]
            \item ${\bf L_{1|1}^1} =\left<e_1|e_2;[e_1,e_2]=e_2\right>$;

            \item ${\bf L_{1|1}^2}=\left<e_1|e_2; s(e_2)=e_1 \right>$;

        \columnbreak        
            \item ${\bf L_{1|1}^3} =\left<e_1|e_2\right>$ \textup{(}abelian\textup{)};

        \end{enumerate}
    \end{multicols}    
    \item \underline{$\sdim(\fg)=(2|0)$}: $\fg=\left<e_1, e_2|0\right>$.
        \begin{multicols}{2}
        \begin{enumerate}[$(1)$]
            \item ${\bf L_{2|0}^1} =\left<e_1,e_2|0\right>; [e_1,e_2]=e_2$ ;

        \columnbreak             
            \item ${\bf L_{2|0}^2} =\left<e_1,e_2|0\right>$ \textup{(}abelian\textup{)};
           
        \end{enumerate}
        \end{multicols}
    
    \end{itemize} 
\end{Proposition}

\subsection{Cohomology of Lie superalgebras in characteristic 2}\label{cohoLie2}  In this section, we will review the cohomology theory for Lie superalgebras in characteristic $2$. The theory's full description can be found at \cite{BM23}, whereas its first description can be found at \cite{BGLL}. In fact, only 1-cocycles and 2-cocycles are needed. The description of cocycles of higher order can be found at \cite{BM23}.

Let $\fg$ be a Lie superalgebra in characteristic $2$ and let $M$ be a $\fg$-module. 

A 1-cocycle on $\fg$ with values in $M$ is a linear map $\varphi: \fg \rightarrow M$ such that
\begin{equation}
\begin{array}{llll}
d_{\ce}^1(\varphi)(x,z)&:=&x\cdot \varphi(z)+z \cdot \varphi(x)+\varphi([x,z])=0, & \forall x,z\in \fg; \\[2mm] 
\delta^1(\varphi)(x)&:=&x\cdot \varphi(x)+\varphi(s(x))=0,&\forall x\in \fg_\od.
 \end{array}
\end{equation}
The space of 1-cocycles on $\fg$ with values in $M$ is denoted by $XZ^1(\fg;M)$. We also use the notation $\fd^1(\varphi):=(d^1_{\ce}(\varphi),\delta^1(\varphi))$.

A 2-cocycle on $\fg$ with values in $M$ consists of a pair $(\alpha, \gamma)$  such that: 
\[
\begin{array}{l}
\alpha:\fg\wedge \fg \rightarrow M \text{ is a bilinear map}, \text{ and } \gamma:\fg_{\od} \rightarrow M \text{ is a squaring}
\end{array}
\]
such that $\alpha(x,y)=\gamma(x+y)+\gamma(x)+\gamma(y)$ for all $x,y\in \fg_{\od}$. Moreover, we require that 
\begin{equation}\label{2cocycle}
\begin{array}{llll}
d_{\ce}^2(\alpha)(x,y,z)&:=&\underset{x,y,z}{\circlearrowleft}\bigl(x \cdot \alpha(y,z)+\alpha([x,y],z)\bigl)=0, &\forall x,y,z\in \fg; \\[2mm] 
\delta^2(\alpha,\gamma)(x,z)&:=&x\cdot \alpha(x,z)+z \cdot \gamma(x)+\alpha(s(x),z)+\alpha([x,z],x)=0,&\forall x\in \fg_\od,~\forall z\in \fg,
 \end{array}
\end{equation}

The space of 2-cocycles on $\fg$ with values in $M$ is denoted by $XZ^2(\fg;M)$.  We also use the notation $\fd^2(\alpha,\gamma):=(d^1_{\ce}(\alpha),\delta^2(\alpha,\gamma))$.
A complete description of $n$-cocycles for $n\geq 3$ can be found in \cite{BM23}. If $\alpha$ is an even map, then $\text{Im}(\gamma)\subset M_\ev$ whereas if $\alpha$ is odd, then $\text{Im}(\gamma)\subset M_\od$. Therefore, we can define a graduation on the space $XZ^2(\fg;M)$ by setting
$$|(\alpha,\gamma)|:=|\alpha|,~\forall (\alpha,\gamma)\in XZ^2(\fg;M). $$

\subsection{Bilinear forms on Lie superalgebras} \label{bill}

We will briefly review some concepts related to bilinear forms over a superspace in characteristic 2. For more details, see \cite[Sec. 2.2 and 2.4]{BB2}.  Let $\fg=\fg_\ev\oplus \fg_\od$ be a Lie superalgebra in characteristic $2$. A bilinear form $\w$ on $\fg$ with values in $\mathbb K$ is called

\begin{enumerate}[$(i)$]
    \item \emph{non-degenerate} if $\ker(\w):=\{x\in \fg,~\w(x,y)=0~\forall y\in \fg\}=\{0\}$;
    \item \emph{ortho-orthogonal} if $\w$ is even;
    \item \emph{periplectic} if $\w$ is odd;
    \item \emph{closed} if the following cocycle conditions are satisfied:
        \begin{align}
            \underset{x,y,z}{\circlearrowleft}\w([x,y],z)&=0,~\forall x,y,z\in \fg;\\    
            \w(s(x),y)&=\w(x,[x,y]),~\forall x\in \fg_\od,~\forall y\in \fg.\label{cocy2}
        \end{align}
\end{enumerate}

Following \cite{BB2}, an \emph{even} bilinear form on $\fg$ is called \emph{$\od$-antisymmetric} if
\begin{equation}
    \w(x,y)=\w(y,x),~\forall x,y\in \fg \text{ such that } |x|=|y|;~\text{and } \w(x,x)=0~\forall x\in \fg_\ev.    
\end{equation}
An \emph{odd} bilinear form on $\fg$ is called \emph{$\od$-antisymmetric} if
\begin{equation}
    \w(x,y)=\w(y,x),~\forall x,y\in \fg.
\end{equation}
Note that in the case where the form $\w$ is odd, the condition $\w(x,x)=0$ always holds for all $x\in \fg$. For more details, see \cite[Section 2.2]{BB2}.\\

A Lie superalgebra $\fg$ is called {\it quasi-Frobenius} if it is equipped with a $\od$-antisymmetric non-degenerate closed form $\omega$. We denote such an algebra by $(\fg, \omega)$. A quasi-Frobenius Lie superalgebra $(\fg, \omega )$ is called {\it ortho-orthogonal quasi-Frobenius} (resp.  {\it  periplectic quasi-Frobenius}) if the form $\omega$ is even (resp. odd) on $\fg$. \\

Suppose that the form $\w$ is homogeneous, non-degenerate and closed. Let $I\subset \fg$ be a homogeneous ideal and let $I^{\perp}$ be its orthogonal with respect to $\w$ , which is homogeneous as well. Moreover, suppose  that $[I,I^{\perp}]=0$. Then the space $I^{\perp}$ is an ideal. Let us just check that it is closed under squaring. Indeed, let $u\in I^{\perp}_\od$. Then, using \eqref{cocy2}, we have
$$ \w\bigl(s(u),v\bigl)=\w\bigl(u,[u,v]\bigl)=0~~\forall v\in I.$$ Therefore, the quotient $\fg/I^{\perp}$ is well defined.\\
An ideal $I\subset L$ is called \textit{Lagrangian} if $I=I^\perp$. A homogeneous Lagrangian ideal of a quasi-Frobenius Lie superalgebra is always abelian.

\subsection{Left-symmetric superalgebras in characteristic 2}\label{prelie}
\label{LSSA} A left-symmetric superalgebra $(V,\trr)$ in characteristic $p=2$ is a vector superspace $V=V_\ev\oplus V_\od$ endowed with a bilinear product $\trr:V\times V\rightarrow V$ satisfying
\begin{equation*}
\begin{array}{lrlll}
(i) &x\trr (y\trr z)+(x\trr y)\trr z& = & y\trr (x\trr z)+(y\trr x)\trr z,&\forall x,y,z\in V; \\[2mm]
(ii) &x\trr(x\trr y)&=&(x\trr x)\trr y,&\forall x\in V_\od,~\forall y\in V.\\[2mm]
\end{array}
\end{equation*}
Any associative superalgebra is a left-symmetric superalgebra.\\
A linear map $\phi:(V,\trr)\rightarrow(W,\blacktriangleright)$ is called \textit{morphism of left-symmetric superalgebras} if $$\phi(x\trr y)=\phi(x)\blacktriangleright\phi(x),~\forall x,y\in V.$$

\sssbegin{Proposition}\label{LSSA-Lie}
    Let $(V,\trr)$ be a left-symmetric superalgebra. Then, $(\fg(V),[\cdot,\cdot],s)$ is a Lie superalgebra with $\fg(V)=V$ as superspaces and
    \begin{align}\label{LSSA-Lie1}
        [x,y]&:=x\trr y+y\trr x,~\forall x\in V_\ev,\forall y\in V;\\\label{LSSA-Lie2}
        s(x)&:=x\trr x,~\forall x\in V_\od.
    \end{align}
    A left-symmetric product $\trr$ on a Lie superalgebra $(V,[\cdot,\cdot],s)$ is called compatible with the Lie superalgebra structure if Conditions \eqref{LSSA-Lie1} and \eqref{LSSA-Lie2} are satisfied.
\end{Proposition}
\begin{proof}
   Obviously, $s(\lambda x)=\lambda^2 s(x)$ for any $x\in V_\od$ and any $\lambda\in \K$, and the map 
   $$(x,y)\mapsto s(x+y)-s(x)-s(y)=x\trr y+y\trr x$$ is bilinear for all $x,y\in V_\od.$
   Moreover, let $x\in V_\od$ and $y\in V$. We have
    \begin{equation*}
    \begin{array}{llll}
    [s(x),y]&=&(x\trr x)\trr y+y\trr(x\trr x)&\\
    &=&(x\trr y)\trr x+x\trr(x\trr y)+(y\trr x)\trr x+x\trr(y\trr x)&\\
    &=&x\trr(x\trr y+y\trr x)+(x\trr y+y\trr x)\trr x&\\
    &=&[x,x\trr y+y\trr x]=[x,[x,y]].&\qed
    \end{array}
    \end{equation*}\noqed
\end{proof}

 \sssbegin{Proposition}[\text{\cite[Theorem 2.5.3]{BB2}}] Let $(\fg,[\cdot,\cdot],s)$ be a Lie superalgebra equipped with an invertible derivation $D$. Let $x\trr y:=D^{-1}\bigl([x,D(y)]\bigl),~\forall x,y\in \fg$. Then, $\trr$ is a left-symmetric product compatible with the Lie structure in the sense of Proposition \ref{LSSA-Lie}. 
 \end{Proposition}

\sssbegin{Example} Consider the five-dimensional Lie superalgebra $\fg=\foo_{I\Pi}^{(1)}(1|2)$ spanned by $\{s(x),s(y),h|x,y\}$ with brackets $[h,x]=x,~[h,y]=y$ and $[x,y]=h$. Since $XH^1(\fg,\fg)=0$, all the derivations are inner (see \cite[Lemma 2.7]{BGLL}). Consider the inner derivation
$$D_{\lambda}:=\lambda_1\ad_h+\lambda_2\ad_{s(x)}+\lambda_3\ad_{s(y)}+\lambda_4\ad_{x}+\lambda_5\ad_{y},\quad \lambda_i\in\K.$$
This derivation is not invertible for any choice of the parameters $\lambda_i\in\K$. Therefore, the superalgebra $\foo_{I\Pi}^{(1)}(1|2)$ does not admit any invertible derivation.
\end{Example}

 \sssbegin{Example} Consider the Lie superalgebra $\fg=\foo^{(1)}_{II}(1|2)$ of superdimension $\sdim(\fg)=(3|2)$ of \textit{symmetric} $(1|2)\times(1|2)$ matrices with supertrace zero. We have $XH^1(\fg,\fg)=0$, all the derivations are inner (see \cite[Lemma 2.7]{BGLL}), and any linear combination of them is non-invertible.
 \end{Example}

\sssbegin{Example}[The Hamiltonian superalgebra $\fh_{\Pi}(0|4)$]\label{hamilton} Consider the Hamiltonian superalgebra $\fh_{\Pi}(0|4)$, see \cite[5.4]{BB1} and \cite{LeD}. As a vector space it can be considered as follows, where $\xi_1,\xi_2,\eta_1,\eta_2$ are odd indeterminates:
$$
\fh_{\Pi}(0|4)\simeq\Span\{H_f\; |\; f\in \Kee[\xi,\eta]\}\simeq \Kee[\xi,\eta]/\Kee\cdot 1,$$
where  $$
H_f= \frac{\partial f}{\partial {\xi_1}} \frac{\partial  }{\partial {\eta_1}} +\frac{\partial f}{\partial {\eta_1}} \frac{\partial  }{\partial {\xi_1}}+\frac{\partial f}{\partial {\xi_2}} \frac{\partial  }{\partial {\eta_2}} +\frac{\partial f}{\partial {\eta_2}} \frac{\partial  }{\partial {\xi_2}}.$$
The Lie bracket $[H_f, H_g]=H_{\{f,g\}}$ is given by the Poisson bracket:
$$
\{f,g\}:= \frac{\partial f}{\partial {\xi_1}} \frac{\partial  g}{\partial {\eta_1}} +\frac{\partial f}{\partial {\eta_1}} \frac{\partial g }{\partial {\xi_1}}+\frac{\partial f}{\partial {\xi_2}} \frac{\partial  g}{\partial {\eta_2}} +\frac{\partial f}{\partial {\eta_2}} \frac{\partial g }{\partial {\xi_2}},
$$
and the squaring is zero.
Let $D:=D_2+D_3+D_5+D_7+\ad_{\eta_1+\eta_2}$, where

$$\begin{array}{ll}
&D_2=\xi _1\otimes \eta _1^*+( \xi _1\, \xi
   _2)\otimes  (\xi _2\, \eta _1)^*+  (\xi _1\,
   \eta _2) \otimes (\eta _1\, \eta _2 )^*+  (\xi
   _1\, \xi _2\, \eta _2) \otimes \left(\xi _2\, \eta _1\, \eta _2\right)^*;\\[2mm]
   
    &D_3= \xi _2 \otimes  \eta _2^* +  (\xi _1\, \xi
   _2) \otimes (\xi _1\, \eta _2)^* + ( \xi _2\,
   \eta _1 )\otimes (\eta _1\, \eta _2)^*+ ( \xi
   _1\, \xi _2\, \eta _1 )\otimes \left (\xi _1\, \eta _1\, \eta _2\right)^*;\\[2mm]
   
   &D_5=\eta _2 \otimes\xi _2^*+ (\xi _1\, \eta
   _2)  \otimes(\xi _1\, \xi _2)^*+ (\eta _1\,
   \eta _2) \otimes (\xi _2\, \eta _1)^*+ (\xi
   _1\, \eta _1\, \eta _2)  \otimes\left(\xi _1\, \xi _2\, \eta _1\right)^*;\\[2mm]
   
&D_7= (\xi _1\, \xi _2\, \eta _1)\otimes  \xi
   _2^*+(\xi _1\, \xi _2\, \eta _2) \otimes \xi
   _1^*+(\xi _1\, \eta _1\, \eta _2)  \otimes\eta
   _2^*+(\xi _2\, \eta _1\, \eta _2)\otimes  \eta _1^*.
\end{array}$$
Then, $D$ is a derivation of the derived Lie superalgebra $\fh^{(1)}_{\Pi}(0|4)$ (see \cite[Claim 5.4.1]{BB1}). Moreover, the derivation $D$ is invertible. Therefore, it defines a left-symmetric structure on $\fh^{(1)}_{\Pi}(0|4)$. Jacobson proved in \cite{J2} that Lie algebras with invertible derivations must be nilpotent in characteristic zero. In this example, we obtain an invertible derivation of a simple Lie superalgebra in characteristic two.
\end{Example}

\noindent{\bf Open problems:} 
\begin{enumerate}[$(i)$]
\item Investigate whether the  
Hamiltonian type Lie superalgebras  $\fh^{(1)}_{II}(m;\underline{N}|n),~\fh^{(1)}_{I\Pi}(m;\underline{N}|n),$ $\fh^{(1)}_{\Pi I}(m;\underline{N}|n),\fh^{(1)}_{\Pi\Pi}(m ;\underline{N}|n)$ and $\fle^{(1)}(m ;\underline{N}|m)$ (see \cite{ LeD}), admit an invertible derivation or not. 

\item Determine whether they admit left-symmetric structures.
\end{enumerate}

\sssbegin{Proposition} $($see \cite[Section 2.5]{BB2}$)$
    Let $(\fg,[\cdot,\cdot],s)$ be a Lie superalgebra equipped with a closed non-degenerate (ortho-orthogonal or periplectic) bilinear form $\w$. Define a product $\trr: \fg\times \fg\rightarrow \fg$ by
    \begin{equation}
        \w(x\trr y,z):=\w(y,[x,z]),~\forall x,y,z\in \fg.
    \end{equation}
    Then, $\trr$ is a left-symmetric structure on $\fg$ that is compatible with the Lie structure (in the sense of Proposition \ref{LSSA-Lie}).
\end{Proposition}


\sssbegin{Proposition}[See also \text{\cite{Bu4}}]\label{invertiblecocycle} 
     Let $(\fg,[\cdot,\cdot],s)$ be a Lie superalgebra. There exists a left-symmetric product $\trr$ on $\fg$ compatible with the Lie structure (in the sense of Proposition \ref{LSSA-Lie}) if and only if there exists a $\fg$-module $M$ satisfying $\sdim(M)=\sdim(\fg)$ such that $XZ^1(\fg;M)$ contains an invertible element.
\end{Proposition}

\begin{proof}
    Let $M$ be a $\fg$-module $M$ satisfying $\sdim(M)=\sdim(\fg)$ and let $\varphi\in XZ^1(\fg;M)$ be an invertible 1-cocycle. We denote the action of $\fg$ on $M$ by $\rho:\fg\rightarrow\End(M)$. We consider the map $f:\fg\rightarrow\End(M)$ defined by
    \begin{equation}
        f(x):=\varphi^{-1}\circ\rho(x)\circ\varphi,~\forall x\in \fg.
    \end{equation}
    \textbf{Claim.} The map $f$ defines a representation of the Lie superalgebra $\fg$.\\
   We prove the identity involving the squaring. Let $x\in \fg_\od$. We have
    $$f\bigl(s(x)\bigl)=\varphi^{-1}\circ \rho(s(x))\circ\varphi=\varphi^{-1}\circ \rho(x)^2\circ\varphi=\varphi^{-1}\circ\rho(x)\circ\varphi\circ \varphi^{-1}\circ\rho(x)\circ\varphi=f(x)^2. $$
    Let us define a product on $\fg$ by $x\trr y:=f(x)(y),~\forall x,y\in \fg$.
    Let $x\in \fg_\od$. Recall that since $\varphi$ is a 1-cocycle, we have $\rho(x)\circ \varphi(x)=\varphi\bigl(s(x)\bigl)$. Therefore,
    $$x\trr x=\varphi^{-1}\circ\rho(x)\circ\varphi(x)=\varphi^{-1}\circ\varphi\circ s(x)=s(x). $$ Similarly, we obtain that $[y,z]=y\trr z+z\trr y,~\forall y,z\in \fg$. Moreover, for any $y\in \fg_\od$, we have
    \begin{align*}
        x\trr(x\trr y)&=f(x)\bigl(\varphi^{-1}\circ\rho(x)\circ\varphi(y)\bigl)\\
        &=\varphi^{-1}\circ\rho(x)^2\circ\varphi(y)\\
        &=\varphi^{-1}\circ\rho(s(x))\circ\varphi(y)\\
        &=\varphi^{-1}\circ\rho(x\trr x)\circ\varphi(y)=(x\trr x)\trr y.
    \end{align*}
    proving identity $(i)$ of Definition \ref{LSSA} is routine. Reciprocally, suppose that $\trr$ is a left-symmetric structure on $(\fg,[\cdot,\cdot],s)$ satisfying \eqref{LSSA-Lie1} and \eqref{LSSA-Lie2}. Then, take $M=\fg$, the $\fg$-module structure being given by $\fg\rightarrow \End(\fg),~x\mapsto (x\trr\cdot)$. One can also check that $\id$ is an invertible element of $XZ^1(\fg;M)$.
\end{proof}

\sssbegin{Corollary}
[See also \text{\cite{Bu4}}]    Let $(\fg,[\cdot,\cdot],s)$ be a Lie superalgebra of total dimension $\dim(\fg)=n$. Suppose that there exists a left-symmetric product $\trr$ on $\fg$ compatible with the Lie structure (in the sense of Proposition \ref{LSSA-Lie}). Then, $(\fg,[\cdot,\cdot],s)$ admits a faithful representation of dimension $n+1$.
\end{Corollary}

\begin{proof}
    Using Proposition \ref{invertiblecocycle}, there exists a $\fg$-module $M$ of dimension $n$ such that $XZ^1(\fg;M)$ contains an invertible cocycle denoted by $\varphi$. We also denote the action of $\fg$ on $M$ by $\rho:\fg\rightarrow\End(M)$. We consider
    $$\pi:\fg\rightarrow \fgl(M\times \K),~ \pi(x)(m,\lambda):=\bigl(\rho(x)(m)+\lambda\varphi(x),0\bigl),~\forall x\in \fg,\forall m\in M,\forall \lambda\in \K.$$
    Let $x\in \fg_\od$ and $m\in M$.
    \begin{align*}
        \pi\circ s(x)(m,\lambda)&=\bigl(\rho(s(x))(m)+\lambda\varphi(s(x)),0\bigl)\\
         &=\bigl(\rho(x)^2(m)+\lambda\rho(x)(\varphi(x)),0\bigl)\\
         &=\bigl(\rho(x)\bigl(\rho(x)(m)+\lambda\varphi(x)\bigl),0\bigl)=\pi(x)^2(m,\lambda).
    \end{align*} A similar computation shows that $\pi([x,y])=\pi(x)\circ\pi(y)+\pi(y)\circ\pi(x),~\forall x,y\in \fg$. Therefore, $\pi$ is a representation.
    Let $x\in\Ker(\pi)$. We have
    $$(0,0)=\pi(x)(0,1)=\bigl(\varphi(x),0\bigl). $$ Therefore, $x\in \Ker(\varphi)$ and it follows that $x=0$ since $\varphi$ is invertible. Thus, the representation $\pi$ is faithful.
\end{proof}

\sssbegin{Proposition}[See also \text{\cite{Bu4}}]\label{propLS2} 
   Let $(\fg,[\cdot,\cdot],s)$ be a Lie superalgebra.  There exists a left-symmetric product $\trr$ on $\fg$ compatible with the Lie structure (in the sense of Proposition \ref{LSSA-Lie}) if and only if there exists an isomorphism of vector superspaces $\varphi:\fg\rightarrow \fg$ and a linear map $\pi:\fg\rightarrow\fgl(\fg)$ such that the map
   $$\Psi: \fg\rightarrow\fgl(\fg)\ltimes \fg,~x\mapsto \bigl(\pi(x),\varphi(x)\bigl)$$ is a morphism of Lie superalgebras. In that case, $\pi$ is a representation, $\varphi\in XZ^1(\fg;\fg)$, where $\fg$ is seen as an $\fg$-module using $\pi$, and $x\trr y=(\varphi^{-1}\circ\pi(x)\circ\varphi)(y),~\forall x,y\in \fg.$
\end{Proposition}

\begin{proof}
    Recall that the Lie superalgebra structure on the semidirect product $\fgl(\fg)\ltimes \fg$ is given in Proposition \ref{semidirectLie}. Suppose that there exists $\varphi,\pi$ and $\Psi$ as in Proposition \ref{propLS2}. We will only show the identities involving the squaring. Since $\Psi$ is a morphism of Lie superalgebras, we have
    $$\bigl(\pi\circ s(x),\varphi\circ s(x)\bigl)=\Psi(s(x))=s_{\ltimes}\bigl(\Psi(x)\bigl)=\bigl(\pi(x)^2,\pi(x)\circ \varphi(x)\bigl),~\forall x\in \fg_\od.$$ It follows that $\pi\circ s(x)=\pi(x)^2,~\forall x\in \fg_\od$ and $\varphi\circ s(x)=\pi(x)\circ \varphi(x)$. Thus, the map $\pi$ is a representation and $\varphi$ is a 1-cocycle. We can then apply Proposition \ref{invertiblecocycle} to construct the left-symmetric product $\trr$. 
    
    Reciprocally, if $\fg$ is equipped with a left-symmetric product $\trr$, we define $$\pi(x)(y):=x\trr y~\forall x,y\in \fg \text{ and } \varphi:=\id.$$
\end{proof}

\subsubsection{Analogs of $\mathcal{O}$-operators in characteristic $2$} Let $(\fg,[\cdot,\cdot],s)$ be a Lie superalgebra and let $M$ be a $\fg$-module. An $\mathcal{O}$\textit{-operator} on $\fg$ with respect to $M$ is an even map $R:M\rightarrow \fg$ satisfying
    \begin{align}
        [R(m),R(n)]&=R\bigl(R(m)\cdot n+m\cdot R(n)\bigl),~\forall m,n\in M;\\
        s\circ R(m)&=R\bigl(R(m)\cdot m),~\forall m\in M_\od.
    \end{align}
    In the particular case where $M=\fg$ and the action is given by the adjoint representation, the above identities become
    \begin{align}
        [R(x),R(y)]&=R\bigl([R(x),y]+[x, R(y)]\bigl),~\forall x,y\in \fg;\\
        s\circ R(x)&=R\bigl([R(x),x]),~\forall x\in \fg_\od.
    \end{align}
    In that case, the $\mathcal{O}$-operator $R$ is called \textit{Rota-Baxter operator}.\\
    
\noindent\textbf{Example.} Let $(\fg,[\cdot,\cdot],s)$ be a Lie superalgebra and let $D$ be an even invertible derivation of $\fg$. Then, $R=D^{-1}$ is a Rota-Baxter operator on $\fg$.

\sssbegin{Proposition}[See also \text{\cite{Bai2}}]
    Let $(\fg,[\cdot,\cdot],s)$ be a Lie superalgebra equipped with a Rota-Baxter operator $R:\fg\rightarrow \fg$. Then, $(\fg,\trr)$ is a left-symmetric superalgebra, with $x\trr y:=[R(x),y],~\forall x,y\in \fg$.
\end{Proposition}

\begin{proof}
    We prove the identity involving the squaring. Let $x\in \fg_\od$ and $y\in \fg$. 
    $$x\trr(x\trr y)=[R(x),[R(x),y]]=[s\circ R(x),y]=\bigl[R\bigl([R(x),x]\bigl),y  \bigl]=(x\trr x)\trr y.   $$
\end{proof}

\subsection{Queerification}\label{queer} This section aims at adapting the queerification process (see \cite{BLLSq}) to left-symmetric algebras. Recall that an algebra $(V,\trr)$ is called \textit{left-alternative} if it satisfies the identity
\begin{equation}
    x\trr(x\trr y)=(x\trr x)\trr y,~\forall x,y\in V.
\end{equation}
It has been shown in \cite[Proposition 2.5.1]{BB2} that every left-alternative algebra is left-symmetric in characteristic 2. The converse does not hold, see Section \ref{classifLSSA2} for several counterexamples in dimension 2 and Remark \ref{pasha} for a counterexample in dimension 10.

\subsubsection{Restricted Lie algebras in characteristic 2} \label{def:restrictedLie} Let $\fg$ be  a finite-dimensional Lie algebra over $\K$. Following \cite{J}, a map $(\cdot)^{[2]}:\fg\rightarrow \fg, ~x\mapsto x^{[2]}$ is called a~\textit{$2$-structure} on $\fg$ and $\fg$ is said to be {\it restricted} if 
\begin{equation}\label{RRR}
\begin{array}{lrll}
(i)&\text{ $\ad_{x^{[2]}}$}&=&\text{$(\ad_x)^2$ for all $x\in \fg$;}\\[2pt]

(ii)&\text{ $(\lambda x)^{[2]}$}&=&\text{$\lambda^2 x^{[2]}$ for all $x\in \fg$ and for all $\lambda \in \mathbb{K}$;}\\[2pt]

(iii)&\text{ $(x+y)^{[2]}$}&=&\text{$x^{[2]}+y^{[2]}+[x,y]$, for all $x,y\in \fg.$}\\
\end{array}
\end{equation}
It has be shown in \cite[3.2]{BLLSq} that every restricted Lie algebra in characteristic 2 can be queerified. Indeed, let $\bigl(\fg,[\cdot,\cdot](\cdot)^{[2]}\bigl)$ be a restricted Lie algebra. Consider the superspace $\fq(\fg):=\fq(\fg)_\ev\oplus\fq(\fg)_\od$, where $\fq(\fg)_\ev:=\fg$ and $\fq(\fg)_\od:=\Pi(\fg)$. Then, $\bigl(\fq(\fg),[\cdot,\cdot]_\fq,s_\fq\bigl)$ is a Lie superalgebra with
    \begin{equation}\label{queerLie}
        [x,y]_\fq:=[x,y],~[x,\Pi(y)]_\fq:=\Pi([x,y]),~s_\fq\bigl(\Pi(x)\bigl):=x^{[2]},~\forall x,y\in\fg.
    \end{equation}
Note that
\sssbegin{Lemma}\label{LemQ1}
    An algebra $(V,\trr)$ is left-alternative if and only if $(\fg(V),[\cdot,\cdot],(\cdot)^{[2]})$ is a restricted Lie algebra, where $\fg(V)=V$ as vector spaces, $[x,y]:=x\trr y-y\trr x,~\forall x,y\in V$ and $x^{[2]}:=x\trr x,~\forall x\in V.$
\end{Lemma}

\begin{proof}
Let $x\in V$. We show that the map $x\mapsto x\trr x$ is a 2-structure on $\fg(V)$. Conditions $(ii)$ and $(iii)$ of Definition \ref{def:restrictedLie} are readily satisfied. Moreover, for any $y\in V$, we have
\begin{align*}
    [x,[x,y]]-[x\trr x,y]&=x\trr(x\trr y)+x\trr(y\trr x)+(x\trr y)\trr x+(y\trr x)\trr x\\&\quad-\bigl((x\trr x)\trr y+y\trr(x\trr x)\bigl)\\
    &=x\trr(x\trr y)-(x\trr x)\trr y.
\end{align*}
Therefore, Condition $(i)$ of Definition \ref{def:restrictedLie} is satisfied if and only if $(V,\trr)$ is left-alternative. 

\end{proof}

\subsubsection{Queerification of left-symmetric algebras} Let $(V,\trr)$ be a left-symmetric algebra and consider the superspace $\fq(V):=\fq(V)_\ev\oplus\fq(V)_\od$, where $\fq(V)_\ev:=V$ and $\fq(V)_\od:=\Pi(V)$. Define a product on $\fq(V)$ by
\begin{equation}\label{queerLSA}
    x\trr_\fq y:=x\trr y;~x\trr_\fq \Pi(y):=\Pi(x)\trr_\fq y:=\Pi(x\trr y),~\Pi(x)\trr_\fq \Pi(y):=x\trr y,~\forall x,y\in V.
\end{equation}

\sssbegin{Lemma}\label{LemQ2}
    The superspace $\fq(V)$ endowed with the product $\trr_\fq$ is a left-symmetric superalgebra if and only if $(V,\trr)$ is a left-alternative algebra.
\end{Lemma}
\begin{proof}
Suppose that $(\fq(V),\trr_\fq)$ is a left-symmetric superalgebra. Let $x,y\in V$. We have
$$\Pi(x)\trr_\fq\bigl(\Pi(x)\trr_\fq y\bigl)=\bigl(\Pi(x)\trr_\fq\Pi(x)\bigl)\trr_\fq y. $$ By definition of the product $\trr_\fq$, this implies that 
$x\trr(x\trr y)=(x\trr x)\trr y,~\forall x,y\in V.$ Thus, $(V,\trr)$ is left-alternative. Reciprocally, given a left-alternative algebra $(V,\trr)$, it is routine to show that its queerification $(\fq(V),\trr_\fq)$ is left-symmetric.
\end{proof}

\sssbegin{Proposition}\label{lediagrammecommute} Let $(V,\trr)$ be a left-alternative algebra. Then, $(V,\trr)$ can be queerified into a left-symmetric superalgebra and the following diagram commutes.
\begin{center}
\begin{tikzcd}
\text{left-alternative algebra} \arrow[rr, "\fq"] \arrow[d, "\fg(\cdot)"'] &  & \text{left-symmetric superalgebra} \arrow[d, "\fg(\cdot)"] \\
\text{restricted Lie algebra} \arrow[rr, "\fq"']          &  & \text{Lie superalgebra}          
\end{tikzcd}
\end{center}
\end{Proposition}
    \begin{proof}
        Lemmas \ref{LemQ1} and \ref{LemQ2} together with \eqref{queerLie} and \eqref{queerLSA} imply that $\fg(\fq(V)=\fq(\fg(V))$ for any left-alternative algebra $(V,\trr)$.
    \end{proof} 

\sssbegin{Remark}\label{pasha} Here we provide an answer to an open problem posed in \cite[5.(7)(i)]{BB2}. P. Zusmanovich provided us with an example of a 10-dimensional left-symmetric algebra that is not left-alternative. This algebra is given in the basis $\{e_1,e_2,\cdots,e_{10}\}$ by the non-zero products
\begin{align*}
e_1e_1&=e_5;~e_1e_2=e_4;~e_1e_3=e_7+e_8+e_9;~e_1e_4=e_6;\\
~e_2e_1&=e_3;~e_2e_5=e_7;~e_3e_1=e_8;~e_4e_1=e_9;~e_5e_2=e_{10}.\\
\end{align*}
We also provide other examples in Section \ref{classifLSSA2}.
\end{Remark}

\subsection{Classification of left-symmetric superalgebras in dimension 2}\label{classifLSSA2} For each Lie superalgebra $\fg$ of Proposition \ref{classif2}, we classify up to isomorphism all the non-zero left-symmetric structures that are compatible with the bracket and the squaring of $\fg$ in the sense of Proposition \ref{LSSA-Lie}. We also indicate whether those algebras are left-alternative or not. In the case of dimension 2, left-alternativity is equivalent to associativity. The computations were performed with the SuperLie package of the computer algebra system Mathematica (see \cite{G}).

\subsubsection{The case where $\sdim(\fg)=(1|1)$} See Proposition \ref{classif2} for the definition of the superalgebras.
\begin{center}
    \begin{tabular}{|c|c|c|c|}
    \hline
          $\fg$ & Left-symmetric product on $\fg$& Conditions    &Left-alternative\\\hline
             \multirow{4}{*}{$\bf{L^1_{1|1}}$}& $e_1e_2=e_2$& None & No \\\cline{2-4} 
            & $e_1e_1=\eps e_1;~e_1e_2=e_2$& $\eps\neq 0,1$& No \\\cline{2-4} 
            & $e_1e_1= e_1;~e_1e_2=e_2$&None& Yes \\\cline{2-4} 
          & $e_1e_1=\eps e_1;~e_1e_2=(1+\eps)e_2;~e_2e_1=\eps e_2$ & $\eps\neq0$& Yes iff $\varepsilon=1$ \\\hline

         \multirow{2}{*}{$\bf{L^2_{1|1}}$} & $e_2e_2=e_1$&None   & Yes\\ \cline{2-4} 

        & $e_1e_1=e_1;~e_1e_2=e_2;~e_2e_1=e_2;~e_2e_2=e_1$ &None  & Yes\\ \hline 
    
        \multirow{2}{*}{$\bf{L^3_{1|1}}$} & $e_1e_1=e_1$&None & Yes \\\cline{2-4} 
        & $ e_1e_1=e_1;~e_1e_2=e_2;e_2e_1=e_2$& None& Yes \\\hline

    \end{tabular}
   \\~\\ The case where $\sdim(\fg)=(1|1)$.
\end{center}

\subsubsection{The case where $\fg=\bf{L^1_{2|0}}$} See Proposition \ref{classif2} for the definition of the superalgebra.

 \begin{center}
     \begin{tabular}{|c|c|c|}
     \hline
          Left-symmetric product on $\fg$ & Conditions & Left-alternative \\\hline
           
            $e_1e_1=\eps e_1;~e_1e_2=e_2$ & $\eps \neq 0, 1$& No \\\hline
            $e_1e_1=e_1+\eps e_2;~e_1e_2=e_2$ & $\eps \neq 0$& No \\\hline
            $e_1e_1=\eps e_2;~e_1e_2=e_2$ & None & No\\\hline
            $e_1e_1=e_1;~e_1e_2=e_2$ & None & Yes \\\hline
            $e_1e_1=(1+\eps)e_1;~e_1e_2=\eps e_2;~e_2e_1=(1+\eps)e_2;~e_2e_2=e_1$ & None& No \\\hline
            $e_1e_1=(1+\eps)e_1;~e_1e_2=\eps e_2;~e_2e_1=(1+\eps)e_2;~e_2e_2=0$ & None & Yes iff $\eps=0$ \\\hline
            $e_1e_1=e_1;~e_1e_2=e_1;~e_2e_1=e_1+e_2;~e_2e_2=e_2$ & None & No\\\hline
           $e_1e_1=e_1;~e_1e_2=e_1;~e_2e_1=e_1+e_2;~e_2e_2=\eps e_1+e_2$ & $\eps \neq 0$ & No\\\hline
            $e_1e_1=e_1+e_2;~e_1e_2=e_1;~e_2e_1=e_1+e_2;~e_2e_2=e_2$ & None & No\\\hline
            
            
            $e_1e_1=(1+\eps)e_1+\eps(1+\eps)e_2;~e_1e_2=e_1+\eps e_2;$&\multirow{2}{*}{$\eps \neq 0, 1$}&\multirow{2}{*}{No}\\
            $e_2e_1=e_1+(1+\eps)e_2;~e_2e_2=(1+\eps)^{-1}e_1+e_2$ &  &\\\hline
          
    \end{tabular}
    \\~\\ The case where $\fg=\mathbf{L^1_{2|0}}$.
\end{center}

\subsubsection{The case where $\fg=\bf{L^2_{2|0}}$} In that case, $\fg$ is abelian. Therefore, any compatible left symmetric product on $\fg$ has to be commutative and associative (see \cite[Lemma 3.1]{Bu3}). Therefore, we can use \cite[Theorem 1]{GK}.

\begin{center}
	\begin{tabular}{|c|c|c|}
		\hline
		Left-symmetric product on $\fg$ & Conditions& Left-alternative \\\hline

		$e_1e_1=e_1$ & None &Yes\\\hline
		$e_1e_1=e_2$ & None& Yes\\\hline
		$e_1e_1=e_1;~e_2e_2=e_2$ & None& Yes \\\hline
		$e_1e_1=e_1;~e_1e_2=e_2;~e_2e_1=e_2$ & None & Yes  \\\hline
		$e_1e_1=e_1;~e_1e_2=e_2;~e_2e_1=e_2;~e_2e_2=e_1+e_2$ & None &Yes\\\hline
	\end{tabular}
	\\~\\ The case where $\fg=\mathbf{L^2_{2|0}}$.
\end{center}

\section{Lagrangian extensions for Lie superalgebras}\label{sectionlagrange}

\subsection{The language of connections}

We will need the concept of a connection on a Lie (super)algebra in order to write down a Lagrangian extension. We define connections as non-associative products on a superalgebra with purely algebraic definitions. See \cite{BC, BM, MKSV} for the non-super case. 

\subsubsection{Connections, curvature and torsion}\label{cct} Let $(\fg,[\cdot,\cdot],s)$ be a Lie superalgebra. A \textit{connection} on $\fg$ is an {\bf even} linear map 
\[
\nabla:\fg\rightarrow\End(\fg), \qquad x\mapsto \nabla_x.
\] 

The \textit{torsion of the connection} $\nabla$ is given by a pair of maps $(T,U)$, where 
\[
T: \fg \times \fg \rightarrow \fg, \text{ and } U: \fg_{\bar 1}  \rightarrow \fg, \text{ are defined as follows:}
\]
\begin{align}
    T(x,y):=\nabla_{x}(y)+\nabla_y(x)+[x,y],~\forall x,y\in \fg; \; U(x):=\nabla_x(x)+s(x),~\forall x\in \fg_\od.
\end{align}
The connection $\nabla$ is called \textit{torsion-free} if $(T,U)=(0,0)$.

The \textit{curvature of the connection} $\nabla$ is given by a pair of maps $(R,S)$, where
\[
R: \fg \times \fg \rightarrow \End(\fg), \text{ and } S:\fg_{\bar 1}  \rightarrow \End(\fg), \text{ are defined as follows:}
\]
\begin{align}
    R(x,y):=\nabla_{[x,y]}+[\nabla_x,\nabla_y],~\forall x,y\in \fg; \;     S(x):=\nabla_{s(x)}+\nabla_x^2,~\forall x\in \fg_\od.
\end{align}

The connection $\nabla$ is called \textit{flat} if $(R,S)=(0,0)$. In that case, the map $x\mapsto\nabla_x$ is a representation of the Lie superalgebra $\fg$ into $\End(\fg)$.

It is worth mentioning that the data of a connection $\nabla$ on $\fg$ is equivalent to the existence of a bilinear operation $\fg\times \fg \rightarrow \fg$ given by $x\ast y:=\nabla_x(y),~\forall x,y\in \fg$. Moreover, the operation $*$ is a left-symmetric product compatible with the bracket $[\cdot,\cdot]$ and the squaring $s$ on $\fg$ (in the sense of Proposition \ref{LSSA-Lie}) if and only if the connection $\nabla$ is flat and torsion-free.

\subsubsection{Covariant derivatives.} Let $(\fg,[\cdot,\cdot],s)$ be a Lie superalgebra and let $\nabla$ be a connection on $\fg$. For $z,w\in \fg$, the \textit{covariant derivation of the curvature} is defined by
\begin{align}
    (\nabla_z R)(x,y)w:=&\nabla_z\bigl(R(x,y)w\bigl)+R\bigl(\nabla_z(x),y\bigl)w\\\nonumber~&+R\bigl(x,\nabla_z(y)\bigl)w+R(x,y)(\nabla_zw),~\forall x,y\in \fg.\\
    (\nabla_z S)(x)w:=&\nabla_z\bigl(S(x)w\bigl)+R\bigl(\nabla_z(x),x\bigl)w+S(x)(\nabla_z w),~\forall x\in \fg_\od.
\end{align}
Similarly, the \textit{covariant derivation of the torsion} is defined by
\begin{align}
    (\nabla_z T)(x,y):=&\nabla_z\bigl(T(x,y)\bigl)+T\bigl(\nabla_z(x),y\bigl)+T\bigl(x,\nabla_z(y)\bigl),~\forall x,y\in \fg.\\
    (\nabla_z U)(x):=&\nabla_z\bigl(U(x)\bigl)+T\bigl(\nabla_z(x),x\bigl),~\forall x\in \fg_\od.
\end{align}

A connection $\nabla$ on $\fg$ is called \textit{parallel} if $(\nabla_z T,\nabla_z U)=(0,0),~\forall z\in \fg.$

\sssbegin{Remark}
Over a field of characteristic zero, it is well known (see e.g. \cite{MKL})
 that there is a one-to-one correspondence between Lie algebras equipped with a flat parallel connection and post-Lie algebras. This remains true for Lie superalgebras in characteristic $2$, see Proposition \ref{postliecorrespondence}. We initiate the study of post-Lie superalgebras in characteristic $2$ in Appendix \ref{postlie}.
\end{Remark}

\subsubsection{Dual connections} Let $(\fg,[\cdot,\cdot],s)$ be a Lie superalgebra and let $\nabla$ be a connection on $\fg$. The \textit{dual connection} of $\nabla$ is defined by
\begin{equation}
        \onabla_x(y)=\nabla_y(x)+[x,y],~\forall x,y\in \fg.
\end{equation}
Notice that we have $\overline{\onabla}=\nabla$. We denote the curvature and the torsion of $\onabla$ by $(\overline{R},\overline{S})$ and $(\overline{T},\overline{U})$, respectively.



\sssbegin{Lemma}
    Let $(\fg,[\cdot,\cdot],s)$ be a Lie superalgebra and let $\nabla=0$ be the zero connection on $\fg$. Then, $(\overline{R},\overline{S})=(0,0)$ and $\onabla$ is the adjoint representation of $\fg$.
\end{Lemma}

\begin{proof}
    Let $x,y,z\in \fg$.
    \begin{align*}
        \overline{R}(x,y)(z)&=\onabla_{[x,y]}(z)+\onabla_x\circ\onabla_y(z)+\onabla_y\circ\onabla_x(z)\\
        &=[[x,y],z]+\onabla_x([y,z])+\onabla_y([x,z])\\
        &=[[x,y],z]+[[y,z],x]+[[z,x],y]\\&=0.
    \end{align*}
    Let $x\in \fg_\od$ and $z\in \fg$. Then
        $\overline{S}(x)(z)=\onabla_{s(x)}(z)+\onabla_x^2(z)=[s(x),z]+[x,[x,z]]=0.$
        Therefore, $\onabla$ is a representation.
\end{proof}

\sssbegin{Proposition}\label{propA}
Let $(\fg,[\cdot,\cdot],s)$ be a Lie superalgebra and let $\nabla$ be a connection on $\fg$. Then
\begin{align}
    (\nabla_z T)(x,y)&=\overline{R}(x,y)(z)+R(y,z)(x)+R(z,x)(y),~\forall x,y,z\in \fg;\\\label{id2}
    (\nabla_z U)(x)&=\overline{S}(x)(z)+R(z,x)(x),~\forall x\in \fg_\od,~\forall z\in \fg.
\end{align}
\end{Proposition}
\begin{proof}
We prove Identity \eqref{id2}. Let $z\in \fg$ and $x\in \fg_\od$. We have
    \begin{align*}  
        (\nabla_zU)(x)&=\nabla_z(U(x))+T(\nabla_z(x),x)\\
                        &=\nabla_z\bigl(\nabla_x(x)+s(x)\bigl)+\nabla_{\nabla_z(x)}(x)+\nabla_x\circ\nabla_z(x)+[\nabla_z(x),x]\\
                        &=\nabla_z\circ\nabla_x(x)+\nabla_z\circ s(x)+\onabla_x\circ\nabla_z(x)+\nabla_x\circ\nabla_z(x)\\
                        &=\nabla_z\circ\nabla_x(x)+\nabla_x\circ\nabla_z(x)+\nabla_z\circ s(x)+\onabla_x\circ\bigl(\onabla_x(z)+[z,x]\bigl)\\
                        &=\nabla_z\circ\nabla_x(x)+\nabla_x\circ\nabla_z(x)+\onabla_{s(x)}(z)+[z,s(x)]+\onabla_x^2(z)+\nabla_{[z,x]}(x)+[x,[z,x]]\\
                        &=\onabla_{s(x)}(z)+\onabla_x^2(z)+\nabla_{[z,x]}(x)+[\nabla_z,\nabla_x](x)\\
                        &=\overline{S}(x)(z)+R(z,x)(x).
    \end{align*}

\end{proof}

\sssbegin{Corollary}
Let $(\fg,[\cdot,\cdot],s)$ be a Lie superalgebra and let $\nabla$ be a flat connection on $\fg$. Then
    \begin{equation}
        (\overline{R},\overline{S})=(0,0)\text{ if and only if } (\nabla_z T,\nabla_z U)=(0,0),\forall z\in \fg.
    \end{equation}
\end{Corollary}

\sssbegin{Proposition}\label{propB}
        Let $(\fg,[\cdot,\cdot],s)$ be a Lie superalgebra and let $\nabla$ be a connection on $\fg$. Then
        \begin{align}
            \underset{x,y,z}{\circlearrowleft}T(x,T(y,z))=&\underset{x,y,z}{\circlearrowleft}R(x,y)(z)+\underset{x,y,z}{\circlearrowleft}\overline{R}(x,y)(z),~\forall x,y,z\in \fg;\\\label{propB2}
            T\bigl(x,T(x,y)\bigl)+T\bigl(U(x),y\bigl)=&~S(x)(y)+R(y,x)(x)\\\nonumber&+\overline{S}(x)(y)+\overline{R}(y,x)(x),~\forall y\in \fg,~\forall x\in \fg_\od.
        \end{align}
\end{Proposition}

\begin{proof}
Let $y\in \fg$ and $x\in \fg_\od.$ We have
\begin{align*}
T\bigl(x,T(x,y)\bigl)=&\nabla_x^2(y)+\nabla_x\nabla_y(x)+\nabla_x([x,y])\\
                        &~+\nabla_{\nabla_x(y)}(x)+\nabla_{\nabla_y(x)}(x)+\nabla_{[x,y]}(x)\\
                        &~+[x,\nabla_x(y)]+[x,\nabla_y(x)]+[x,[x,y]],                      
\end{align*}
and
\begin{align*}
    T\bigl(U(x),y\bigl)=&\nabla_{\nabla_x(x)}(y)+\nabla_{s(x)}(y)+\nabla_y\circ\nabla_x(x)+\nabla_y(s(x))+[\nabla_x(x),y]+[s(x),y].
\end{align*} On the other hand,
\begin{align*}
S(x)y+R(y,x)x&=\nabla_{s(x)}(y)+\nabla^2_x(y)+\nabla_{[y,x]}(x)+\nabla_y\circ\nabla_x(x)+\nabla_x\circ\nabla_y(x),~\text{ and }\\
\overline{S}(x)y+\overline{R}(y,x)x&=\onabla_{s(x)}(y)+\onabla^2_x(y)+\onabla_{[y,x]}(x)+\onabla_y\circ\onabla_x(x)+\onabla_x\circ\onabla_y(x).
\end{align*}
Using the identities
\begin{align*}
\nabla_{\nabla_x(x)}(y)&=\nabla_{\overline{\nabla}_x(x)}(y)=\onabla_y\circ\onabla_x(x)+[\onabla_x(x),y];\\
\nabla_{\nabla_y(x)}(x)&=\nabla_{[x,y]}(x)+\onabla_x\onabla_x(y)+[\onabla_x(y),x];\\
\nabla_{\nabla_x(y)}(x)&=\nabla_{[x,y]}(x)+\onabla_x\onabla_x(y)+[\onabla_x(y),x],
\end{align*}
as well as
$$ \onabla_{s(x)}(y)=\nabla_y(s(x))+[s(x),y]~\text{ and }~\onabla_{[y,x]}(x)=\nabla_x([y,x])+[[y,x],x],$$ we obtain
\begin{align}
            T\bigl(x,T(x,y)\bigl)+T\bigl(U(x),y\bigl)=&S(x)(y)+R(y,x)(x)+\overline{S}(x)(y)+\overline{R}(y,x)(x)\nonumber\\
            &+[\onabla_y(x),x]+[x,\onabla_x(y)]+[\onabla_x(y),x]\label{pro1}\\
            &+[\nabla_y(x),x]+[y,\onabla_x(x)]+[\onabla_x(x),y].\label{pro2}
\end{align}
We have $\eqref{pro1}+\eqref{pro2}=0,$ thus the conclusion.

\end{proof}

\sssbegin{Corollary}
    Let $(\fg,[\cdot,\cdot],s)$ be a Lie superalgebra and let $\nabla$ be a connection on $\fg$. Suppose that $(R,S)=(\overline{R},\overline{S})=(0,0)$. Then, $\bigl(\fg,T(\cdot,\cdot),U(\cdot)\bigl)$ is a Lie superalgebra.                                                                 
\end{Corollary}

\sssbegin{Corollary}[First Bianchi identity]                                               
    Let $(\fg,[\cdot,\cdot],s)$ be a Lie superalgebra and let $\nabla$ be a connection on $\fg$. Then
    \begin{align}
        \underset{x,y,z}{\circlearrowleft}(\nabla_zT)(x,y)=&\underset{x,y,z}{\circlearrowleft}R(z,x)y+\underset{x,y,z}{\circlearrowleft}T(z,T(x,y)),~\forall x,y,z\in \fg;\\\label{Bianchi1.2}
        (\nabla_yU)(x)+(\nabla_xT)(y,x)=&~R(y,x)x+S(x)y\\\nonumber&~+T(y,U(x))+T(x,T(y,x)),~\forall y\in \fg,~\forall x\in \fg_\od.
    \end{align}
\end{Corollary}

\begin{proof}
We prove Identity \eqref{Bianchi1.2}. Let $y\in \fg,~x\in \fg_\od$. Using Proposition \ref{propA}, we have
$$(\nabla_yU)(x)+(\nabla_xT)(y,x)=\overline{S}(x)y+R(y,x)x+\overline{R}(y,x)x+R(x,y)x.$$ Applying \eqref{propB2} to the right-hand side gives \eqref{Bianchi1.2}.
\end{proof}

\sssbegin{Proposition}[Second Bianchi identity]                                               
    Let $(\fg,[\cdot,\cdot],s)$ be a Lie superalgebra and let $\nabla$ be a connection on $\fg$. Then
    \begin{align}
        \underset{x,y,z}{\circlearrowleft}(\nabla_zR)(x,y)=&\underset{x,y,z}{\circlearrowleft}R(z,T(x,y)),~\forall x,y,z\in \fg;\\\label{Bianchi2.2}
        (\nabla_yS)(x)+(\nabla_xR)(x,y)=&~R(x,T(x,y))+R(y,U(x)),~\forall y\in \fg,~\forall x\in \fg_\od.
    \end{align}
\end{Proposition}

\begin{proof}
Direct computations.
\end{proof}

\sssbegin{Proposition}
Let $(V,\ast)$ is a superspace in characteristic $2$ endowed with a bilinear map $\ast:V\times V\rightarrow V$. Recall that the associator is given by
    \begin{equation}
        \Asso(x,y,z):=(x*y)*z-x*(y*z)=\nabla_{\nabla_x(y)}(z)-\nabla_x\circ\nabla_y(z),~\forall x,y,z\in V,
    \end{equation} where $\nabla_x(y)=x*y,~\forall x,y\in V.$\\
Let us denote $[x,y,z]:=\Asso(x,y,z)-\Asso(y,x,z),~\forall x,y,z\in V.$ Then we have
\begin{align}
    [x,y,z]&=R(x,y)(z)+T(x,y)*z=R(x,y)(z)+\nabla_{T(x,y)}(z),~\forall x,y,z\in V;\\
    \Asso(x,x,z)&=S(x)(z)+U(x)*z=S(x)(z)+\nabla_{U(x)}(z),~\forall x\in V_\od,~\forall z\in V.\label{asso-xxz}
\end{align}
\end{Proposition}

\begin{proof}
    We prove Identity \eqref{asso-xxz}. Let $x\in V_\od$ and $z\in V$. We have
    \begin{align*}
    S(x)z+\nabla_{U(x)}(z)&=\nabla_{s(x)}(z)+\nabla^2_x(z)+\nabla_{\nabla_x(x)+s(x)}(z)\\
    &=\nabla^2_x(x)+\nabla_{\nabla_x(x)}(x)\\&=x*(x*z)+(x*x)*z.
    \end{align*}
\end{proof}

\subsection{Lagrangian extensions} To our best knowledge, Lagrangian extensions were introduced by Bordemann (see \cite{Bo}) under the name $T^*$-extensions, dealing with symmetric bilinear forms. The notion was then studied by Baues and Cortes (see \cite{BC}), who name it Lagrangian extensions, in the context of Lie algebras admitting a flat connection and the bilinear form is antisymmetric. They also introduced the notion of Lagrangian cohomology and showed that equivalence classes of Lagrangian extensions are captured by this cohomology. The superization is due to Bouarroudj and Maeda (see \cite{BM}), who used the terminology $T^*$-extensions (resp. $\Pi T^*$-extensions) in the case where the form is even (resp. odd). In this paper, we use the terminology ``Lagrangian extension" to mention both $T^*$-extensions and $\Pi T^*$-extensions.

 Let $(\fg,[\cdot,\cdot],s)$ be a Lie superalgebra equipped with a non-degenerate closed ortho-orthogonal or periplectic bilinear form $\w$ and $I\subset \fg$ be a homogeneous ideal such that $[I,I^{\perp}]=0.$ Let $\fh:=\fg/I^{\perp}$. There is a non-degenerate bilinear pairing between $\fh$ and $I$ by declaring
 \begin{equation}\label{pairing}
    \w_\fh(x,u):=\w(\tilde{x},u),~~\forall x\in \fh,~\forall u\in I,
 \end{equation}
where $\tilde{x}$ is a lift of $x$ to $\fg$. This expression is well-defined, see \cite{BM}. Notice that
\begin{equation}
    \widetilde{[u,v]}-[\tilde{u},\tilde{v}]\in I^\perp,\quad\text{and}\quad \widetilde{s(x)}-s(\tilde{x})\in (I^\perp)_\ev,\quad \forall u,v\in \fg,~\forall x\in \fg_\od.
\end{equation}
Moreover, it is always possible to find a lift that preserves the parity.

The following proposition has been proved in \cite{BC} for Lie algebras, and in \cite{BM} for Lie superalgebras but $p\not =2$. Here, we complete the picture for $p=2$.  This will be useful in Section \ref{sec:converse}.

\sssbegin{Proposition}\label{prop:ppty}
 Let $(\fg,[\cdot,\cdot],s)$ be a Lie superalgebra equipped with a non-degenerate closed ortho-orthogonal or periplectic bilinear form $\w$ and $I\subset \fg$ be a homogeneous ideal such that $[I,I^{\perp}]=0.$ Let $\fh:=\fg/I^{\perp}$ and let $\w_\fh$ be the pairing defined in Eq. \eqref{pairing}. If $x\in \fh$, we will denote by $\tilde{x}$ a lift of $x$ to $\fg$. 
    \begin{enumerate}[(i)]
        \item The map $\ad_I:\fh\rightarrow \End(I),\quad x\mapsto\ad_{\tilde{x}}$ is a representation, which induces representations of the spaces $I^*$ and $\Pi(I^*)$ given by
            \begin{align}\label{mod1}
                \ad^*_I: \fh&\rightarrow \End(I^*),~\quad\quad \text{ where }\quad  \ad_I^*(x)(\xi):=\xi\circ\ad_{\tilde{x}};\\\label{mod2}
                \Pi\ad^*_I: \fh&\rightarrow \End(\Pi(I^*)),\quad \text{ where }\quad  \Pi\ad_I^*(x)(\Pi(\xi)):=\Pi\circ\xi\circ\ad_{\tilde{x}};
            \end{align}
        \item In the case where $\w$ is ortho-orthogonal, there is an isomorphism
        $$\varphi:\fh\rightarrow I^*,\quad x\mapsto \w_\fh(x,\cdot),$$ and the map $\varphi$ is a one-cocycle in $XZ^1_\ev(\fh,I^*)$, where the module structure is given by \eqref{mod1};
        \item In the case where $\w$ is periplectic, there is an isomorphism
        $$\psi: \fh\rightarrow \Pi(I^*),\quad x\mapsto \Pi\circ \w_\fh(x,\cdot),
        $$ and the map $\psi$ is a one-cocycle in $XZ^1_\ev(\fh,\Pi(I^*))$, where the module structure is given by \eqref{mod2};
        \item There is a torsion-free flat connection $\nabla$ on $\fh$ defined by
            \begin{equation}
                \w_\fh\bigl(\nabla_x(y),u\bigl):=\w\bigl(\tilde{y},[\tilde{x},u]\bigl),~\forall x,y\in \fh,~\forall u\in I.
            \end{equation}
    \end{enumerate}
 
\end{Proposition}

\begin{proof} We will only prove identities involving the squaring. For identities involving the brackets, see \cite{BC,BM}.
\begin{enumerate}[$(i)$]
    \item The map $\ad_I$ is well-defined. Let $x\in \fh,~y\in I$. Since $[I,I^\perp]=0$, we have 
    $$\ad_I\bigl(\widetilde{s(x)}\bigl)(y)=[s(\tilde{x})+(I^\perp)_\ev,y]=[s(\tilde{x}),y]=\ad_I(\tilde{x})(y).$$
Let $\xi\in I^*$. We have
$$\ad^*_I\bigl(s(x)\bigl)(\xi)=\xi\circ\ad_I\bigl(\widetilde{s(x)}\bigl)=\xi\circ\ad_I^2(\tilde{x})=\ad^*_I\bigl(\xi\circ\ad_I(\tilde{x}) \bigl)=(\ad^*_I)^2(x)(\xi).$$

    \item Let $x\in \Ker(\varphi).$ Then, $\w_H(x,v)=0~\forall v\in I$. Therefore, $\tilde{x}\in I^\perp$ and $x=0\in \fh=\fg/I^\perp$. 
    Since $\dim(\fh)=\codim(I^\perp)=\dim(I^*)$ the map $\varphi$ is an isomorphism. Next, let $x\in \fh$ and $v\in I$.
    \begin{align*}
            \bigl(x\cdot\varphi(x)\bigl)(v)+\varphi\bigl(s(x)\bigl)(v)&=\varphi(x)\circ\ad_{\tilde{x}}(v)+\w_H\bigl(s(x),v\bigl)\\
            &=\w\bigl(\tilde{x},[\tilde{x},v]\bigl)+\w\bigl(\widetilde{s(x)},v\bigl)\\
            &=\w\bigl(\tilde{x},[\tilde{x},v]\bigl)+\w\bigl(s(\tilde{x}),v\bigl)=0.
    \end{align*} The 1-cocycle condition involving the bracket follows as in \cite{BM}. We conclude that $\varphi$ is a 1-cocycle.
    \item Similar to $(ii)$.
    \item The same computation as in  \cite{BM} shows that the maps $T=0$ and $R=0$. Only the maps $S$ and $U$ will be computed here. Let $x\in H_\od,~y\in \fh,~u\in I$.
        \begin{align*}
            \w_\fh\bigl(S(x)y,u\bigl)&:=\w_\fh\bigl(\nabla_{s(x)}(y)+\nabla^2_x(y),u\bigl)\\
            &=\w\bigl(\tilde{y},[\widetilde{s(x)},u]\bigl)+\w\bigl(\widetilde{\nabla_x(y)},[\tilde{x},u]\bigl)\\
           & =\w\bigl(\tilde{y},[s(\tilde{x})+I,u]\bigl)+\w_\fh\bigl(\nabla_x(y),[\tilde{x},u]\bigl)\\
            &=\w\bigl(\tilde{y},[s(\tilde{x}),u]\bigl)+\w\bigl(\tilde y,[\tilde{x},[\tilde{x},u]]\bigl)=0.
        \end{align*}
        We deduce that $S(x)y=0,~\forall x,y\in \fh$. Similarly, we have
        \begin{align*}
           \w_\fh\bigl(U(x),u\bigl)&=\w_\fh\bigl(\nabla_x(x)+s(x),u\bigl)\\
           &=\w\bigl(\tilde{x},[\tilde{x},u]\bigl)+\w\bigl(\widetilde{s(x)},u\bigl)=0.
        \end{align*}   
        
        Therefore, $U(x)=0,~\forall x\in H_\od$.
\end{enumerate}

\end{proof}

\subsubsection{Strong polarization} Let $(\fg,\w)$ be a quasi-Frobenius Lie superalgebra. A strong polarization of $(\fg,\w)$ is a decomposition $\fg=\fa\oplus N$ as vector superspaces, where $\fa$ is a homogeneous Lagrangian ideal of $\fg$ and $N$ is a Lagrangian subspace. Recall that a subspace $\fl\in\fg$ is called Lagrangian if $\fl=\fl^{\perp}$, with $\fl^{\perp}$ the orthogonal with respect to $\w$. A strongly polarized quasi-Frobenius Lie superalgebra will be denoted by $(\fg,\w,\fa,N)$.

A morphism of strongly polarized quasi-Frobenius Lie superalgebra is a Lie superalgebras morphism $\phi:(\fg,\w,\fa,N)\rightarrow(\fg',\w',\fa',N')$ such that $\w=\phi^*\w'$, $\phi(\fa)=\fa'$ and $\phi(N)=N'$.\\

Let $(\fh,[\cdot,\cdot]_\fh,s_\fh)$ be a Lie superalgebra endowed with a torsion-free flat connection $\nabla:\fh\rightarrow\End(\fh)$. Recall that since $\nabla$ is flat, it is a representation of $\fh$. Following \cite{BM}, we define the dual representations $\rho:\fh\rightarrow\End(\fh^*)$ and $\chi:\fh\rightarrow\End(\Pi(\fh^*))$ as follows:
\begin{equation}\label{module1}
    \rho:\fh\rightarrow\End(\fh^*),\quad x\mapsto \rho(x),\quad \text{where } \rho(x)(\xi):=\xi\circ\nabla_x,~\forall x\in \fh,~\forall \xi\in \fh^*.
\end{equation}
\begin{equation}\label{module2}
    \chi:\fh\rightarrow\End(\Pi(\fh^*)),\quad x\mapsto \chi(x):=\Pi\circ\rho(x)\circ\Pi.
\end{equation}

\sssbegin{Lemma}
    The maps $\rho$ and $\chi$ are indeed representations.
\end{Lemma}
\begin{proof}
We will only prove identities involving the squaring. See \cite[Lemma 6.1.1]{BM}, for a proof of the other identities. Let $x\in \fh_\od$ and $\xi\in \fh^*$.
Since the connection is flat, we have
$$\rho(s_\fh(x))(\xi)=\xi\circ\nabla_{s_\fh(x)}=\xi\circ\nabla^2_x=\rho(x)^2(\xi).$$ Therefore, $\rho$ is a representation. Furthermore,
$$\chi(s_\fh(x))=\Pi\circ\rho(s_\fh(x))\circ\Pi=\Pi\circ\rho(x)^2\circ\Pi=\Pi\circ\rho(x)\circ\Pi\circ \Pi\circ\rho(x)\circ\Pi=\chi(x)^2, $$ so $\chi$ is a representation as well.
\end{proof}

Let $(\alpha,\gamma)\in XZ^2(\fh,\fh^*)_\ev$ (resp. $(\beta,\theta)\in XZ^2(\fh,\Pi(\fh^*))_\ev$) be 2-cocycles. We will construct an ortho-orthogonal (resp. a periplectic) Lie superalgebra structure on $\fg:=\fh\oplus\fh^*$ (resp. on $\fg:=\fh\oplus\Pi(\fh^*)$). The spaces $\fh$ and $\Pi(\fh^*)$ are $\fh$-modules using the representations \eqref{module1} and \eqref{module2}.

\underline{On the space $\fg:=\fh\oplus\fh^*$.} The brackets and squaring are defined as follows:
\begin{align}
    [x,y]_\fg&:=[x,y]_\fh+\alpha(x,y),\quad [x,\xi]_\fg:=\rho(x)(\xi),\quad\forall x,y\in \fh,~\forall \xi\in \fh^*.\\
    s_\fg(x+\xi)&:=s_\fh(x)+\gamma(x)+\rho(x)(\xi),\quad\forall x\in \fh_\od,~\forall\xi\in\fh^*_\od.
\end{align}  We define an ortho-orthogonal form as follows:
\begin{equation}\label{formw}
    \w(x+\xi,y+\zeta)=\xi(y)+\zeta(x),\quad \forall x+\xi,y+\zeta\in \fg.
\end{equation}

\underline{On the space $\fg:=\fh\oplus\Pi(\fh^*)$.} The brackets and squaring are defined as follows:
\begin{align}
    [x,y]_\fg&:=[x,y]_\fh+\beta(x,y),\quad [x,\Pi(\xi)]_\fg:=\chi(x)(\Pi(\xi)),\quad\forall x,y\in \fh,~\forall \Pi(\xi)\in \Pi(\fh^*).\\
    s_\fg(x)&:=s_\fh(x)+\theta(x)+\chi(x)(\Pi(\xi)),\quad\forall x\in \fh_\od,~\forall\Pi(\xi)\in\Pi(\fh^*_\od).
\end{align}
We define a periplectic form as follows:
\begin{equation}\label{formk}
    \kappa(x+\Pi(\xi),y+\Pi(\zeta))=\xi(y)+\zeta(x),\quad \forall x+\xi,y+\zeta\in \fg.
\end{equation}

Let us prove Jacobi Identities involving the squaring. Let $x\in \fh_\od,~y\in\fh$ and let $\xi\in\fh^*$. Recall that the second part of the 2-cocycle condition \eqref{2cocycle} for the pair $(\alpha,\gamma)$ reads
\begin{equation}\label{2cocyh}
    \alpha(x,y)\circ\nabla_x+\gamma(x)\circ\nabla_y=\alpha(s_\fh(x),y)+\alpha(x,[x,y]_\fh).
\end{equation}
We have
\begin{align*}
    [x,[x,y]_\fg]_\fg&=[x,[x,y]_\fh]_\fh+\alpha(x,[x,y]_\fh)+\alpha(x,y)\circ\nabla_x\\
                    &=[s_\fh(x),y]_\fh+\alpha(s_\fh(x),y)+\gamma(x)\circ\nabla_y\quad (\text{using\eqref{2cocyh}})\\
                    &=[s_\fh(x)+\gamma(x),y]_\fg=[s_\fg(x),y]_\fg.
\end{align*}
Furthermore, since $\nabla$ is flat, we have $\nabla^2_x=\nabla_{s_\fh(x)}~\forall x\in\fh_\od$. Therefore, 
$$[x,[x,\xi]_\fg]_\fg=[x,\xi\circ\nabla_x]_\fg=\xi\circ\nabla^2_x= \xi\circ\nabla_{s_\fh(x)}=[s_\fh(x),\xi]_\fg=[s_\fh(x)+\gamma(x),\xi]_\fg=[s_\fg(x),\xi]_\fg.  $$

\sssbegin{Theorem}\label{sensdirect} Let $(\fh,\nabla)$ be a Lie superalgebra equipped with a flat and torsion-free connection $\nabla$ and let $(\alpha,\gamma)\in XZ^2(\fh,\fh^*)_\ev$ (resp. $(\beta,\theta)\in XZ^2(\fh,\Pi(\fh^*))_\ev$) be 2-cocycles.
    \begin{enumerate}
        \item The form $\w$ on $\fg=\fh\oplus\fh^*$ defined in \eqref{formw}  is closed if and only if
            \begin{align}\label{condalphagamma1}
                \underset{x,y,z}{\circlearrowleft}\alpha(x,y)(z)&=0,~\forall x,y,z\in\fh;\\\label{condalphagamma2}
                \gamma(x)(y)+\alpha(x,y)(x)&=0,~\forall x\in \fh_\od,~\forall y\in \fh.
            \end{align}
            In this case, one can canonically define a strongly polarized quasi-Frobenius Lie superalgebra $(\fg,\w,\fh^*,\fh)$, where $\fg=\fh\oplus\fh^*$, called $T^*$-extension of $(\fh,\nabla)$. 
             \item The form $\kappa$ on $\fg=\fh\oplus\Pi(\fh^*)$ defined in \eqref{formk} is closed if and only if
            \begin{align}\label{condbetatheta1}
                \underset{x,y,z}{\circlearrowleft}\beta(x,y)(z)&=0,~\forall x,y,z\in\fh;\\
                \theta(x)(y)+\beta(x,y)(x)&=0,~\forall x\in \fh_\od,~\forall y\in \fh.\label{condbetatheta2}
            \end{align}
             In this case, one can canonically define a strongly polarized quasi-Frobenius Lie superalgebra $(\fg,\kappa,\Pi(\fh^*),\fh)$, where $\fg=\fh\oplus\Pi(\fh^*)$, called $\Pi T^*$-extension of $(\fh,\nabla)$.
    \end{enumerate}
\end{Theorem}

\begin{proof}
    We will prove Eq. \eqref{condalphagamma2}. Let $x+\xi\in\fg_\od$ and $y+\zeta\in\fg.$
    \begin{align*}
        &\w\bigl(s_\fg(x+\xi),y+\zeta\bigl)+\w\bigl(x+\xi,[x+\xi,y+\zeta]_\fg\bigl)\\
        =&~\w(s_\fh(x)+\gamma(x)+\xi\circ\nabla_x,y+\zeta)+\alpha(x,y)(x)+\xi\circ\nabla_y(x)+\zeta\circ\nabla_x(x)+\xi([x,y]_\fh)\\
        =&~\gamma(x)(y)+\underline{\xi\circ\nabla_x(y)}+\overline{\zeta\bigl(s_\fh(x)\bigl)}+\alpha(x,y)(x)+\underline{\xi\circ\nabla_y(x)}+\overline{\zeta\circ\nabla_x(x)}+\underline{\xi\bigl([x,y]_\fh\bigl)}.
    \end{align*}
    Since $\nabla$ is torsion free, the underlined and overlined terms cancel. Thus the conclusion.
\end{proof}

\subsubsection{Converse of Theorem \ref{sensdirect}}\label{sec:converse} We will show that every strongly polarized quasi-Frobenius Lie superalgebra can be obtained as a $T^*$-extension or $\Pi T^*$-extension of a smaller Lie superalgebra equipped with a flat torsion-free connection. Let $(\fg,\w_\fg,\fa,N)$ be a strongly polarized quasi-Frobenius Lie superalgebra. We recall that $\fg=\fa\oplus N$ as vector superspaces, where $\fa$ is a Lagrangian ideal and $N$ is a Lagrangian subspace. Let $\fh:=\fg/\fa$.

    Since $\fa$ is Lagrangian, there is a flat torsion free connection $\nabla$ on $\fh$ given by Proposition \ref{prop:ppty} and defined by
    \begin{equation}\label{inducednabla}
        \w_\fh(\nabla_u(v),a)=\w_\fg(\tilde{v},[\tilde{u},a]),~\forall u,v\in\fh,~\forall a\in \fa.
    \end{equation}
We define an $\fh$-module structure on $\fa$ by
\begin{equation}\label{adrep}
    \ad_{\fh,\fa}:\fh\rightarrow \text{End}(\fa,\fa),\quad u\mapsto[\tilde{u},\cdot].
\end{equation}
We denote $p_\fa:\fg\rightarrow\fa$ and $p_N:\fg\rightarrow N$ the natural projections. We define
\begin{align}
    \tilde{\alpha}(u,v)&:=p_\fa\bigl([p_N(\tilde{u}),p_N(\tilde{v})]_\fg\bigl),\quad\forall u,v\in\fh\\
    \tilde{\gamma}(u)&:=p_\fa\bigl(s_\fg(p_N(\tilde{u}))\bigl),\quad\forall u\in \fh_\od.
\end{align}

\sssbegin{Lemma}\label{lemmacocycle}
    We have $(\tilde{\alpha},\tilde{\gamma})\in XZ^2_\ev(\fh,\fa)$.
\end{Lemma}

\begin{proof}
    Recall that (see \cite[Lemma 6.1.5]{BM})
    \begin{equation}\label{615}
        u\cdot\tilde{\alpha}(v,w)+\tilde{\alpha}(u,[v,w]_\fg)=p_{\fa}\bigl(p_N(\tilde{u}),[p_N(\tilde{v}),p_N(\tilde{w})]_\fg\bigl),\quad \forall u,v,w\in \fh.
    \end{equation}
    Let $u\in\fh_\od$ and $w\in\fh$. We have
    \begin{align}  \tilde{\alpha}\bigl(s_\fh(u),w\bigl)=p_\fa\bigl(\bigl[p_N(\widetilde{s_\fh(u)}),p_N(\tilde{w})\bigl]_\fg\bigl)=p_\fa\bigl(\bigl[p_N(s_\fg(\tilde{u})), p_N(\tilde{w})\bigl]_\fg\bigl)
    \end{align} and
    \begin{align}
        w\cdot\tilde{\gamma}(u)&=\bigl[ \tilde{w},p_\fa\bigl(s_\fg(p_N(\tilde{u}))\bigl)\bigl]_\fg=p_\fa\bigl(\bigl[p_N(\tilde{w}),p_\fa(s_\fg(p_N(\tilde{u})))\bigl]_\fg\bigl).
    \end{align}
Therefore, 
  \begin{align}
        w\cdot\tilde{\gamma}(u)+\tilde{\alpha}\bigl(s_\fh(u),w\bigl)&=p_\fa\bigl(\bigl[p_N(\tilde{w}),p_\fa(s_\fg(p_N(\tilde{u})))+p_N(\widetilde{s_\fh(u)}) \bigl]_\fg\bigl)\nonumber\\ &=p_\fa\bigl(\bigl[p_N(\tilde{w}),\underset{=s_\fg(p_N(\tilde{u}))}{\underbrace{p_\fa(s_\fg(p_N(\tilde{u})))+p_N(s_\fg(p_N(\tilde{u})))}}\\&\quad\text{ }+\underset{=0}{\underbrace{p_N(s_\fg(p_\fa(\tilde{u}))+[p_\fa(\tilde{u}),p_N(\tilde{u})]_\fg)}}\bigl]_\fg\bigl)\nonumber\\\label{bibis}
        &=p_\fa\bigl(\bigl[p_N(\tilde{w}),s_\fg(p_N(\tilde{u}))\bigl]_\fg\bigl)
        =p_{\fa}\bigl(p_N(\tilde{u}),[p_N(\tilde{u}),p_N(\tilde{w})]_\fg\bigl).
    \end{align} Using \eqref{615} and \eqref{bibis}, we obtain that
    \begin{equation}
         u\cdot\tilde{\alpha}(u,w)+\tilde{\alpha}(u,[u,w]_\fg)+ w\cdot\tilde{\gamma}(u)+\tilde{\alpha}\bigl(s_\fh(u),w\bigl)=0,~\forall u\in \fh_\od,~\forall w\in \fh.
    \end{equation}
\end{proof}

\sssbegin{Theorem}\label{conversethm}
Let $(\fg,\w,\fa,N)$ be a strongly polarized quasi-Frobenius Lie superalgebra and let $(\fh,\nabla)$ be the quotient algebra equipped with the flat torsion-free connection described in  \eqref{inducednabla}. 
    \begin{enumerate}[$(i)$]
        \item If $\w$ is ortho-orthogonal, there exists $(\alpha,\gamma)\in XZ^2_\ev(\fh,\fh^*)$ satisfying \eqref{condalphagamma1} and \eqref{condalphagamma2} such that $(\fg,\w,\fa,N)$ is isomorphic to the $T^*$-extension of $(\fh,\nabla)$ by $(\alpha,\gamma)$.
        \item If $\w$ is periplectic, there exists $(\beta,\theta)\in XZ^2_\ev(\fh,\Pi(\fh^*))$ satisfying \eqref{condbetatheta1} and \eqref{condbetatheta2} such that $(\fg,\w,\fa,N)$ is isomorphic to the $\Pi T^*$-extension of $(\fh,\nabla)$ by $(\beta,\theta)$.
    \end{enumerate}
\end{Theorem}

\begin{proof}
    In the proof, we will denote an ortho-orthogonal (resp. periplectic) form on $\fg$ by $\w$ (resp. $\kappa$). Consider the maps
    \begin{equation}
        i_\w:\fa\rightarrow \fh^*,~i_\w(a)(u):=\w(a,\tilde{u}),\quad  i_\kappa:\fa\rightarrow \fh^*,~i_\kappa(a)(u):=\kappa(a,\tilde{u}).
    \end{equation}
    Recall (see \cite[Lemma 6.1.7]{BM}) that the representations $\ad_{\fh,\fa}$ (see Eq. \eqref{adrep}) and $\rho$ (resp. $\chi$) are equivalent, i.e., for all $u\in\fh$, we have
    \begin{align}\label{eqrepr}
        \rho(u)\circ i_\w&= i_\w\circ \ad_{\fh,\fa}(u),\\
        \chi(u)\circ\Pi\circ i_\kappa &= \Pi\circ i_\kappa\circ \ad_{\fh,\fa}(u).
    \end{align}
    Let $p_\fh:\fg\rightarrow\fh=\fg/\fa$ the canonical projection. We denote by $p_\fh:N\rightarrow\fh$ the isomorphism of vector superspaces given by the restriction of $p$ on the subspace $N$. We will prove the Theorem in the ortho-orthogonal case first. Consider the pair $(\alpha,\gamma)$ given by  
    \begin{equation}\label{cocyclestrandardortho}
    \alpha=i_\w\circ\tilde{\alpha} \text{ and } \gamma=i_\w\circ\tilde{\gamma}.
    \end{equation} 
    Using \eqref{eqrepr} and Lemma \ref{lemmacocycle}, we have that $(\alpha,\gamma)\in XZ^2_\ev(\fh,\fh^*)$. Any $x\in\fg$ can be written $x=x_N+x_\fa$, with $x_N=p_N(x)\in N$ and $x_\fa=p_\fa(x)\in \fa$.  Let us define a map $$\Phi:=p_\fh\oplus i_\w:\fg\rightarrow\fh\oplus\fh^*,~(x=x_N+x_\fa)\longmapsto p_\fh(x_N)+\w(x_\fa,\cdot).$$ We will show that $\Phi$ is an isomorphism of Lie superalgebras. We only present the computations involving the squarings. The remaining computations can be found in \cite[Theorem 6.1.4]{BM}. Let $x_N\in N_\od$. Notice that $x_N+p_N\bigl(\widetilde{p_\fh(x_N)}\bigl)\in\fa$. Then,
    \begin{align*}
        \Phi\bigl(s_\fg(x_N)\bigl)&=\Phi\bigl(p_N(s_\fg(x_N))+p_\fa(s_\fg(x_N))\bigl)\\
        &=p_\fh\bigl(p_N(s_\fg(x_N))\bigl)+i_\w\bigl(p_\fa(s_\fg(x_N))\bigl)\\
        &=p_\fh(s_\fg(x_N))+i_\w\bigl(p_\fa\bigl(s_\fg(p_N(\widetilde{p_\fh(x_N)})+\fa)    \bigl)\bigl)\\
        &=s_\fh(p_\fh(x_N))+\gamma(p_\fh(x_N))\\
        &=s_{\fh\oplus\fh^*}\bigl(\Phi(x_N)\bigl).
    \end{align*}
    Now, let $x_\fa\in\fa_\od$. Notice that since $\fa$ is a Lagrangian ideal, $i_\w(s_\fg(x_\fa))=0$. Moreover, $s_{\fh\oplus\fh^*}(i_\w(x_\fa))=0$ by definition. Therefore, $$\Phi\bigl(s_\fg(x_\fa)\bigl)+s_{\fh\oplus\fh^*}\bigl(\Phi(x_\fa)\bigl)=i_\w(s_\fg(x_\fa))+s_{\fh\oplus\fh^*}(i_\w(x_\fa))=0.$$ Now, for all $x\in\fh_\od$, we have
    \begin{align*}
        \Phi\bigl(s_\fg(x)\bigl)&=\Phi\bigl(s_\fg(x_N+x_\fa)\bigl)\\
        &=\Phi\bigl(s_\fg(x_N)+s_\fg(x_\fa)+[x_N,x_\fa]_\fg\bigl)\\
&=s_{\fh\oplus\fh^*}\bigl(\Phi(x_N)\bigl)+s_{\fh\oplus\fh^*}\bigl(\Phi(x_\fa)\bigl)+\bigl[\Phi(x_N),\Phi(x_\fa)\bigl]_{\fh\oplus\fh^*}\\
&=s_{\fh\oplus\fh^*}\bigl(\Phi(x_N+x_\fa)\bigl)=s_{\fh\oplus\fh^*}\bigl(\Phi(x)\bigl).
\end{align*}
This implies that $\Phi$ commutes with the squarings. Therefore, $\Phi$ is an isomorphism of Lie superalgebras. Thus, $\fg$ is isomorphic to the $T^*$-extension $\fh\oplus\fh^*$.

For the periplectic case, let us define 
\begin{equation}\label{cocyclestandardperi}
\beta:=\Pi\circ i_\kappa\circ\tilde{\alpha} \text{ and } \theta:=\Pi\circ i_\kappa\circ\tilde{\gamma}.
\end{equation}
We have $(\beta,\theta)\in XZ^2_\ev(\fh,\Pi(\fh^*))$. The isomorphism is then given by
\begin{equation}
    \Psi=p_\fh\oplus(\Pi\circ i_\kappa):\fg\rightarrow \fh\oplus\Pi(\fh^*).
\end{equation}
\end{proof}
The tuple $(\fh,\nabla,\alpha,\gamma)$ (resp. $(\fh,\nabla,\beta,,\theta)$) constructed above is called the \textit{even} (resp. \textit{odd}) extension tuple associated to $(\fg,\w,\fa,N)$, where $\w$ is ortho-orthogonal (resp. periplectic).

\subsection{Functoriality of the correspondence}\label{functor} Following \cite{BC,BM}, we will study the functorial properties of the constructions. The followings Lemmas are straightforward adaptation of \cite[Lemmas 6.2.1 and 6.2.2]{BM}.

\sssbegin{Lemma}
    Let $(\fg,\w)$ \textup{(}resp. $(\fg',\w')$\textup{)} be an ortho-orthogonal or periplectic quasi-Frobenius Lie superalgebra. Let $\Phi:(\fg,\w)\rightarrow(\fg',\w')$ be an isomorphism of Lie superalgebras preserving the forms, and let $\fa$ be an abelian ideal of $\fg$ and let $(\fh=\fg/\fa,\nabla)$ \textup{(}resp. $(\fh'=\fg'/\fa',\nabla')$\textup{)} be the associated flat Lie superalgebra, where $\nabla$  \textup{(}resp. $\nabla'$\textup{)} is given by Equation \eqref{inducednabla} and $\fa':=\Phi(\fa)$. Finally, let $\Phi_\fh:\fh\rightarrow\fh'$ be the induced map on the quotient spaces. Then, $\nabla'=(\Phi_\fh)_*\nabla$, where $(\Phi_\fh)_*\nabla:=\nabla\circ\Phi_\fh^{-1}$ is the push-forward of $\nabla$.
\end{Lemma}

An\textit{ isomorphism of strongly polarized quasi-Frobenius Lie superalgebras} $\Phi:(\fg,\w,\fa,N)\rightarrow(\fg',\w',\fa',N')$ is an isomorphism of Lie superalgebras $\Phi:(\fg,\w)\rightarrow(\fg',\w')$ preserving the forms such that $\Phi(\fa)=\fa'$ and $\Phi(N)=N'$.

\sssbegin{Lemma}\label{pushforwardcocycle}
    Let $\Phi:(\fg,\w,\fa,N)\rightarrow(\fg',\w',\fa',N')$ be an isomorphism of strongly polarized quasi-Frobenius Lie superalgebras. 
    \begin{enumerate}[(i)]
    \item \underline{In the case where $\w$ is ortho-orthogonal}, denote $(\alpha,\gamma)$ \textup{(}respectively $(\alpha',\gamma')$) the 2-cocycle associated with $(\fg,\w,\fa,N)$ (respectively $(\fg',\w',\fa',N')$\textup{)} as given in Equation \eqref{cocyclestrandardortho}. Then,
    \begin{equation}
        (\alpha',\gamma')=\bigl((\Phi_\fh)_*\alpha,(\Phi_\fh)_*\gamma\bigl).
    \end{equation}

    \item \underline{In the case where $\w$ is periplectic}, denote $(\beta,\theta)$ \textup{(}respectively $(\beta',\theta')$\textup{)} the 2-cocycle associated with $(\fg,\w,\fa,N)$ (respectively $(\fg',\w',\fa',N')$) as given in Equation \eqref{cocyclestandardperi}. Then,
    \begin{equation}
        (\beta',\theta')=\bigl((\Phi_\fh)_*\beta,(\Phi_\fh)_*\theta\bigl).
        \end{equation}
    \end{enumerate}
\end{Lemma}

\sssbegin{Corollary}
    There is a one-to-one correspondence between isomorphism classes of strongly polarized Lie superalgebras and isomorphism classes of extension tuples.
\end{Corollary}

\subsection{Change of strong polarization} Let $(\fg,\w,\fa,N)$ and $(\fg,\w,\fa,N')$ be two strong polarizations of the same Lagrangian extension. 
\begin{enumerate}[(i)]
\item \underline{In the case where $\w$ is ortho-orthogonal}, denote $(\alpha,\gamma)$ (resp. $(\alpha',\gamma')$) the extension 2-cocycle associated with $(\fg,\w,\fa,N)$ (resp. $(\fg,\w,\fa,N')$).
\item \underline{In the case where $\w$ is periplectic}, denote $(\beta,\theta)$ (resp. $(\beta',\theta')$) the extension 2-cocycle associated with $(\fg,\w,\fa,N)$ (resp. $(\fg,\w,\fa,N')$).
\end{enumerate}
The aim of this Section is to show that $(\alpha,\gamma)$ (resp. $(\beta',\theta')$) is cohomologous to $(\alpha',\gamma')$ (resp. $(\beta',\theta')$). The following Lemma is analogous to \cite[Lemma 6.3.1]{BM}.

\sssbegin{Lemma}\label{lem:choice} Let $(\fg,\w,\fa,N)$ and $(\fg,\w,\fa,N')$ be two strong polarizations of the same Lagrangian extension with corresponding extension cocycles given respectively by $(\alpha,\gamma)$ and $(\alpha',\gamma')$ in the ortho-orthogonal case ($(\beta,\theta)$ and $(\beta',\theta')$ in the periplectic case).
    \begin{enumerate}[(i)]
\item \underline{In the case where $\w$ is ortho-orthogonal}, there exists $\sigma\in \Hom(\fh,\fh^*)$ satisfying $$\sigma(x)(y)=\sigma(y)(x),~\forall x,y\in \fh,$$ such that $(\alpha',\gamma')=(\alpha,\gamma)+\fd^1\sigma$.
\item \underline{In the case where $\w$ is periplectic}, there exists $\mu\in \Hom(\fh,\Pi(\fh^*))$ satisfying $$\mu(x)(y)=\mu(y)(x),~\forall x,y\in \fh,$$ such that $(\beta',\theta')=(\beta,\theta)+\fd^1\mu$.
\end{enumerate}
\end{Lemma}

\begin{proof}
Let us prove the Lemma in the case where $\w$ is ortho-orthogonal. Denote by $p_\fa:\fa\oplus N\rightarrow\fa$ and  $p_\fa':\fa\oplus N'\rightarrow\fa$ the canonical projections. Let us consider the even map $\tau:=p_\fa-p_\fa'$. Notice that $\fa\subset\Ker(\tau)$. Moreover, denote by $p_N:\fa\oplus N\rightarrow N$ and  $p_N':\fa\oplus N'\rightarrow N'$ the projections on the second component. Let $x\in\fg$. Since $x=p_\fa(x)+p_N(x)=p_\fa'(x)+p_N'(x)$, we have $p_N'=p_N+\tau$. Similarly to \cite[Lemma 6.3.1]{BM}, we have
    \begin{equation}
        \w\bigl(\tau(n),m\bigl)=\w\bigl(p_N'(n),p_N(m)\bigl)= \w\bigl(n,\tau(m)\bigl),~\forall n,m\in N.
    \end{equation}
Recall that, using the same notations as Theorem \ref{conversethm}, we have $(\alpha,\gamma)=(i_{\w}\circ\tilde{\alpha},i_{\w}\circ\tilde{\gamma}),$
where      
    \begin{equation} \tilde{\alpha}:=p_\fa\bigl([p_N(\tilde{u}),p_N(\tilde{v})]_\fg\bigl),~\forall u,v\in\fh\text{ and }\tilde{\gamma}(u):=p_\fa\bigl(s_\fg(p_N(\tilde{u}))\bigl),~\forall u\in \fh_\od.
    \end{equation}
    For $u\in\fh,$ let $\tilde{\tau}(u):=\tau(\tilde{u})$. First, notice that
    \begin{equation}
        s_\fg\bigl(p_N'(\tilde{u})\bigl)=s_\fg(p_N(\tilde{u}))+[\tilde{u},\tau(\tilde{u})],~\forall u\in\fh.
    \end{equation} Then, we have
    \begin{align*}
        \gamma'(u)&=\w\bigl(p_\fa'\circ s_\fg(p_N'(\tilde{u})),\cdot\bigl)\\
        &=\w\bigl(p_\fa'(s(p_N(\tilde{u}))+[\tilde{u},\tau(\tilde{u})]),\cdot\bigl)\\
        &=\gamma(u)+\w\bigl(p_\fa([\tilde{u},\tau(\tilde{u})])+\tau\circ s(p_N(\tilde{u})),\cdot\bigl)\\
        &=\gamma(u)+\rho(u)\circ i_\w\circ\tau(\tilde{u})+i_\w\circ\tau(\widetilde{s(u)})\\
        &=\gamma(u)+\rho(u)\circ i_\w\circ\tilde{\tau}(u)+i_\w\circ\tilde{\tau}(s(u))\\
        &=\gamma(u)+\delta^1(i_\w\circ\tilde{\tau})(u).
    \end{align*} Thus, we take $\sigma:=i_\w\circ\tilde{\tau}$. Moreover, it has been shown in \cite[Lemma 6.3.1]{BM} that $\alpha'=\alpha+d^1_{\ce}\sigma$. Therefore, we have  $(\alpha',\gamma')=(\alpha,\gamma)+\fd^1\sigma$.
\end{proof}

\subsection{Equivalence of Lagrangian extensions} In this section, we introduce the Lagrangian cohomology in characteristic 2 and show that it captures Lagrangian extensions of Lie superalgebras equipped with a flat torsion-free connection. A Lagrangian extension of a Lie superalgebra equipped with a flat torsion-free connection $(\fh,\nabla)$ is a quasi-Frobenius Lie superalgebra $(\fg,\w)$ together with a short exact sequence of Lie superalgebras
\begin{equation}
        0\longrightarrow\fa{\longrightarrow}\fg{\longrightarrow}\fh\longrightarrow 0,
\end{equation} where the image of $\fa$ is a Lagrangian ideal of $\fg$. An \textit{isomorphism of Lagrangian extensions of $\fh$} is a Lie isomorphism $\Phi:(\fg,\w)\rightarrow(\fg',\w')$ satisfying $\w(x,y)=\w'\bigl(\Phi(x),\Phi(y)\bigl),~\forall x,y\in \fg$ such that the following diagram commutes:
\begin{equation}\label{diag-iso}
    \begin{tikzcd}
0 \arrow[r] & \fa \arrow[r] \arrow[d, "\Phi|_{\fa}"'] & \fg \arrow[r] \arrow[d, "\Phi" description] & \fh \arrow[r] \arrow[d, "="] & 0 \\
0 \arrow[r] & \fa' \arrow[r]                          & \fg' \arrow[r]                              & \fh \arrow[r]                & 0
\end{tikzcd}
\end{equation}
Note that such an isomorphism must satisfy $\Phi|_\fh=\id$.

\subsubsection{Lagrangian cohomology} Let $(\fh,\nabla)$ a Lie superalgebra equipped with a flat torsion-free connection $\nabla$. We construct the second cohomology space of the Lagrangian cohomology in characteristic 2 (see \cite{BC,BM} for the characteristic zero case). This cohomology describe all ortho-orthogonal $T^*$-extensions and periplectic $\Pi T^*$-extensions of $(\fh,\nabla)$, see Theorem \ref{correspondence-thm}.
 We define the first and second Lagrangian cochain spaces as
 \begin{align*}
    XC^1_L(\fh,\fh^*)&:=\bigl\{\psi\in XC^1(\fh,\fh^*),~\psi(x)(y)=\psi(y)(x),~\forall x,y\in \fh  \bigl\},\\
    XC^2_L(\fh,\fh^*)&:=\bigl\{(\alpha,\gamma)\in XC^2(\fh,\fh^*), (\alpha,\gamma) \text{ satisfies } \eqref{condalphagamma1} \text{ and } \eqref{condalphagamma2} \bigl\}.\\
    XC^1_L(\fh,\Pi(\fh^*))&:=\bigl\{\psi\in XC^1(\fh,\Pi(\fh^*)),~\psi(x)(y)=\psi(y)(x),~\forall x,y\in \fh  \bigl\},\\
    XC^2_L(\fh,\Pi(\fh^*))&:=\bigl\{(\beta,\theta)\in XC^2(\fh,\Pi(\fh^*)), (\beta,\theta) \text{ satisfies } \eqref{condbetatheta1} \text{ and } \eqref{condbetatheta2} \bigl\}.
 \end{align*}
In that case, the $\fh$-module structure on $\fh^*$ (resp. $\Pi(\fh^*)$) is given by the map $\rho$ as defined in Eq. \eqref{module1} (resp. the map $\chi$ defined in Eq. \eqref{module2}).

\sssbegin{Lemma}
     Let $(\fh,\nabla)$ be a Lie superalgebra equipped with a flat torsion-free connection $\nabla$. Then $\fd^1(XC^1_L(\fh,\fh^*))\subseteq XC^2_L(\fh,\fh^*)$ and $\fd^1(XC^1_L(\fh,\Pi(\fh^*)))\subseteq XC^2_L(\fh,\Pi(\fh^*))$. Moreover, we have $\fd^2\circ\fd^1=(0,0)$ in both cases.
\end{Lemma}
\begin{proof}
    Let $\psi\in XC^1_L(\fh,\fh^*)$, $x\in\fh_\od$ and $y\in\fh.$ We have
    \begin{align*}
        \delta^1\psi(x)+d_{\ce}^1\psi(x,y)(x)=&~\psi(x)\circ\nabla_x(y)+\psi(s_\fh(x))(y)+\psi(y)\circ\nabla_x(x)\\&~~~+\psi(x)\circ\nabla_y(x)+\psi([x,y])(x)\\
        =&~\underline{\psi(x)\circ\nabla_x(y)}+\psi(s_\fh(x))(y)+\psi(y)\circ s_\fh(x)\\&~~~+\underline{\psi(x)\circ\nabla_y(x)}+\underline{\psi(x)([x,y])}\\   
        =&~0,
    \end{align*} where the underlined terms cancel because $\nabla$ is torsion-free and the others because of the definition of the space $XC^1_L(\fh,\fh^*)$. Therefore, the pair $(d_{\ce}^1\psi,\delta^1\psi)$ satisfies Condition \eqref{condalphagamma2}. Similarly, we show that $d_{\ce}^1$ satisfy  Condition \eqref{condalphagamma1}. Thus, $\fd^1\psi=(d_{\ce}^1\psi,\delta^1\psi)\in XC^2_L(\fh,\fh^*)$. Let us show that $\fd^2\circ\fd^1=(0,0)$. 
    \begin{align*}
        \delta^2\bigl(d^1_{\ce}\psi,\delta^1\psi\bigl)=&~\rho(x)\circ d^1_{\ce}\psi(x,y)+\rho(y)\circ\delta^1\psi(x)+d^1_{\ce}\psi(s_\fh(x),y)+d^1_{\ce}\psi([x,y],y)\\
        =&~\psi(y)\circ\nabla_x^2+\psi(x)\circ\nabla_y\circ\nabla_x+\psi([x,y])\circ\nabla_x\\
        &~~~+\psi(x)\circ\nabla_x\circ\nabla_y+\psi(s_\fh(x))\circ\nabla_y\\
        &~~~+\psi(y)\circ\nabla_{s_\fh(x)}+\psi(s_\fh(x))\circ\nabla_y+\psi([s_\fh(x),y])\\
        &~~~+\psi(x)\circ\nabla_{[x,y]}+\psi([x,y])\circ\nabla_x+\psi\bigl( [[x,y],x]\bigl)=0,
    \end{align*}
    since the connection $\nabla$ is flat. Similarly, we have $d^2_{\ce}\circ d^1_{\ce}\psi=0.$
\end{proof} 
 We denote by 
    \begin{align*}
        XZ^2_L(\fh,\fh^*)&:=\bigl\{(\alpha,\gamma)\in XC^2_L(\fh,\fh^*),~\fd^2(\alpha,\gamma)=(0,0)\bigl\}\\
        XZ^2_L(\fh,\Pi(\fh^*))&:=\bigl\{(\beta,\theta)\in XC^2_L(\fh,\Pi(\fh^*)),~\fd^2(\beta,\theta)=(0,0)\bigl\}.
    \end{align*} Thus, we define
    \begin{align*}
    XH^2_L(\fh,\fh^*)&:=XZ^2_L(\fh,\fh^*)/\fd^1\bigl(XC^1_L(\fh,\fh^*)\bigl);\\
    XH^2_L(\fh,\Pi(\fh^*))&:=XZ^2_L(\fh,\Pi(\fh^*))/\fd^1\bigl(XC^1_L(\fh,\Pi(\fh^*))\bigl),
    \end{align*}
    the second cohomology groups of the Lagrangian cohomology.

\sssbegin{Theorem}\label{correspondence-thm}
   Let $(\fh,\nabla)$ be a Lie superalgebra equipped with a flat torsion-free connection $\nabla$. There is a one-to-one correspondence \textup{(}up to isomorphism in the sense of \eqref{diag-iso}\textup{)} between 
    \begin{enumerate}[$(i)$]
        \item $T^*$-extensions of $\fh$ and $XH^2_L(\fh,\fh^*)_\ev$, in the ortho-orthogonal case;
        \item $\Pi T^*$-extensions of $\fh$ and $XH^2_L(\fh,\Pi(\fh^*))_\ev$, in the periplectic case.
   \end{enumerate}
   Moreover, two Lagrangian extensions of $(\fh,\nabla)$ are isomorphic \textup{(}in the sense of \eqref{diag-iso}\textup{)} if and only if they correspond to the same extension 2-cocycle in $XH^2_L(\fh,\fh^*)_\ev$ or $XH^2_L(\fh,\Pi(\fh^*))_\ev$.
\end{Theorem}

\begin{proof}
We prove the Theorem in the ortho-orthogonal case. The proof follows \cite[Theorem 6.4.3]{BM}. Let $(\fg,\w,\fa)$ be a Lagrangian extension of $(\fh,\nabla)$. Consider a strong polarization $(\fg,\w,\fa,N)$. By Theorem \ref{sensdirect}, $(\fg,\w,\fa,N)$ is a $T^*$-extension of $(\fh,\nabla)$ by means of a 2-cocycle $(\alpha,\gamma)\in XH_L^2(\fh,\fh^*)_\ev$. By Lemma \ref{lem:choice}, the 2-cocycle $(\alpha,\gamma)$ does not depend on the choice of the strong polarization $N$. Let $\Phi:(\fg,\w,\fa)\rightarrow(\fg',\w',\fa')$ be an isomorphism of Lagrangian extensions (in the sense of \eqref{diag-iso}) and denote by $(\alpha',\gamma')\in XH_L^2(\fh,\fh^*)_\ev$ the 2-cocycle corresponding to $(\fg',\w',\fa')$. Then,  for any choice of strong polarization $N$, the map $\Phi:(\fg,\w,\fa,N)\rightarrow(\fg',\w',\fa',\Phi(N))$ is an isomorphism of strongly polarized ortho-orthogonal Lie superalgebras (see Section \ref{functor}). By Lemma \ref{pushforwardcocycle}, we have $\alpha'=\Phi_*\alpha$ and $\gamma'=\Phi_*\gamma$. Therefore, isomorphic extensions correspond to cohomologous extension cocycles. 

Let $(\fg,\w,\fa)$ and $(\fg',\w',\fa')$ be two $T^*$-extensions of $(\fh,\nabla)$ sharing the same extension cocycle. By Theorem \ref{sensdirect}, it is enough to show that any two strongly polarized ortho-orthogonal Lie superalgebras produce isomorphic extensions if $(\alpha,\gamma)=(\alpha',\gamma')\in XH_L^2(\fh,\fh^*)_\ev$, that is, $(\alpha',\gamma')=(\alpha,\gamma)+\fd^1\sigma$ (see Lemma \ref{lem:choice}). The map
$$\Phi:(\fg,\w,\fa)\rightarrow (\fg',\w',\fa'),\quad (u,\xi)\mapsto (u,\xi+\sigma(u))$$ is an isomorphism, since we have
\begin{align*}
    s_\fg\bigl(\Phi(u+\xi)\bigl)&=\bigl(s_\fh(u),\gamma(u)+\rho(u)\xi+\rho(u)\sigma(u)\bigl)\\
    &=\bigl(s_\fh(u),\gamma'(u)+\rho(u)\sigma(u)+\sigma(s_\fh(u))\bigl)=\Phi\bigl(s_{\fg'}(u+\xi)\bigl).
\end{align*}
\end{proof}

\subsection{Examples of Lagrangian extensions}\label{fewewamples} Hereafter we list a few examples.

\subsubsection{The Lie superalgebra $\bf{L_{1|1}^1}$} Consider the solvable Lie superalgebra $\fh:=\bf{L_{1|1}^1}$ of superdimension $(1|1)$ given in the basis $(e|f)$ by the bracket $[e,f]=f$ and the squaring $s=0$. Let  $\eps\in\K$ be a parameter. We define a flat torsion-free connection $\nabla^{\eps}$ on $\fh$ by
\begin{equation}
    \nabla^\eps_e(e)=(1+\eps)e,\quad \nabla^\eps_e(f)=\eps f,\quad \nabla^\eps_f(e)=(1+\eps)f,\quad \nabla^\eps_f(f)=0.
\end{equation}
The $\fh$-module structure on $\fh^*$ is given by
\begin{equation}
    e\cdot e^*=(1+\eps)e^*,\quad f\cdot e^*=0,\quad e\cdot f^*=\eps f^*,\quad f\cdot f^*=(1+\eps)e^*.
\end{equation}
Let $\lambda_1,\lambda_2,\lambda_3\in\K$ and consider the 2-cochains $(\alpha_1,\gamma_1), (\alpha_2,\gamma_2), (0,\gamma_3)$, with 
$$\alpha_1=e^*\otimes e^*\wedge f^*,\gamma_1(f)=\lambda_1 e^*;~\alpha_2=f^*\otimes e^*\wedge f^*,\gamma_2(f)=\lambda_2 e^*;~\gamma_3(f)=\lambda_3 e^*.$$

\textbf{Claim 1.} The cochain $(\alpha_1,\gamma_1)$ is a $2$-cocycle if and only if  $\eps=1$ or ($\eps\neq 1$ and $\lambda_1=0$).
\noindent\textit{Proof of Claim 1.} The second cocycle condition (see Eq. \eqref{2cocycle}) is given by
\begin{equation}
f\cdot\alpha_1(f,z)+z\cdot\gamma_1(f)=\alpha_1\bigl(s(f),z\bigl)+\alpha_1\bigl([z,f],f\bigl),~\forall z\in\fh.
\end{equation}
By plugging $z=e$ in the above equation, we obtain $\lambda_1(1+\eps)e^*=0.$ The case $z=f$ reduces to $0=0$. Similarly, we obtain the follwing statements.

\textbf{Claim 2.} The cochain $(\alpha_2,\gamma_2)$ is a $2$-cocycle if and only if $\eps=1$ or ($\eps\neq 1$ and $\lambda_2=1$).

\textbf{Claim 3.} The cochain $(0,\gamma_3)$ is a $2$-cocycle if and only if $\eps=1$ or $\lambda_3=0$.

\textbf{Claim 4.} The cohomology group $XH^2(\fh,\fh^*)$ is not trivial if and only if $\eps=1$. In that case, a basis of $XH^2(\fh,\fh^*)$ is given by $\{(\alpha_2,0),(0,\gamma_3^1)\}$, where $\alpha_2=f^*\otimes e^*\wedge f^*$ and $\gamma_3^1(f)=e^*$.\\
\noindent\textit{Proof of Claim 4.} First, notice that $\gamma_1,\gamma_2$ and $\gamma_3$ are cohomologous to $\gamma_3^1$.\\
\underline{The case where $\eps\neq1$}. We have
\begin{align*}
d^1_{\ce}(f^*\otimes f^*)(e,f)&=(f^*\otimes f^*)([e,f])+e\cdot(f^*\otimes f^*)(f)+f\cdot(f^*\otimes f^*)(e)=(1+\eps)f^*;\\
\delta^1(f^*\otimes f^*)(f)&=f\cdot(f^*\otimes f^*)(f)+(f^*\otimes f^*)(s(f))=(1+\eps)e^*.
\end{align*}
Therefore, $\alpha_2=d^1_{\ce}\bigl((1+\eps)^{-1}(f^*\otimes f^*)\bigl)$ and $\gamma_2=\delta^1\bigl((1+\eps)^{-1}(f^*\otimes f^*)\bigl)$. Thus, $(\alpha_2,\gamma_2)$ is a coboundary. Similarly, one can show that
$$(\alpha_1,0)=\Bigl(d^1_{\ce}\bigl(f^*\otimes e^*+e^*\otimes f^*\bigl),\delta^1\bigl(f^*\otimes e^*+e^*\otimes f^*\bigl)\Bigl).$$ Therefore, $XH_L^2(\fh,\fh^*)=0.$ \\
\underline{The case where $\eps=1$}. We have
\begin{align*}
d^1_{\ce}(e^*\otimes f^*)(e,f)&=(e^*\otimes f^*)([e,f])+e\cdot(e^*\otimes f^*)(f)+f\cdot(e^*\otimes f^*)(e)=e^*;\\
\delta^1(e^*\otimes f^*)(f)&=(e^*\otimes f^*)(s(f))=0.
\end{align*} Therefore, $(\alpha_1,0)=\bigl(d^1_{\ce}(e^*\otimes f^*+f^*\otimes e^*),\delta^1(e^*\otimes f^*+f^*\otimes e^*)\bigl)$ is a coboundary. Therefore, $XH^2(\fh,\fh^*)$ has dimension $2$ and is spanned by $\{(\alpha_2,0),(0,\gamma_3^1)\}$. In that case, we compute the Lagrangian cohomology. The pairs $(\alpha_2,0)$ and $(0,\gamma_3^1)$ does not satisfy Eqs. \eqref{condalphagamma1} and \eqref{condalphagamma2}. However, we have
$$\gamma_3^1(f)(e)+\alpha_2(f,e)(f)=e^*(e)+f^*(f)=0 \text{ and } \gamma_3^1(f)(f)+\alpha_2(f,f)(f)=0.$$ Therefore, the pair $(\alpha_2,\gamma_3^1)$ is Lagrangian. It follows that the space $XH^2_L(\fh,\fh^*)$ is one-dimensional and spanned by $(\alpha_2,\gamma_3^1)$.

\subsubsection{The Lie superalgebra $\fba(1)$} Consider the nilpotent Lie superalgebra $\fh:=\fba(1)$ of superdimension $(1|2)$ given in the basis $(e_1|e_2,e_3)$ by the bracket $[e_1,e_2]=e_3$ and the squaring $s(e_2)=s(e_3)=0.$ Let $\eps\in\K$ be a parameter. An  example of a flat torsion-free connection $\nabla$ on $\fh$ has the form
\begin{equation}
    \nabla_{e_1}(e_1)=\eps e_1,\quad \nabla_{e_1}(e_2)= e_3.
\end{equation}
The $\fh$-module structure $\rho$ on $\fh^*$ is then given by
\begin{equation}
    \rho(e_1)(e_1^*)=\eps e_1^*,\quad \rho(e_1)(e_3^*)= e_2^*.
\end{equation} Consider the following $2$-cochains.

\begin{equation*}\label{2cocyclee}
\begin{array}{llll}
(\alpha_1=e_1^*\otimes e_1^*\wedge e_2^*,\gamma_1=0);&&(\alpha_2=e_1^*\otimes e_1^*\wedge e_3^*,\gamma_2=0);&\\[2mm] (\alpha_3=0,\gamma_3(e_2)=e_1^*);&&  
(\alpha_4=0,\gamma_4(e_3)=e_1^*);\\[2mm] 
(\alpha_5=e_2^*\otimes e_1^*\wedge e_2^*,\gamma_5=0);& &(\alpha_6=e_2^*\otimes e_1^*\wedge e_3^*,\gamma_6=0);&  \\[2mm] 
(\alpha_7=0,\gamma_7(e_2)=e_2^*);&&
(\alpha_8=0,\gamma_8(e_3)=e_2^*);&   \\[2mm] 

 (\alpha_9=e_3^*\otimes e_1^*\wedge e_2^*,\gamma_9=0  );& &(\alpha_{10}=e_3^*\otimes e_1^*\wedge e_3^*,\gamma_{10}=0);&\\[2mm] 
 (\alpha_{11}=e_2^*\otimes e_2^*\wedge e_3^*,\gamma_{11}(e_2)=e_3^*).  \\[2mm] 
 \end{array}
\end{equation*}

\underline{The case where $\eps=0$.} In that case, the Lagrangian $2$-cocycles space is given by
$$ XZ^2_L(\fh,\fh^*)=\Span\bigl\{(\alpha_1,0);~ (\alpha_2,0);~(\alpha_5,\gamma_3);~(\alpha_{10},\gamma_4);~(\alpha_6+\alpha_9,0);~(\alpha_{11},\gamma_{11})     \bigl\}.$$
Moreover, we have $\fd^1(e_1^*\otimes e_3+e_3^*\otimes e_1)=(\alpha_1,0)$ and $\fd^1(e_3^*\otimes e_3^*)=(\alpha_6+\alpha_9,0)$. Therefore, $(\alpha_1,0)$ and $(\alpha_6+\alpha_9)$ are Lagrangian coboundaries and
$$ XH^2_L(\fh,\fh^*)=\Span\bigl\{(\alpha_2,0);~(\alpha_5,\gamma_3);~(\alpha_{10},\gamma_4);~(\alpha_{11},\gamma_{11})\bigl\}.$$

\underline{The case where $\eps=1$.} In that case, the Lagrangian $2$-cocycles space is given by
$$ XZ^2_L(\fh,\fh^*)=\Span\bigl\{(\alpha_1,0);~ (\alpha_2,0);~(\alpha_6+\widetilde{\alpha_9},0);~(\alpha_{11},\gamma_{11})     \bigl\}.$$
Moreover, we have $\fd^1(e_1^*\otimes e_2+e_2^*\otimes e_1)=(\alpha_1,0)$, $\fd^1(e_1^*\otimes e_3+e_3^*\otimes e_1+e_1^*\otimes e_2+e_2^*\otimes e_1)=(\alpha_2,0)$ and $\fd^1(e_3^*\otimes e_3^*)=(\alpha_6+\alpha_9,0)$. Therefore, $(\alpha_1,0)$, $(\alpha_2,0)$ and $(\alpha_6+\alpha_9,0)$ are Lagrangian coboundaries and
$$ XH^2_L(\fh,\fh^*)=\Span\bigl\{(\alpha_{11},\gamma_{11})      \bigl\}.$$

\subsubsection{The Lie superalgebra $D^1$} Consider the solvable Lie superalgebra $D^1$ of superdimension $ \sdim=(3|1)$ given in the basis $(e_1,e_2,e_3|e_4)$ (even $|$ odd) by the brackets $[e_2,e_3]=e_1,~ [e_2,e_4]=e_4$ and the squaring $s\equiv 0$.
The general form of a closed bilinear form on $D^1$ is given by
$$\w=c_{1,2}e_1^*\wedge e_2^*+c_{1,3}e_1^*\wedge e_3^*+c_{2,3}e_2^*\wedge e_3^*+c_{2,4}e_2^*\wedge e_4^*,~c_{i,j}\in\K.  $$
The form $\w$ is non-degenerate if and only if $c_{1,3}c_{2,4}\neq 0$. Therefore, there is no homogeneous closed bilinear form on $D^1$. 

\subsubsection{The Lie superalgebra $(2A_{1,1}+2A)^2$} Consider the nilpotent  Lie superalgebra $(2A_{1,1}+2A)^2$ of superdimension $ \sdim=(2|2)$ given in the basis $(e_1,e_2|e_3,e_4)$ by the bracket $[e_3,e_4]=e_1$ and the squaring $s(e_3)=e_1,~s(e_4)=e_2,~s(e_3+e_4)=e_2$. The general form of a closed bilinear form on $(2A_{1,1}+2A)^2$ is given by
$$\w=\lambda e_2^*\wedge e_4^*,~\lambda\in\K.$$ This form $\w$ is always degenerate.
    

\subsubsection{The Lie superalgebra $(C^1_1+A)$} Consider the solvable Lie superalgebra $(C^1_1+A)$ of superdimension $ \sdim=(2|2)$ given in the basis $(e_1,e_2|e_3,e_4)$ by the brackets $[e_1,e_2]=e_2,~[e_1,e_3]=e_3,~[e_3,e_4]=e_2$ and the squaring $s(e_3)=s(e_4)=0,~s(e_3+e_4)=e_2$.  Consider the closed ortho-orthogonal form
$$\w=e_1^*\wedge e_2^*+e_3^*\wedge e_4^*.$$ Let us denote $f:=e_3+e_4$. Then, the superalgebra $((C^1_1+A),\w)$ is a $T^*$-extension of the $(1|1)$-dimensional superalgebra $\fh$ spanned by $(e|f)$ with the bracket $[e,f]_\fh=f$ and the squaring $s_{\fh}\equiv 0$.  The flat torsion-free connection on $\fh$ is given by
$$\nabla_{e_1}(e_1)=e_1,~\nabla_{e_1}(f)=f,~\nabla_{f}(e_1)=f,~\nabla_{f}(f)=0$$ and the $2$-cocycle is given by $(\alpha,\gamma),$ where $\alpha=f^*\otimes e_1^*\wedge f^*$ and $\gamma(f)=e_1^*$. 

\section{Classification of Lagrangian extensions in dimension $4$}\label{sectionclassif4}

This section aims at classifying all $4$-dimensional Lie superalgebras that are either $T^*$-extensions or $\Pi T^*$-extensions of $2$-dimensional Lie superalgebras. In this Section, the base field $\K$ is assumed to be algebraically closed.

\subsection{Left-symmetric structures and cohomology} Recall that there is a one-to-one correspondence between flat torsion-free connections on a Lie superalgebras $(\fh,[\cdot,\cdot],s)$ and left symmetric structures on $\fh$ that are compatible with the Lie superalgebra structure (in the sense of Proposition \ref{LSSA-Lie}). In the tables given in Appendix \ref{appendixcoho}, we list nonisomorphic left-symmetric products on the Lie superalgebras $\fh$ given in Section \ref{classifLSSA2}, as well as the basis for nontrivial Lagrangian $2$-cocycles with values in $\fh^*$, where the action of $\fh$ on $\fh^*$ (respectively $\Pi(\fh^*$)) is given by Equation \eqref{module1} (respectively Equation \eqref{module2}).

\newpage

\subsection{Building the extensions} In this section, we list all the 4-dimensional Lie superalgebras obtained as Lagrangian extensions of the 2-dimensional Lie superalgebras

Whenever a coefficient $\lambda$ appears, it should be understood as $\forall \lambda\in\K$.

\subsubsection{Superalgebras of $\sdim=(2|2)$ obtained as $\Pi T^*$-extensions of $\bf{L^1_{0|2}}$} The brackets are given in the basis $(\Pi(e_1^*),\Pi(e_2^*)|e_1,e_2)$.\\
\scriptsize{}
$
\begin{array}{lllll}
    \bf{L^l_{2|2}}:& [e_1,e_2]=\Pi(e_1^*);&s(e_1)=\Pi(e_2^*).&&\\
    \bf{L^m_{2|2}}:& [e_1,e_2]=\Pi(e_2^*);&s(e_2)=\Pi(e_1^*).&&\\
    \bf{L^{l+m}_{2|2}}:& [e_1,e_2]=\Pi(e_1^*)+\Pi(e_2^*);&s(e_1)=\Pi(e_2^*);&s(e_2)=\Pi(e_1^*);&s(e_1+e_2)=\Pi(e_1^*)+\Pi(e_2^*).\\
 \end{array}
$
\normalsize{}

\subsubsection{Superalgebras of $\sdim=(2|2)$ obtained as $T^*$-extensions of $\bf{L^1_{1|1}}$} The brackets are given in the basis $(e_1,e_1^*|e_2,e_2^*)$.\\

\scriptsize{}
$
\begin{array}{lllll}
    \bf{L^a_{2|2}}:& [e_1,e_2]=e_2;&[e_1,e_2^*]=e_2^*.&& \\
    \bf{\widetilde{L^a_{2|2}}}:& [e_1,e_2]=e_2+e_2^*;&[e_1,e_2^*]=e_2^*;&s(e_2)=e_1^*.& \\
    \bf{L^b_{2|2}}:& [e_1,e_2]=e_2+e_2^*;&s(e_2)=e_1^*.&& \\
    \bf{L^c_{2|2}}(\eps):& [e_1,e_2]=e_2;&[e_1,e_1^*]=\eps e_1^*;&[e_1,e_2^*]=e_2^*,&\eps\neq 0,1. \\
    \bf{L^d_{2|2}}(\eps):& [e_1,e_2]=e_2;&[e_1,e_1^*]=\eps e_1^*;&[e_1,e_2^*]=(1+\eps)e_2^*;&\eps\neq0.\\
     & s(e_2+\lambda e_2^*)=\lambda\eps e_1^*;&[e_2,e_2^*]=\eps e_1^*.&& \\

 \end{array}
$

\normalsize{}
\subsubsection{Superalgebras of $\sdim=(2|2)$ obtained as $\Pi T^*$-extensions of $\bf{L^1_{1|1}}$}

The brackets are given in the basis $(e_1,\Pi(e_2^*)|e_2,\Pi(e_1^*))$.\\

\scriptsize{}
$
\begin{array}{lllll}
    \bf{L^n_{2|2}}:& [e_1,e_2]=e_2;&[e_1,\Pi(e_2^*)]=\Pi(e_2^*).&& \\
        \bf{L^o_{2|2}}(\eps):& [e_1,e_2]=e_2;&[e_1,\Pi(e_1^*)]=\eps\Pi(e_1^*);&[e_1,\Pi(e_2^*)]=\Pi(e_2^*),&\eps\neq 0,1. \\
        \bf{L^p_{2|2}}:& [e_1,e_2]=e_2+\Pi(e_1^*);&[e_1,\Pi(e_1^*)]=\Pi(e_1^*);&[e_1,\Pi(e_2^*)]=\Pi(e_2^*).& \\
        \bf{L^q_{2|2}}(\eps):& [e_1,e_2]=e_2;&[e_1,\Pi(e_1^*)]=\eps\Pi(e_1^*);&[e_1,\Pi(e_2^*)]=(1+\eps)\Pi(e_2^*);& \eps\neq0\\&[e_2,\Pi(e_2^*)]=\eps\Pi(e_1^*).&&\\
 \end{array}
$

\normalsize{}
\subsubsection{Superalgebras of $\sdim=(2|2)$ obtained as $T^*$-extensions of $\bf{L^2_{1|1}}$}

The brackets are given in the basis $(e_1,e_1^*|e_2,e_2^*)$.\\

\scriptsize{}
$
\begin{array}{lllll}
    \bf{L^e_{2|2}}:& [e_2,e_1^*]=e_2^*;&s(e_2)=e_1.&&\\
    \bf{L^f_{2|2}}:& [e_1,e_1^*]=e_1^*;&[e_1,e_2^*]=e_2^*;&[e_2,e_1^*]=e_2^*;&\\
    &s(e_2)=e_1;& s(e_2+\lambda e_2^*)=e_1+\lambda e_1^*;&[e_2,e_2^*]=e_1^*.&\\
    \bf{L^g_{2|2}}:& s(e_2)=e_1.&&&\\

 \end{array}
$
\normalsize{}
\subsubsection{Superalgebras of $\sdim=(2|2)$ obtained as $\Pi T^*$-extensions of $\bf{L^2_{1|1}}$}

The brackets are given in the basis $(e_1,\Pi(e_2^*)|e_2,\Pi(e_1^*))$.\\

\scriptsize{}
$
\begin{array}{lllll}
  \bf{L^r_{2|2}}:& s(e_2)=e_1;&s(e_2+\lambda\Pi(e_1^*))=e_1+\lambda\Pi(e_2^*);&[e_2,\Pi(e_1^*)]=\Pi(e_2^*).&\\
    \bf{L^s_{2|2}}:&[e_1,e_2]=\Pi(e_1^*);& [e_1,\Pi(e_1^*)]=\Pi(e_1^*);&[e_1,\Pi(e_2^*)]=\Pi(e_2^*);&\\
    &[e_2,\Pi(e_1^*)]=\Pi(e_2^*);&[e_2,\Pi(e_2^*)]=\Pi(e_1^*);&s(e_2)=e_1+
    \Pi(e_2^*).&\\& s(e_2+\lambda\Pi(e_1^*))=e_1+(1+\lambda)\Pi(e_2^*).&

 \end{array}
$
\normalsize{}
\subsubsection{Superalgebras of $\sdim=(2|2)$ obtained as $T^*$-extensions of $\bf{L^3_{1|1}}$}

The brackets are given in the basis $(e_1,e_1^*|e_2,e_2^*)$.\\
\scriptsize{}
$
\begin{array}{lllll}
    \bf{L^h_{2|2}}:& [e_1,e_1^*]=e_1^*;&&&\\
    \bf{L^i_{2|2}}:& [e_1,e_1^*]=e_1^*;&[e_1,e_2^*]=e_2^*;&[e_2,e_2^*]=e_1^*;&s(e_2+\lambda e_2^*)=\lambda e_1^*.\\
    \bf{L^j_{2|2}}:& \text{abelian.}&&&\\
    \bf{L^k_{2|2}}:& [e_1,e_2]=e_2^*;&s(e_2)=e_1^*.&&\\
 \end{array}
$

\normalsize{}
\newpage
\subsubsection{Superalgebras of $\sdim=(2|2)$ obtained as $\Pi T^*$-extensions of $\bf{L^3_{1|1}}$}

The brackets are given in the basis $(e_1,\Pi(e_2^*)|e_2,\Pi(e_1^*))$.\\

\scriptsize{}
$
\begin{array}{lllll}
  \bf{L^t_{2|2}}:& [e_1,\Pi(e_1^*)]=\Pi(e_1^*).&&\\
   \bf{L^u_{2|2}}:& [e_1,\Pi(e_1^*)]=\Pi(e_1^*);&[e_1,\Pi(e_2^*)]=\Pi(e_2^*);&[e_2,\Pi(e_2^*)]=\Pi(e_1^*).\\
    \bf{L^v_{2|2}}:&[e_1,e_2]=\Pi(e_1^*);& [e_1,\Pi(e_1^*)]=\Pi(e_1^*);&[e_1,\Pi(e_2^*)]=\Pi(e_2^*);&\\
    & [e_2,\Pi(e_2^*)]=\Pi(e_1^*).&&\\
    \bf{L^w_{2|2}}:& [e_1,e_2]=\Pi(e_1^*).&&\\

 \end{array}
$

\normalsize{}
\subsubsection{Superalgebras of $\sdim=(2|2)$ obtained as $\Pi T^*$-extensions of $\bf{L^1_{2|0}}$}

The brackets are given in the basis $(e_1,e_2|\Pi(e_1^*),\Pi(e_2^*))$.\\

\scriptsize{}
$
\begin{array}{lllll}
    
    \bf{L^x_{2|2}}(\eps):& [e_1,e_2]=e_2;&[e_1,\Pi(e_1^*)]=\eps \Pi(e_1^*);&\eps\neq 0,1\\&[e_1,\Pi(e_2^*)]=\Pi(e_2^*).&\\~\\
   
    \bf{L^y_{2|2}}(\eps): &[e_1,e_2]=e_2;&[e_1,\Pi(e_1^*)]=\Pi(e_1^*);&\eps\neq 0\\&[e_1,\Pi(e_2^*)]=\eps \Pi(e_1^*)+\Pi(e_2^*),&\\~\\
   
    \bf{L^z_{2|2}}:& [e_1,e_2]=e_2;&[e_1,\Pi(e_2^*)]=
    \eps \Pi(e_1^*)+\Pi(e_2^*).\\~\\
   
    \bf{L^{aa}_{2|2}}: &[e_1,e_2]=e_2;&[e_1,\Pi(e_1^*)]=
     \Pi(e_1^*);\\&[e_1,\Pi(e_2^*)]= \Pi(e_2^*).\\~\\
     
    \bf{L^{bb}_{2|2}}(\eps):& [e_1,e_2]=e_2;&[e_1,\Pi(e_1^*)]=
     (1+\eps)\Pi(e_1^*);\\&[e_1,\Pi(e_2^*)]=\eps \Pi(e_2^*);&[e_2,\Pi(e_1^*)]=\Pi(e_2^*);\\&[e_2,\Pi(e_2^*)]=(1+\eps)\Pi(e_1^*).&&\\~\\
     
    \bf{L^{cc}_{2|2}}(\eps):& [e_1,e_2]=e_2;&[e_1,\Pi(e_1^*)]=
     (1+\eps)\Pi(e_1^*);&\eps\neq1\\&[e_1,\Pi(e_2^*)]=\eps \Pi(e_2^*);&[e_2,\Pi(e_2^*)]=(1+\eps)\Pi(e_1^*).\\~\\
     
    \bf{L^{dd}_{2|2}}:& [e_1,e_2]=e_2;&[e_1,\Pi(e_1^*)]=
     \Pi(e_1^*)+\Pi(e_2^*);\\&[e_2,\Pi(e_1^*)]=\Pi(e_1^*);&[e_2,\Pi(e_2^*)]=\Pi(e_1^*)+\Pi(e_2^*).\\~\\
     
    \bf{L^{ee}_{2|2}}(\eps):& [e_1,e_2]=e_2;&[e_1,\Pi(e_1^*)]=
     \Pi(e_1^*)+\Pi(e_2^*);&\eps\neq0\\&[e_2,\Pi(e_1^*)]=\Pi(e_1^*)+\eps \Pi(e_2^*);&[e_2,\Pi(e_2^*)]=\Pi(e_1^*)+\Pi(e_2^*).\\~\\
     
    \bf{L^{ff}_{2|2}}:& [e_1,e_2]=e_2;&[e_1,\Pi(e_1^*)]=
     \Pi(e_1^*)+\Pi(e_2^*);\\&[e_1,\Pi(e_2^*)]= \Pi(e_1^*);&[e_2,\Pi(e_1^*)]=\Pi(e_1^*);\\&[e_2,\Pi(e_2^*)]=\Pi(e_1^*)+\Pi(e_2^*).&&\\~\\
     
     
    \bf{L^{hh}_{2|2}}(\eps):& [e_1,e_2]=e_2;&[e_1,\Pi(e_1^*)]=
     (1+\eps)\Pi(e_1^*)+\Pi(e_2^*);&\eps\neq0,1;\\&[e_2,\Pi(e_1^*)]=\Pi(e_1^*)+(1+\eps)^{-1}\Pi(e_2^*);&[e_2,\Pi(e_2^*)]=(1+\eps)\Pi(e_1^*)+\Pi(e_2^*);&&\\
     &[e_1,\Pi(e_2^*)]= \eps(1+\eps)\Pi(e_1^*)+\eps \Pi(e_2^*).&&&\\~\\

     \bf{L^{ii}_{2|2}}:& [e_1,e_2]=e_2.&&&\\

 \end{array}
$

\normalsize{}

\subsubsection{Superalgebras of $\sdim=(2|2)$ obtained as $\Pi T^*$-extensions of $\bf{L^2_{2|0}}$} The brackets are given in the basis $(e_1,e_2|\Pi(e_1^*),\Pi(e_2^*))$.\\
\scriptsize{}
$
\begin{array}{llllll}

    \bf{L^{ll}_{2|2}}:& [e_1,\Pi(e_1^*)]=\Pi(e_1^*).&&&\\
    \bf{L^{nn}_{2|2}}:& [e_1,\Pi(e_2^*)]=\Pi(e_1^*).&&&\\

    \bf{L^{pp}_{2|2}}:& [e_1,\Pi(e_1^*)]=\Pi(e_1^*);&[e_2,\Pi(e_2^*)]=\Pi(e_2^*).&&\\
   
    \bf{L^{rr}_{2|2}}:& [e_1,\Pi(e_1^*)]=\Pi(e_1^*);&[e_1,\Pi(e_2^*)]=\Pi(e_2^*);&[e_2,\Pi(e_2^*)]=\Pi(e_1^*).\\~\\
   
    \bf{L^{tt}_{2|2}}:& [e_1,\Pi(e_1^*)]=\Pi(e_1^*);&[e_1,\Pi(e_2^*)]=\Pi(e_2^*);&[e_2,\Pi(e_1^*)]=\Pi(e_2^*);\\&[e_2,\Pi(e_2^*)]=\Pi(e_1^*)+\Pi(e_2^*).\\
   
    
    \bf{L^{xx}_{2|2}}:&\text{ abelian.}
 \end{array}
$
\normalsize{}
\newpage
\subsubsection{Superalgebras of $\sdim=(0|4)$ obtained as $ T^*$-extensions of $\bf{L^1_{0|2}}$}The brackets are given in the basis $(0|e_1,e_2,e_1^*,e_2^*)$.\\

\scriptsize{}
$
\begin{array}{lllll}
    \bf{L^a_{0|4}}:& \text{abelian}.&&\\

 \end{array}
$



\subsubsection{Superalgebras of $\sdim=(4|0)$ obtained as $ T^*$-extensions of $\bf{L^1_{2|0}}$}
The brackets are given in the basis $(e_1,e_1^*,e_2,e_2^*|0)$.\\

\scriptsize{}
$
\begin{array}{lllll}
    \bf{L^a_{4|0}}(\eps):& [e_1,e_2]=e_2+e_2^*;&[e_1,e_1^*]=\eps e_1^*;&[e_1,e_2^*]=e_2^*,& \eps\neq 0,1.\\
    
    \bf{L^b_{4|0}}(\eps):& [e_1,e_2]=e_2;&[e_1,e_1^*]=\eps e_1^*;&[e_1,e_2^*]=e_2^*,&\eps\neq 0,1.\\
    
    \bf{L^c_{4|0}}(\eps):& [e_1,e_2]=e_2+e_2^*;&[e_1,e_1^*]=e_1^*;&[e_1,e_2^*]=\eps e_1^*+e_2^*,&\eps\neq 0,1.\\
    
    \bf{L^d_{4|0}}(\eps): &[e_1,e_2]=e_2;&[e_1,e_1^*]=e_1^*;&[e_1,e_2^*]=\eps e_1^*+e_2^*,&\eps\neq 0.\\
    
    \bf{L^e_{4|0}}(\eps): &[e_1,e_2]=e_2+e_2^*;&[e_1,e_2^*]=
    \eps e_1^*+e_2^*.&\\
    
    \bf{L^f_{4|0}}(\eps):& [e_1,e_2]=e_2;&[e_1,e_2^*]=
    \eps e_1^*+e_2^*.\\
    
    \bf{L^g_{4|0}}:& [e_1,e_2]=e_2+e_1^*;&[e_1,e_1^*]=
     e_1^*;&[e_1,e_2^*]= e_2^*.\\
     
    \bf{L^h_{4|0}}: &[e_1,e_2]=e_2+e_2^*;&[e_1,e_1^*]=
     e_1^*;&[e_1,e_2^*]= e_2^*.\\

         \bf{L^{g+h}_{4|0}}: &[e_1,e_2]=e_2+e_1^*+e_2^*;&[e_1,e_1^*]=
     e_1^*;&[e_1,e_2^*]= e_2^*.\\
     
    \bf{L^i_{4|0}}: &[e_1,e_2]=e_2;&[e_1,e_1^*]=
     e_1^*;&[e_1,e_2^*]= e_2^*.\\~\\
     
    \bf{L^j_{4|0}}(\eps):& [e_1,e_2]=e_2;&[e_1,e_1^*]=
     (1+\eps)e_1^*;&[e_1,e_2^*]=\eps e_2^*;\\&[e_2,e_1^*]=e_2^*;&[e_2,e_2^*]=(1+\eps)e_1^*.&&\\~\\
     
    \bf{L^k_{4|0}}(\eps):& [e_1,e_2]=e_2;&[e_1,e_1^*]=
     (1+\eps)e_1^*;&[e_1,e_2^*]=\eps e_2^*;\\&[e_2,e_2^*]=(1+\eps)e_1^*.&&\\~\\
     
    \bf{L^l_{4|0}}:& [e_1,e_2]=e_2;&[e_1,e_1^*]=
     e_1^*+e_2^*;&[e_2,e_1^*]=e_1^*;\\&[e_2,e_2^*]=e_1^*+e_2^*.&&\\~\\
     
    \bf{L^m_{4|0}}(\eps):& [e_1,e_2]=e_2;&[e_1,e_1^*]=
     e_1^*+e_2^*;&[e_2,e_1^*]=e_1^*+\eps e_2^*,&\eps\neq0;\\&[e_2,e_2^*]=e_1^*+e_2^*.&&\\~\\
     
    \bf{L^n_{4|0}}:& [e_1,e_2]=e_2;&[e_1,e_1^*]=
     e_1^*+e_2^*;&[e_1,e_2^*]= e_1^*;\\&[e_2,e_1^*]=e_1^*;&[e_2,e_2^*]=e_1^*+e_2^*.&&\\~\\
     
    
    \bf{L^p_{4|0}}(\eps):& [e_1,e_2]=e_2;&[e_1,e_1^*]=
     (1+\eps)e_1^*+e_2^*;&&\eps\neq0,1;\\&[e_2,e_1^*]=e_1^*+(1+\eps)^{-1}e_2^*;&[e_2,e_2^*]=(1+\eps)e_1^*+e_2^*;&&\\
     &[e_1,e_2^*]= \eps(1+\eps)e_1^*+\eps e_2^*.&&&\\~\\
     
    \bf{L^q_{4|0}}:& [e_1,e_2]=e_2+e_1^*.&&&\\
    
    \bf{L^r_{4|0}}:& [e_1,e_2]=e_2+e_2^*.&&&\\

\bf{L^{q+r}_{4|0}}:& [e_1,e_2]=e_2+e_1^*+e_2^*.&&&\\
 \end{array}
$
\normalsize{}
\subsubsection{Superalgebras of $\sdim=(4|0)$ obtained as $ T^*$-extensions of $\bf{L^2_{2|0}}$} The brackets are given in the basis $(e_1,e_1^*,e_2,e_2^*|0)$.\\
\scriptsize{}
$
\begin{array}{llllll}
    \bf{L^s_{4|0}}:& [e_1,e_2]=e_1^*.&&&\\
    \bf{L^t_{4|0}}:& [e_1,e_2]=e_2^*.&&&\\

    \bf{L^{s+t}_{4|0}}:& [e_1,e_2]=e_1^*+e_2^*.&&&\\
    \bf{L^u_{4|0}}:& [e_1,e_1^*]=e_1^*.&&&\\
    \bf{L^v_{4|0}}:&[e_1,e_2]=e_2^*;& [e_1,e_1^*]=e_1^*.&&\\
    \bf{L^w_{4|0}}:& [e_1,e_2^*]=e_1^*.&&&\\
    \bf{L^x_{4|0}}:&[e_1,e_2]=e_2^*;& [e_1,e_2^*]=e_1^*.&&\\
    \bf{L^y_{4|0}}:& [e_1,e_1^*]=e_1^*;&[e_2,e_2^*]=e_2^*.&&\\
    \bf{L^z_{4|0}}:&[e_1,e_2]=e_2^*;& [e_1,e_1^*]=e_1^*;&[e_2,e_2^*]=e_2^*.&\\
    \bf{L^{aa}_{4|0}}:& [e_1,e_1^*]=e_1^*;&[e_1,e_2^*]=e_2^*;&[e_2,e_2^*]=e_1^*.\\
    \bf{L^{bb}_{4|0}}:&[e_1,e_2]=e_1^*;& [e_1,e_1^*]=e_1^*;&[e_1,e_2^*]=e_2^*;&[e_2,e_2^*]=e_1^*.\\
    \bf{L^{cc}_{4|0}}:& [e_1,e_1^*]=e_1^*;&[e_1,e_2^*]=e_2^*;&[e_2,e_1^*]=e_2^*;&[e_2,e_2^*]=e_1^*+e_2^*.\\
    \bf{L^{dd}_{4|0}}:&[e_1,e_2]=e_1^*;& [e_1,e_1^*]=e_1^*;&[e_1,e_2^*]=e_2^*;&[e_2,e_1^*]=e_2^*;&[e_2,e_2^*]=e_1^*+e_2^*.\\
    \bf{L^{gg}_{4|0}}:&\text{ abelian.}
 \end{array}
$

\normalsize{}

\subsection{Isomorphisms vs. symplectomorphisms} The following Propositions gives isomorphisms and symplectomorphisms between the superalgebras of the above lists. Our main tool is the SuperLie package (see \cite{G}) of the computer algebra system Mathematica and the invariants tables given in Appendix \ref{invariants}.  Note that those symplectomorphisms are not isomorphisms in the sense of \eqref{diag-iso} since their restriction to $\fh$ is never the identity.

\sssbegin{Proposition}[Extensions of superdimension $(2|2)$]

  We have the following Lie superalgebras isomorphisms:
$$\bf{L^{tt}_{2|2}}\cong\bf{L^{pp}_{2|2}};~\bf{L^{f}_{2|2}}\cong\bf{L^{s}_{2|2}};~ \bf{L^{q}_{2|2}}(\eps)\cong\bf{L^{cc}_{2|2}}(\nu)\iff \eps\nu=1;~
        \bf{L^{bb}_{2|2}}(\eps)\cong\bf{L^{dd}_{2|2}}(\nu)\cong\bf{L^{ff}_{2|2}}~(\eps\neq1); $$
        $$\bf{L^{u}_{2|2}}\cong\bf{L^{v}_{2|2}}\cong\bf{L^{ee}_{2|2}}(1)\cong\bf{L^{hh}_{2|2}}(\eps);\bf{L^{l}_{2|2}}\cong\bf{L^{r}_{2|2}};~\bf{L^{i}_{2|2}}\cong\bf{L^{d}_{2|2}}(1);$$
        $$\bf{L^{x}_{2|2}}(\eps)\cong\bf{L^{o}_{2|2}}(\nu)\iff \eps=\nu;~\bf{L^{y}_{2|2}}(\eps)\cong\bf{L^{p}_{2|2}}.$$
     We have the following symplectic Lie superalgebras symplectomorphisms:
$$ \bf{L^{q}_{2|2}}(1)\cong\bf{L^{rr}_{2|2}};~ \black\bf{L^{ll}_{2|2}}\cong\bf{L^{t}_{2|2}};~\bf{L^{nn}_{2|2}}\cong\bf{L^{w}_{2|2}}; $$
$$ \bf{L^{cc}_{0}}\cong\bf{L^{u}_{2|2}};~\bf{L^{ii}_{2|2}}\cong\bf{L^{h}_{2|2}};~\bf{L^{l}_{2|2}}\cong\bf{L^{m}_{2|2}}\cong\bf{L^{l+m}_{2|2}}. $$

\end{Proposition}

\sssbegin{Proposition}[Extensions of superdimension $(4|0)$]
    We have the following Lie superalgebras isomorphisms:
        $$\bf{L^{d}_{4|0}}(\eps)\cong\bf{L^{g}_{4|0}}\cong\bf{L^{h}_{4|0}}\cong\bf{L^{g+h}_{4|0}};~\bf{L^{j}_{4|0}}(\eps)\cong\bf{L^{n}_{4|0}},\eps\neq 1;~\bf{L^{s}_{4|0}}\cong\bf{L^{s+t}_{4|0}}; $$
        $$\bf{L^{j}_{4|0}}(1)\cong\bf{L^{k}_{4|0}}(0)\cong\bf{L^{p}_{4|0}}(\eps)\cong\bf{L^{bb}_{4|0}};~\bf{L^{q}_{4|0}}\cong\bf{L^{u}_{4|0}};~\bf{L^{y}_{4|0}}\cong\bf{L^{z}_{4|0}};~\bf{L^{y}_{4|0}}\cong\bf{L^{dd}_{4|0}}. $$
    We have the following symplectic Lie superalgebras symplectomorphisms:
    $$ \bf{L^{d}_{4|0}}(\eps)\cong\bf{L^{g}_{4|0}};~ \bf{L^{h}_{4|0}}\cong\bf{L^{g+h}_{4|0}};~\bf{L^{l}_{4|0}}\cong\bf{L^{m}_{4|0}(\eps)} \iff \eps\neq 1;~\bf{L^{f}_{4|0}}(\eps)\cong\bf{L^{k}_{4|0}}(1);$$
    $$\bf{L^{j}_{4|0}}(1)\cong\bf{L^{m}_{4|0}}(1);~ \bf{L^{q}_{4|0}}\cong\bf{L^{r}_{4|0}}\cong\bf{L^{q+r}_{4|0}}; \bf{L^{s}_{4|0}}\cong\bf{L^{t}_{4|0}}\cong\bf{L^{w}_{4|0}};~\bf{L^{y}_{4|0}}\cong\bf{L^{cc}_{4|0}};~\bf{L^{z}_{4|0}}\cong\bf{L^{dd}_{4|0}}. $$

\end{Proposition}

\newpage

\section{Appendix: Post-Lie superalgebras in characteristic $2$}\label{postlie}

In this Appendix, we initiate the study of post-Lie superalgebras in characteristic $2$. Over a field of characteristic zero, Vallette introduced the post-Lie operad as the Koszul dual of the operad of commutative triassociative algebras (\cite{Va}). They were independently described in the context of numerical integration methods in \cite{MKW,MKL}. For a survey on those structures, see \cite{CEM}.

\subsection{Post-Lie superalgebras}\label{defpost1}
Let $(V,\{\cdot,\cdot\},s)$ be a Lie superalgebra in characteristic $2$. It is called \textit{post-Lie superalgebra} if it is endowed with a bilinear product $\ast:V\times V\rightarrow V$ satisfying
\begin{equation*}
\begin{array}{lrlll}
(i) &\{x,y\}\ast z& = & (y\ast x)\cdot z+y\ast(x\ast z)+(x\ast y)\ast z+x\ast(y\ast z),&\forall x,y,z\in V; \\[2mm]
(ii) &s(x)\ast y&=&x\ast(x\ast y)+(x\ast x)\ast y,&\forall x\in V_\od,~\forall y\in V;\\[2mm]
(iii) &x\ast\{y,z\}&=&\{x\ast y,z\}+\{y,x\ast z\},&\forall x,y,z\in V;\\[2mm]
(iv)&x\ast s(y)&=&\{x\ast y,y\},&\forall x\in V_\od,~\forall y\in V.\\[2mm]
 \end{array}
\end{equation*}

A post-Lie superalgebra will be denoted $(V,\{\cdot,\cdot\},s,\ast).$ Note that a post-Lie product over an abelian Lie superalgebra in nothing more that a left-symmetric product.\\

\noindent\textbf{Remarks.}
\begin{enumerate}
\item Axiom $(i)$ is equivalent to
$ \{x,y\}\ast z=\Asso(y,x,z)+\Asso(x,y,z),~\forall x,y,z\in V,  $ where $\Asso(u,v,w)=(u\ast v)\ast w+u\ast(v\ast w),~\forall u,v,w\in V.$

\item Axiom $(ii)$ is equivalent to $s(x)\ast y=\Asso(x,x,y),~\forall x\in V_\od,~\forall y\in V;$
\item Axioms $(iii)$ and $(iv)$ together are equivalent to saying that $\fl_x$ is a derivation for all $x\in V$ (in the sense of \eqref{derivation}), where $\fl_x:V\rightarrow V,~y\mapsto x\ast y$ is the left multiplication.
\end{enumerate}

\sssbegin{Proposition}\label{propassoc}
 Let $(V,\{\cdot,\cdot\},s,\ast)$ be a post-Lie superalgebra. Let us define a new bracket and a new squaring on $V$ as follows:
\begin{center}
    $\begin{cases}
        [x,y]:=x\ast y+y\ast x+\{x,y\},~\forall x,y\in V;\\
        \tilde{s}(x):=x\ast x+s(x),~\forall x\in V_\od.
     \end{cases}$
\end{center}
Then $(V,[\cdot,\cdot],\tilde{s})$ is a Lie superalgebra, called \emph{associated Lie superalgebra.}
 
\end{Proposition}
\begin{proof}
    We only prove the Jacobi identity involving the squaring. Let $x\in V_\od$ and $y\in V$. We have
    \begin{align*}
        [x,[x,y]]=&~[x,x\ast y]+[x,y\ast x]+[x,\{x,y\}]\\
                 =&~x\ast(x\ast y)+(x\ast y)\ast x+\{x,x\ast y\}\\
                 &+x\ast(y\ast x)+(y\ast x)\ast x+\{x,y\ast x\}\\
                 &+x\ast\{x,y\}+\{x,y\}\ast x+\{x,\{x,y\}\}.
                 \end{align*}
On the other hand, 
\begin{align*}
[\tilde{s}(x),y]=&~(x\ast x)\ast y+s(x)\ast y+y\ast(x\ast x)+y\ast s(x)+\{x\ast x,y\}+\{s(x),y\}\\
                =&~x\ast(x\ast y)+y\ast(x\ast x)+\{y\ast x,x\}+\{x\ast x,y\}+\{s(x),y\}.
    \end{align*}     
    Therefore,
    $$[x,[x,y]]+[\tilde{s}(x),y]=\{x\ast x,y\}+\{x,x\ast y\}+x\ast\{x,y\}=0.  \qed $$
    \noqed
\end{proof}

\sssbegin{Proposition}
 Let $(V,\{\cdot,\cdot\},s,\ast)$ be a post-Lie superalgebra and $(V,[\cdot,\cdot],\tilde{s})$ be its associated Lie superalgebra. Then
 \begin{align}
    [x,y]\ast z&=x\ast(y\ast z)+y\ast(x\ast z),~\forall x,y,z\in V;\\
    \tilde{s}(x)\ast y&=x\ast(x\ast y),~\forall x\in V_\od,~\forall y\in V.
 \end{align}
 
\end{Proposition}

\begin{proof}
We only prove the identity involving the squaring. Let $x\in V_\od$ and $y\in V$. We have
$$\tilde{s}(x)\ast y=\bigl(x\ast x+s(x)\bigl)\ast y=(x\ast x)\ast y+x\ast (x\ast y)+(x\ast x)\ast y=x\ast(x\ast y).\qed $$
\noqed 
\end{proof}

\sssbegin{Corollary}
   For all $x\in V$, the left multiplication $\fl_x$ is a representation of $(V,[\cdot,\cdot],\tilde{s})$.
\end{Corollary}

Following \cite{BD}, we will use the following notation. We will use 
$\fn:=(\fn,\{\cdot,\cdot\},s)$ for the Lie superalgebra $(V,\{\cdot,\cdot\},s)$ and $\fg:=(\fg,[\cdot,\cdot],\tilde{s})$ for the Lie superalgebra $(V,[\cdot,\cdot],\tilde{s})$ (see Prop \ref{propassoc}). As super vector spaces, we have $V=\fn=\fg$.

\subsubsection{An equivalent definition of post-Lie superalgebras}\label{defpost2} Let $(\fn,\fg)$ be a pair of Lie superalgebras sharing the same super vector space $V$. A post-Lie structure on $(\fn,\fg)$ is a bilinear map $\ast:V\times V\rightarrow V$ satisfying
\begin{equation}
\begin{array}{lrlll}
(i) &[x,y]-\{x,y\}& = & x\ast y+y\ast x,&\forall x,y\in V; \\[2mm]
(ii) &\tilde{s}(x)-s(x)& = & x\ast x,&\forall x\in V_\od \\[2mm]
(iii) &[x,y]\ast z& = & x\ast(y\ast z)+y\ast(x\ast z),&\forall x,y,z\in V; \\[2mm]
(iv) &\tilde{s}(x)\ast z& = & x\ast(x\ast z),&\forall x\in V_\od,~\forall z\in V; \\[2mm]
(v) &x\ast\{y,z\}& = & \{x\ast y,z\}+\{y,x\ast z\}&\forall x,y,z\in V; \\[2mm]
(vi) &x\ast s(y)& = & \{x\ast y,y\}&\forall x\in V,~\forall y\in V_\od. \\[2mm]
 \end{array}
\end{equation}
Definitions \ref{defpost1} and \ref{defpost2} are equivalent.

\sssbegin{Proposition}\label{semidirect}
    Let $(\fn,\fg)$ be a pair of Lie superalgebras sharing the same super vector space $V$. Let $\ast:V\times V\rightarrow V$ be a bilinear map. Then, the map $\ast$ is post-Lie with respect to $(\fn,\fg)$ (in the sense of Definition \ref{defpost2}) if and only if 
    \begin{enumerate}
        \item $\fl_x\in\Der(\fn)~\forall x\in V;$
        \item  The map 
                $\varphi:\fg  \rightarrow\Der(\fn)\rtimes\fn,~x\mapsto(\fl_x,x)$
        is a morphism of Lie superalgebras.
    \end{enumerate}
\end{Proposition}

\begin{proof}
The superspace $\Der(\fn)\rtimes\fn$ is a Lie superalgebra with the bracket 
\begin{equation}
    \bigl[(f,x),(g,y)\bigl]:=\bigl([f,g],f(y)+g(x)+\{x,y\}\bigl),~\forall f,g\in \Der(V),~\forall x,y\in V,
\end{equation} and the squaring 
\begin{equation}
    \fs(f,x):=\bigl(f\circ f,s(x)+f(x)\bigl),~\forall f\in \Der(V)_\od,~\forall x\in V_\od.
\end{equation}
We will prove the Jacobi Identity for the square $\fs$. Let $f\in\Der(\fn)_\od,~g\in\Der(\fn),~y\in\fn$ and $x\in \fn_\od$.
\begin{align*}
    \bigl[(f,x),\bigl[(f,x),(g,y)\bigl]\bigl]&=\bigl[(f,x),\bigl([f,g],f(y)+g(x)+\{x,y\}\bigl)\bigl]\\
    &=\Bigl([f,[f,g]],f\bigl(f(y)+g(x)+\{x,y\}\bigl)+[f,g](x)+\bigl\{x,f(y)+g(x)+\{x,y\}\bigl\}\Bigl)\\
    &=\Bigl([f\circ f,g],f\circ f(y)+g\bigl(s(x)\bigl)+g\circ f(x)+\{s(x),y\}+\{f(x),y\}\Bigl)\\
    &=\bigl[(f\circ f,s(x)+f(x)),(g,y)\bigl]\\
    &=\bigl[\fs(f,x),(g,y)\bigl].
\end{align*}

Suppose that $\ast$ defines a post-Lie structure on the pair $(\fn,\fg)$. Then, the left multiplication $\fl_x$ is a derivation of $\fn$ for all $x\in V$ by definition. Moreover, consider the map $$\varphi:\fg  \rightarrow\Der(\fn)\rtimes\fn,~x\mapsto(\fl_x,x).$$ We prove that $\varphi$ commutes with the squarings. First, notice that $\fl_x^2+\fl_{x\ast x}=\fl_{s(x)}~\forall x\in V_\od.$ Then, for all $x\in V_\od$, we have
$$\varphi\bigl(\tilde{s}(x)\bigl)=(\fl_{s(x)},s(x))+(\fl_{x\ast x},x\ast x)=\bigl(\fl_x^2,s(x)+\fl_x(x)\bigl)=\fs(\fl_x,x)=\fs\bigl(\varphi(x)\bigl).$$

Conversely, suppose that Conditions (1) and (2) of Proposition \ref{semidirect} are satisfied. Since $\fl_x$ is a derivation for all $x\in V$, Conditions (\textit{v}) and (\textit{vi}) of Definition \ref{defpost2} are satisfied. Furthermore, since $\varphi$ is a morphism of Lie superalgebras, we have
$$\bigl(\fl_{[x,y]},[x,y]\bigl)=\varphi\bigl([x,y]\bigl)=\bigl[\varphi(x),\varphi(y)\bigl]=\bigl(\fl_x\circ \fl_y+\fl_y\circ \fl_x,x\ast y+y \ast x+\{x,y\}\bigl)~\forall x,y\in V.$$ Therefore, we obtain Conditions (\textit{i}) and (\textit{iii}) of Definition \ref{defpost2}. Finally, since
$$\bigl(\fl_{\tilde{s}(x)},\tilde{s}(x)\bigl)=\varphi\bigl(\tilde{s}(x)\bigl)=\fs\bigl(\varphi(x)\bigl)=\fs(\fl_x,x)=\bigl(\fl_x^2,s(x)+\fl_x(x)\bigl),~\forall x\in V_\od,  $$ Conditions (\textit{ii}) and (\textit{iv}) of Definition \ref{defpost2} follow immediately. \qed

\noqed

\end{proof}

\sssbegin{Proposition}\label{postliecorrespondence}
    The tuple $(V,\{\cdot,\cdot\},s,*)$ is a post-Lie superalgebra if and only if its associated Lie superalgebra $(V,[\cdot,\cdot],\tilde{s})$ (see Proposition \ref{propassoc}) is equipped with a flat parallel connection.
\end{Proposition}
\begin{proof}
    Let $(V,\{\cdot,\cdot\},s,*)$ be a post-Lie superalgebra. Recall that the associated Lie superalgebra $\fg=(V,[\cdot,\cdot],\tilde{s})$ is given by the bracket $[x,y]=x\ast y+y*x+\{x,y\},~\forall x,y\in V$ and the squaring $\tilde{s}(x)=x*x+s(x),~\forall x\in V_\od$. We define a connection on $\fg$ by $\nabla_x(y):=x*y,~\forall x,y\in V$. Therefore, the torsion is given by
    \begin{align}
        T(x,y)&=x*y+y*x+[x,y]=\{x,y\},~\forall x,y\in V\\
        U(x)&=x*x+\tilde{s}(x)=s(x),~\forall x\in V_\od.
    \end{align} We compute the second part of the covariant derivative of the torsion. Let $z\in V$ and $x\in V_\od$.
    \begin{align*}
        (\nabla_zU)(x)&=z*U(x)+T(z*x,x)\\
        &=z*\bigl(x*x+\tilde{s}(x)\bigl)+(z*x)*x+x*(z*x)+[z*x,x]\\
        &=z*s(x)+\{z*x,x\}=0,~\text{see }\ref{defpost2}.
    \end{align*}
    Moreover, $S(x)(z)=\nabla_{\tilde{s}(x)}(z)+\nabla^2_x(z)=0,~\text{see }\ref{defpost2}.$ One can also show that $R=\nabla_z T=0$. Therefore, the connection $\nabla$ is flat and parallel.
    Reciprocally, consider a flat parallel connection $\nabla$ on $\fg=(V,[\cdot,\cdot],\tilde{s})$. Define a product on $V$ by $x*y:=\nabla_x(y),~\forall x,y\in V$. Using \ref{defpost2}, one can show that $(V,\{\cdot,\cdot\},s,*)$ is a post-Lie superalgebra.
\end{proof}
\newpage

\section{Appendix: Lagrangian cohomology of 2-dimensional Lie superalgebras}\label{appendixcoho}
For the definition of the Lie superalgebras of dimension 2, see Proposition \ref{classif2}.

\subsection{The case where $\fh=\bf{L^1_{0|2}}$} In that case, the only suitable connection is $\nabla=0$. We have
\small
\begin{align*}
    XH^2_L(\fh;\fh^*)&=\Span\bigl\{e_1^*\otimes e_1^*\wedge e_2^*,\gamma(e_1)=e_2^*);(e_2^*\otimes e_1^*\wedge e_2^*,\gamma(e_2)=e_1^*))\bigl\};\\
    XH^2_L(\fh;\Pi(\fh^*))&=\Span\bigl\{(\Pi(e_1^*)\otimes e_1^*\wedge e_2^*,\theta(e_1)=\Pi(e_2^*));(\Pi(e_2^*)\otimes e_1^*\wedge e_2^*,\theta(e_2)=\Pi(e_1^*))\bigl\}.
\end{align*}
\normalsize
\subsection{The case where $\sdim(\fh)=(1|1)$}
\small
\begin{center}
    \begin{tabular}{|c|c|c|c|}
    \hline
          $\fh$ & Left-symmetric product on $\fh$ & $(\alpha,\gamma)\in XH^2_L(\fh;\fh^*)$ &  $(\beta,\theta)\in XH^2_L(\fh;\Pi(\fh^*))$ \\\hline 
          & $0$ & $(e_2^*\otimes e_1^*\wedge e_2^*,\gamma(e_2)=e_1^*)$  & $(\Pi(e_2^*)\otimes e_1^*\wedge e_2^*,\theta(e_2)=\Pi(e_1^*))$ \\\cline{2-4} 
            & $e_1e_2=e_2$ &$(e_2^*\otimes e_1^*\wedge e_2^*,\gamma(e_2)=e_1^*)$  &$(\Pi(e_2^*)\otimes e_1^*\wedge e_2^*,\theta(e_2)=\Pi(e_1^*))$ \\\cline{2-4}  
           $\bf{L^1_{1|1}}$ & $e_1e_1=\eps e_1;~e_1e_2=e_2,~\eps\neq 0,1$ & $0$  &$0$  \\\cline{2-4} 
            & $e_1e_1= e_1;~e_1e_2=e_2$ & $(e_1^*\otimes e_1^*\wedge e_2^*,0)$  & $(\Pi(e_1^*)\otimes e_1^*\wedge e_2^*,0)$  \\\cline{2-4} 
          & $e_1e_1=\eps e_1;~e_1e_2=(1+\eps)e_2;$ & \multirow{2}{*}{$0$}  &\multirow{2}{*}{$0$}   \\
          & $e_2e_1=\eps e_2,~\eps\neq0.$ &  & \\\hline 
        \multirow{4}{*}{$\bf{L^2_{1|1}}$} & $0$ & $0$  & $0$ \\\cline{2-4}   
         & $e_2e_2=e_1$ & $0$  & $0$ \\\cline{2-4}  
        & $e_1e_1=e_1;~e_1e_2=e_2;$& \multirow{2}{*}{$(e_1^*\otimes e_1^*\wedge e_2^*,\gamma(e_2)=e_2^*)$}  & \multirow{2}{*}{$(\Pi(e_1^*)\otimes e_1^*\wedge e_2^*,\theta(e_2)=\Pi(e_2^*))$} \\ 
        & $e_2e_1=e_2;~e_2e_2=e_1;$&  &  \\\hline
        \multirow{5}{*}{$\bf{L^3_{1|1}}$} & \multirow{2}{*}{$0$}& $(e_1^*\otimes e_1^*\wedge e_2^*,0)$  & $(\Pi(e_1)^*\otimes e_1^*\wedge e_2^*,0)$ \\\cline{3-4}
         &  & $(e_2^*\otimes e_1^*\wedge e_2^*,\gamma(e_2)=e_1^*)$  & $(\Pi(e_2^*)\otimes e_1^*\wedge e_2^*,\theta(e_2)=\Pi(e_1^*))$ \\\cline{2-4}
         &$e_1e_1=e_1$&$0$&$0$\\\cline{2-4}
        & $ e_1e_1=e_1;~e_1e_2=e_2;$& \multirow{2}{*}{$(e_1^*\otimes e_1^*\wedge e_2^*,0)$}  & \multirow{2}{*}{$(\Pi(e_1^*)\otimes e_1^*\wedge e_2^*,0)$} \\
        & $e_2e_1=e_2$&   & \\\hline
           
    \end{tabular}
   \\~\\ The case where $\sdim(\fh)=(1|1)$.
\end{center}

\normalsize
\subsection{The case where $\fh=\bf{L^1_{2|0}}$}
\small\begin{center}
    \begin{tabular}{|c|c|c|}
    \hline
          Left-symmetric product & $\alpha\in XH^2_L(\fh;\fh^*)$ &  $\beta\in XH^2_L(\fh;\Pi(\fh^*))$\\\hline 
           
           $e_1e_1=\eps e_1;~e_1e_2=e_2,~\eps\neq 0,1$ & $e_2^*\otimes e_1^*\wedge e_2^*$   &$\Pi(e_2^*)\otimes e_1^*\wedge e_2^*$ \\\cline{1-3} 
           $e_1e_1=e_1+\eps e_2;~e_1e_2=e_2,~\eps\neq 0$ & $e_2^*\otimes e_1^*\wedge e_2^*$   &$\Pi(e_2^*)\otimes e_1^*\wedge e_2^*$ \\\cline{1-3} 
            $e_1e_1=\eps e_2;~e_1e_2=e_2$ & $e_2^*\otimes e_1^*\wedge e_2^*$   &$\Pi(e_2^*)\otimes e_1^*\wedge e_2^*$ \\\cline{1-3} 
            \multirow{2}{*}{$e_1e_1= e_1;~e_1e_2=e_2$} & $e_1^*\otimes e_1^*\wedge e_2^*$   &$\Pi(e_1^*)\otimes e_1^*\wedge e_2^*$ \\\cline{2-3}
             & $e_2^*\otimes e_1^*\wedge e_2^*$   &$\Pi(e_2^*)\otimes e_1^*\wedge e_2^*$ \\\cline{1-3} 
            $e_1e_1= (1+\eps)e_1;~e_1e_2=\eps e_2$; & \multirow{2}{*}{$0$} &\multirow{2}{*}{$0$} \\
            $e_2e_1=(1+\eps)e_2;~e_2e_2=e_1  $ &&\\\hline
            $e_1e_1= (1+\eps)e_1;~e_1e_2=\eps e_2$; & \multirow{2}{*}{$0$} &\multirow{2}{*}{$0$} \\
            $e_2e_1=(1+\eps)e_2;~e_2e_2=0  $ &&\\\hline
            $e_1e_1= e_1;~e_1e_2=e_1$; & \multirow{2}{*}{$0$} &\multirow{2}{*}{$0$} \\
            $e_2e_1=e_1+e_2;~e_2e_2=e_2  $ &&\\\hline
            $e_1e_1= e_1;~e_1e_2=e_1$; & \multirow{2}{*}{$0$} &\multirow{2}{*}{$0$} \\
            $e_2e_1=e_1+e_2;~e_2e_2=\eps e_1+e_2,~\eps\neq0  $ &&\\\hline
            $e_1e_1= e_1+e_2;~e_1e_2=e_1$; & \multirow{2}{*}{$0$} &\multirow{2}{*}{$0$} \\
            $e_2e_1=e_1+e_2;~e_2e_2=e_2 $ &&\\\hline
            $e_1e_1= (1+\eps)e_1+\eps(1+\eps)e_2;~e_1e_2=e_1+\eps e_2$; & \multirow{2}{*}{$0$} &\multirow{2}{*}{$0$} \\
            $e_2e_1=e_1+(1+\eps)e_2;~e_2e_2=(1+\eps)^{-1}e_1+e_2,~\eps\neq 0,1 $ &&\\\hline
            \multirow{2}{*}{$0$} & $e_1^*\otimes e_1^*\wedge e_2^*$   &$ \Pi(e_1^*)\otimes e_1^*\wedge e_2^*$\\\cline{2-3}  &$e_2^*\otimes e_1^*\wedge e_2^*$&$ \Pi(e_2^*)\otimes e_1^*\wedge e_2^*$\\\cline{1-3}
    \end{tabular}
   \\~\\ The case where $\fh=\bf{L^1_{2|0}}$.
\end{center}
\normalsize~\\
\subsection{The case where $\fh=\bf{L^2_{2|0}}$}
\small
\begin{center}
    \begin{tabular}{|c|c|c|}
    \hline
          Left-symmetric product & $\alpha\in XH^2_L(\fh;\fh^*)$ &  $\beta\in XH^2_L(\fh;\Pi(\fh^*))$\\\hline 
          
          \multirow{2}{*}{$0$} & $e_1^*\otimes e_1^*\wedge e_2^*$  &$\Pi(e_1^*)\otimes e_1^*\wedge e_2^*$\\\cline{2-3}
          &$e_2^*\otimes e_1^*\wedge e_2^*$&$\Pi(e_2^*)\otimes e_1^*\wedge e_2^*$\\\hline
          
           $e_1e_1=e_1$&$e_2^*\otimes e_1^*\wedge e_2^*$&$\Pi(e_2^*)\otimes e_1^*\wedge e_2^*$\\\hline
           
           $e_1e_1=e_2$&$e_2^*\otimes e_1^*\wedge e_2^*$&$\Pi(e_2^*)\otimes e_1^*\wedge e_2^*$\\\hline
           
            $e_1e_1=e_1;~e_2e_2=e_2$&$e_2^*\otimes e_1^*\wedge e_2^*$&$\Pi(e_2^*)\otimes e_1^*\wedge e_2^*$\\\hline
            
          $e_1e_1=e_1;~e_1e_2=e_2;~e_2e_1=e_2$&$e_1^*\otimes e_1^*\wedge e_2^*$&$\Pi(e_1^*)\otimes e_1^*\wedge e_2^*$\\\hline
          
    $e_1e_1=e_1;~e_1e_2=e_2;$&\multirow{2}{*}{$e_1^*\otimes e_1^*\wedge e_2^*$}&\multirow{2}{*}{$\Pi(e_1^*)\otimes e_1^*\wedge e_2^*$}\\$e_2e_1=e_2;~e_2e_2=e_1+e_2$&&\\\hline

    \end{tabular}
   \\~\\ The case where $\fh=\bf{L^2_{2|0}}$.
\end{center}
\normalsize

\section{Appendix: table of invariants}\label{invariants}
~\\
\subsection{Invariants for Lie superalgebras of superdimension $(2|2)$}
\small{
\begin{center}
    \begin{tabular}{|c|c|c|c|c|c|c|c|}
    \hline
          $\fg$ & $\fg^{(1)}$ & $\fz(\fg)$ & $XH^1(\fg;\K)$& $XH^2(\fg;\K)$ & $XH^3(\fg;\K)$&$XH^4(\fg;\K)$&$XH^5(\fg;\K)$ \\\hline 
           $\bf{L^a_{2|2}}$ & $(0|2)$ & $(1|0)$  & 2 & 4 & 6 &8 & 10\\\hline 
            $\widetilde{\bf{L^a_{2|2}}}$ & $(1|2)$ & $(1|0)$  & 1 & 2 & 2 &2 &2\\\hline 
           $\bf{L^b_{2|2}}$ &$(1|1)$ & $(1|1)$  & 2 & 2   & 2 & 2 & 2\\\hline 
           $\bf{L^c_{2|2}}(\eps)$ & $(1|2)$ & $0$ & 1  & 3 & 3 & 5 & 5\\
           \hline 
           $\bf{L^d_{2|2}}(\eps),\eps\neq 1$ &$(1|2)$  & $0$ & 1 & 2 & 2 & 2 & 2\\\hline 
           
           
           $\bf{L^d_{2|2}}, \eps=1$ &$(1|1)$  & $0$ & 2 & 3 & 3 & 3 & 3\\\hline 
           $\bf{L^e_{2|2}}$ &$(1|1)$ & $(1|1)$ & 2 & 3  & 3 & 3 & 3\\\hline 
           $\bf{L^f_{2|2}}$ & $(2|1)$ & $0$  & 1& 1  & 1 & 1 & 1\\\hline 
           $\bf{L^g_{2|2}}$ & $(1|0)$ & $(2|1)$ & 3 & 4  & 4 & 4 & 4\\\hline 
           $\bf{L^h_{2|2}}$ & $(1|0)$ & $(0|2)$ & 3 &  5 & 7 & 9 & 11\\\hline 
           $\bf{L^i_{2|2}}$ & $0$ & $(1|1)$  &  2&  3 & 3 & 3 & 3\\\hline 
           $\bf{L^j_{2|2}}$ &$0$ & $(2|2)$  & 4 & 8  & 12 & 16 & 20\\\hline 
           $\bf{L^k_{2|2}}$ & $(1|1)$  & $(1|1)$ & 3 & 4  & 4 & 4 & 4\\\hline 
    \end{tabular}
   \\~\\ Invariants for Lie superalgebras of $\sdim=(2|2)$ obtained as $T^*$-extensions of Lie superalgebras of $\sdim=(1|1)$.
\end{center}}
\normalsize{}

\small{
\begin{center}
    \begin{tabular}{|c|c|c|c|c|c|c|c|}
    \hline
          $\fg$ & $\fg^{(1)}$ & $\fz(\fg)$ & $XH^1(\fg;\K)$& $XH^2(\fg;\K)$ & $XH^3(\fg;\K)$&$XH^4(\fg;\K)$&$XH^5(\fg;\K)$ \\\hline 
           $\bf{L^l_{2|2}}$ & $(2|0)$ & $(2|0)$  &  2 & 2 & 2  &2 &2 \\\hline 
           $\bf{L^m_{2|2}}$ & $(2|0)$& $(2|0)$  & 2 & 2 &2  &2 &2 \\\hline 
           $\bf{L^{l+m}_{2|2}}$ & $(2|0)$& $(2|0)$  & 2 & 2 & 2 & 2& 2\\\hline 
          
    \end{tabular}
   \\~\\ Invariants for Lie superalgebras of $\sdim=(2|2)$ obtained as $\Pi T^*$-extensions of $\bf{L^1_{0|2}}$
\end{center}}
\normalsize{}

\small{
\begin{center}
    \begin{tabular}{|c|c|c|c|c|c|c|c|}
    \hline
          $\fg$ & $\fg^{(1)}$ & $\fz(\fg)$ & $XH^1(\fg;\K)$& $XH^2(\fg;\K)$ & $XH^3(\fg;\K)$&$XH^4(\fg;\K)$&$XH^5(\fg;\K)$ \\\hline 
           $\bf{L^n_{2|2}}$ &$(1|1)$& $(0|1)$  & 2 & 4 & 6 & 8 &10 \\\hline 
           $\bf{L^o_{2|2}}(\eps)$ & $(1|2)$& 0  & 1 & 3 & 3 &5 &5 \\\hline
           
            
           $\bf{L^p_{2|2}}$& $(1|2)$& $0$ & 1 & 3 & 3 & 5&5 \\\hline 
           $\bf{L^q_{2|2}}(\eps),\eps\neq1$ &$(1|2)$ & 0  & 1 & 2 & 3 & 4 & 5\\\hline 
             $\bf{L^q_{2|2}}(1),$ & $(0|2)$& 0  & 2 & 3 & 4 & 5 & 6\\\hline 
           $\bf{L^r_{2|2}}$ & $(2|0)$ & $(2|0)$ & 2 & 2 & 2 & 2& 2\\\hline 
           $\bf{L^s_{2|2}}$ &$(2|1)$&0   & 1 & 1 & 1 & 1&1 \\\hline 
           $\bf{L^t_{2|2}}$& $(0|1)$& $(1|1)$  & 3 & 5 & 7 & 9& 11\\\hline  
           $\bf{L^u_{2|2}}$ &$(1|1)$&0   & 2 & 4 & 6 & 8& 10\\\hline 
           $\bf{L^v_{2|2}}$ &$(1|1)$&0   &2  &4  &6  &8 &10 \\\hline 
               $\bf{L^w_{2|2}}$ &$(0|1)$& $(1|1)$  & 3 & 5 & 7 &9 &11 \\\hline 
          
    \end{tabular}
   \\~\\ Invariants for Lie superalgebras of $\sdim=(2|2)$ obtained as $\Pi T^*$-extensions of Lie superalgebras of $\sdim=(1|1)$.
\end{center}}
\normalsize{}

\small{
\begin{center}
    \begin{tabular}{|c|c|c|c|c|c|c|c|}
    \hline
          $\fg$ & $\fg^{(1)}$ & $\fz(\fg)$ & $XH^1(\fg;\K)$& $XH^2(\fg;\K)$ & $XH^3(\fg;\K)$&$XH^4(\fg;\K)$&$XH^5(\fg;\K)$ \\\hline 
           $\bf{L^x_{2|2}}(\eps)$ & $(1|2)$& 0  & 1 & 3 & 3 & 5 &5 \\\hline 
         $\bf{L^y_{2|2}}(\eps)$ & $(1|2)$& 0   & 1 & 3 & 3 & 5 &5 \\\hline 
          $\bf{L^z_{2|2}}$ & $(1|1)$ &$(0|1)$   & 2 &4 & 6 & 8 &12 \\\hline
          $\bf{L^{aa}_{2|2}}$ & $(1|2)$& 0   & 1 &5 & 5 & 9 &9 \\\hline
        $\bf{L^{bb}_{2|2}}(\eps), \eps\neq1$ & $(1|2)$& 0   & 1 &2 & 3 & 4 &5 \\\hline
        $\bf{L^{bb}_{2|2}}(1)$ &  $(1|1)$& 0   & 2&4 & 6 & 8 &12 \\\hline
        $\bf{L^{cc}_{2|2}}(\eps),\eps\neq0$ &  $(1|2)$& 0   & 1 &2 & 3 & 4 &5 \\\hline
        $\bf{L^{cc}_{2|2}}(0)$ &  $(1|1)$ & 0 & 2 & 4 & 6 & 8 & 10 \\\hline
        
        
        $\bf{L^{dd}_{2|2}}(\eps)$ & $(1|2)$ &0  & 1 &2 & 3 & 4 &5 \\\hline
        $\bf{L^{ee}_{2|2}}(\eps),\eps\neq1$ & $(1|2)$ &  0 & 1 &2 & 3 & 4 &5 \\\hline
        $\bf{L^{ee}_{2|2}}(1)$ &$(1|1)$  &   0& 2 &4 & 6 & 8 &10 \\\hline
        
        $\bf{L^{ff}_{2|2}}$ & $(1|2)$ &0   & 1 &2 & 3 & 4 &5 \\\hline
        
        
        $\bf{L^{hh}_{2|2}}(\eps)$ &  $(1|1)$ & 0   & 2 &4 & 6 & 8 &10 \\\hline
        
        $\bf{L^{ii}_{2|2}}$ & $(1|0)$ &  $(0|2)$ & 3 &5 & 7 & 9 &11 \\\hline
          
    \end{tabular}
   \\~\\ Invariants for Lie superalgebras of $\sdim=(2|2)$ obtained as $\Pi T^*$-extensions the Lie superalgebras $\bf{L_{2|0}^1}$.
\end{center}}
\normalsize{}

\small{
\begin{center}
    \begin{tabular}{|c|c|c|c|c|c|c|c|}
    \hline
          $\fg$ & $\fg^{(1)}$ & $\dim(\fz(\fg))$ & $XH^1(\fg;\K)$& $XH^2(\fg;\K)$ & $XH^3(\fg;\K)$&$XH^4(\fg;\K)$&$XH^5(\fg;\K)$\\\hline 
          
           $\bf{L^{ll}_{2|2}}$ & $(0|1)$ & $(1|1)$  &3  &5 & 7 & 9&11\\\hline
           
           $\bf{L^{nn}_{2|2}}$ &  $(0|1)$ & $(1|1)$  &3  &5 & 7 & 9&11\\\hline
           
           $\bf{L^{pp}_{2|2}}$ &  $(0|2)$ & $0$  &2  &3 & 4 & 5&6\\\hline
           
           $\bf{L^{rr}_{2|2}}$ & $(0|2)$ & 0 &2  &3 & 4 & 5&6\\\hline
           
           $\bf{L^{tt}_{2|2}}$ &  $(0|2)$ & 0  &2  &3 & 4 & 5&6\\\hline
           
           
          
    \end{tabular}
   \\~\\ Invariants for Lie superalgebras of $\sdim=(2|2)$ obtained as $\Pi T^*$-extensions of the Lie superalgebra $\bf{L^2_{2|0}}$.
\end{center}}
\newpage
\subsection{Invariants for Lie superalgebras of superdimension $(4|0)$}

\small{
\begin{center}
    \begin{tabular}{|c|c|c|c|c|c|c|}
    \hline
          $\fg$ & $\fg^{(1)}$ & $\dim(\fz(\fg))$ & $XH^1(\fg;\K)$& $XH^2(\fg;\K)$ & $XH^3(\fg;\K)$&$XH^4(\fg;\K)$ \\\hline 
           $\bf{L^a_{4|0}}(\eps)$ & $3$ & $0$  & 1 & 1& 1 &0 \\\hline 
           
           $\bf{L^b_{4|0}}(\eps)$ & $3$ &$0$   & 1 &1 & 1 & 0 \\\hline   
           $\bf{L^c_{4|0}}(\eps)$ & $3$ &$0$  & 1 & 1 & 1 &0 \\\hline
           
           $\bf{L^d_{4|0}}(\eps)$ & $3$ &$0$ & 1 & 2 & 2 & 0 \\\hline 
          
           $\bf{L^e_{4|0}}(\eps)$ & $2$ & $1$  & 2 & 2& 2 & 1\\\hline 
           $\bf{L^f_{4|0}}(\eps)$ & $2$ & $1$  & 2 & 2& 2 &1 \\\hline 
           $\bf{L^g_{4|0}}$ & $3$ & $0$  & 1 & 2 & 2  &0 \\\hline 
           $\bf{L^h_{4|0}}$ & $3$ & $0$  & 1 & 2& 2 &0  \\\hline 

    $\bf{L^{g+h}_{4|0}}$ & $3$ & $0$  & 1 & 2& 2 &0  \\\hline
           $\bf{L^i_{4|0}}$ & $3$ & $0$  & 1 &3 &3  &0  \\\hline 
           $\bf{L^j_{4|0}}(\eps),\eps\neq1$ & $3$ & $0$& 1 & 0& 1 & 1 \\\hline  $\bf{L^j_{4|0}}(1),$ & $2$ & $0$ & 2 & 2 & 2 & 1  \\\hline 
           $\bf{L^k_{4|0}}(\eps),\eps\neq0,1$ & $3$ & $0$  & 1 & 0 & 1 & 1\\\hline 
           $\bf{L^k_{4|0}}(0)$ & $2$ & $0$ & 2 & 2 & 2 & 1\\\hline  $\bf{L^k_{4|0}}(1)$ & $2$ & $1$  & 2 & 2 & 2 & 1\\\hline 
           $\bf{L^l_{4|0}}$ & $3$ & $0$  & 1 & 0 & 1 &1  \\\hline 
           $\bf{L^m_{4|0}}(\eps),\eps\neq1$ & $3$ & $0$ & 1 & 0 & 1  &1  \\\hline 
            $\bf{L^m_{4|0}}(1)$ &$2$ & $0$  & 2 & 2 & 2  &1  \\\hline 
           $\bf{L^n_{4|0}}$ & $3$ & $0$  & 1 & 0 & 1 & 1 \\\hline 
           
           $\bf{L^p_{4|0}}(\eps)$ & $2$ & $0$  & 2 & 2 & 2 & 1  \\\hline
           
           $\bf{L^q_{4|0}}$ & 1 &  2 & 3 &3  & 1 &  0 \\\hline
           
           $\bf{L^r_{4|0}}$ & 1 &  2 & 3 & 3 & 1 & 0  \\\hline
        
        $\bf{L^{q+r}_{4|0}}$ & 1 &  2 & 3 & 3 & 1 & 0  \\\hline
           
    \end{tabular}
   \\~\\ Invariants for Lie superalgebras of $\sdim=(4|0)$ obtained as $T^*$-extensions of the Lie superalgebra $\bf{L^1_{2|0}}$.
\end{center}}
\normalsize{}

\small{
\begin{center}
    \begin{tabular}{|c|c|c|c|c|c|c|}
    \hline
          $\fg$ & $\fg^{(1)}$ & $\dim(\fz(\fg))$ & $XH^1(\fg;\K)$& $XH^2(\fg;\K)$ & $XH^3(\fg;\K)$&$XH^4(\fg;\K)$\\\hline 
           $\bf{L^s_{4|0}}$ & 1 & 2 &3  &4 & 3 & 1\\\hline
           $\bf{L^t_{4|0}}$ & 1 & 2 &3  &4 & 3 & 1\\\hline
           $\bf{L^{s+t}_{4|0}}$ & 1 & 2 &3  &4 & 3 & 1\\\hline
           $\bf{L^u_{4|0}}$ & 1 & 2 &3  &3 & 1 & 0\\\hline
           $\bf{L^v_{4|0}}$ & 2 & 1 &2  &2 & 1 & 0\\\hline 
           $\bf{L^w_{4|0}}$ &1  & 2 &3  &4 & 3 & 1\\\hline
           $\bf{L^x_{4|0}}$ & 2 & 1 &2  &2 & 2 & 1\\\hline
           $\bf{L^y_{4|0}}$ &2  & 0 &  2&1 & 0 &0 \\\hline
           $\bf{L^z_{4|0}}$ & 2 & 0 &2  &1 & 0 & 0\\\hline
           $\bf{L^{aa}_{4|0}}$ & 2 & 0&2  &2 & 2 & 1\\\hline
           $\bf{L^{bb}_{4|0}}$ & 2 & 0 &2  &2 & 2 & 1\\\hline
           $\bf{L^{cc}_{4|0}}$ & 2 & 0 &2  &1 & 0 & 0\\\hline
           $\bf{L^{dd}_{4|0}}$ & 2 & 0 &2  &1 & 0 & 0\\\hline
    \end{tabular}
   \\~\\ Invariants for Lie superalgebras of $\sdim=(4|0)$ obtained as $T^*$-extensions of the Lie superalgebra $\bf{L^2_{2|0}}$.
\end{center}}


\end{document}